\documentclass{article}
\usepackage{amssymb,latexsym,amsmath,amsthm}

\newtheorem{theorem}{Theorem}[section]
\newtheorem{lemma}[theorem]{Lemma}
\newtheorem{corollary}[theorem]{Corollary}
\newtheorem{proposition}[theorem]{Proposition}

\theoremstyle{definition}
\newtheorem{definition}[theorem]{Definition}

\newtheorem{remark}[theorem]{Remark}

\newtheorem{notation}[theorem]{Notation}

\newcommand{\Coker}{{ \rm Coker }}
\newcommand{\Ext}{{ \rm Ext }}
\newcommand{\Hom}{{ \rm Hom }}
\newcommand{\Ker}{{ \rm Ker }\,}

\newcommand{\Ob}{{ \rm Ob }\,}

\newcommand{\Mod}{{ \rm Mod }}
\newcommand{\GMod}{{ \rm GMod }}

\newcommand{\FSMod}{{ \rm FSMod }}

\newcommand{\ad}{{ \rm ad }}
\newcommand{\tw}{{ \rm tw }}

\renewcommand{\frak}{\mathfrak}

\newcommand{\g}{\hbox{-}}
\newcommand{\uddots}{\mathinner{\mkern1mu\raise1pt\vbox{\kern7pt\hbox{.}}
\mkern2mu\raise4pt\hbox{.}\mkern2mu\raise7pt\hbox{.}\mkern1mu}}

\newcommand{\hueca}[1]{\mathbb{#1}}
\newcommand{\lddots}{
\mathinner{
\mkern1mu\raise1pt}\vbox{\kern7pt\hbox{.}}
\mkern2mu\raise3pt\hbox{.}
\mkern2mu\raise7pt\hbox{.}\mkern1mu}

\renewcommand{\Im}{{\rm Im}\,}

\newcommand{\rightdashmap}[1]{\smash{\mathop{\hbox to 
20pt{-\,-\,-\,\rightarrowfill}}\limits^{#1}}}

\newcommand{\rightmap}[1]{\smash{\mathop{\hbox to 
20pt{\rightarrowfill}}\limits^{#1}}}
\newcommand{\leftmap}[1]{\smash{\mathop{\hbox to 
20pt{\leftarrowfill}}\limits^{#1}}}

\newcommand{\lmapdown}[1]
{\llap{$\vcenter{\hbox{$\scriptstyle#1$}}$}
\Bigg\downarrow}
\newcommand{\rmapdown}[1]{\Bigg\downarrow\rlap{$\vcenter{\hbox{$\scriptstyle#1$}
}$}}

\newcommand{\shortlmapdown}[1]
{\llap{$\vcenter{\hbox{$\scriptstyle#1$}}$}\big\downarrow}
\newcommand{\shortrmapdown}[1]
{\big\downarrow\rlap{$\vcenter{\hbox{$\scriptstyle#1$}}$}}

\newcommand{\shortlmapup}[1]
{\llap{$\vcenter{\hbox{$\scriptstyle#1$}}$}\big\uparrow}
\newcommand{\shortrmapup}[1]
{\big\uparrow\rlap{$\vcenter{\hbox{$\scriptstyle#1$}}$}}

\newcommand{\longrightmap}[1]{\smash{\mathop{\hbox to 
4cm{\rightarrowfill}}\limits^{#1}}}
\newcommand{\longleftmap}[1]{\smash{\mathop{\hbox to 
4cm{\leftarrowfill}}\limits^{#1}}}
\newcommand{\medrightmap}[1]{\smash{\mathop{\hbox to 
2cm{\rightarrowfill}}\limits^{#1}}}
\newcommand{\medleftmap}[1]{\smash{\mathop{\hbox to 
2cm{\leftarrowfill}}\limits^{#1}}}

\newcommand{\idmapdown}[1]
{\hskip-8pt\mathop{\hskip-5pt\raise6pt
\hbox{$\scriptstyle#1$}\hskip-5pt\swarrow}}

\newcommand{\ddmapdown}[1]
{\hskip-5pt\mathop{\searrow\hskip-6pt\raise5pt\hbox{$\scriptstyle#1$}}}

\newcommand{\idmapup}[1]
{\hskip-5pt\mathop{\nwarrow\hskip-6pt \raise5pt\hbox{$\scriptstyle#1$}}}

\newcommand{\ddmapup}[1]
{\hskip-8pt\mathop{\hskip-5pt\raise6pt\hbox{$\scriptstyle#1$}\hskip-5pt\nearrow}
}

\newcommand{\flechypunt}[2]{\ \smash{\mathop{
   \raise 3pt \hbox to 40pt{\rightarrow}\hskip-40pt \lower 3pt
   \hbox to 40pt{\dashrightarrow}}\limits^{#1}_{#2}}\ }

\newcommand{\longequal}{\ \smash{\mathop{
   \raise 5pt \hbox to 35pt{\hrulefill}\hskip-35pt \lower 0pt
   \hbox to 35pt{\hrulefill}}}\ }

\newcommand{\raya}[1]{\ \smash{\mathop{\raise 2pt \hbox to 
10pt{\hrulefill}}\limits^{#1}}\ }

\title{\bf {\Large $(b,\nu)$-algebras and their twisted modules}}

\author{R. Bautista, E. P\'erez, L. Salmer\'on}

\begin{document}
 \date{}
 \maketitle
  
  \renewcommand{\thefootnote}{}


  
  \begin{abstract}
\noindent 
We give an intrinsic characterization of the closure under shifts $\widehat{\cal A}$ of a given strictly unital $A_\infty$-category ${\cal A}$.  We study some arithmetical properties of its higher operations and special conflations in the precategory of cocycles ${\cal Z}({\cal A})$ of its $A_\infty$-category of twisted modules.
We exhibit a structure for ${\cal Z}(\widehat{\cal A})$ similar to a special Frobenius category. We derive  that the cohomology category ${\cal H}(\widehat{\cal A})$ appears as the corresponding stable category and then we review how this implies that ${\cal H}(\widehat{\cal A})$ is a triangulated category.
\end{abstract}

  \section{Introduction} In this work we consider a special kind of algebraic structures  $\hat{Z}$, which we call $(b,\nu)$-algebras over an algebra with enough idempotents $\hat{S}$, arising from a special kind of  $A_\infty$-categories with strict identities and we give a detailed proof of the fact that the cohomology category ${\cal H}(\hat{Z})$ associated to the $A_\infty$-category of its
   twisted modules $\tw(\hat{Z})$ is a triangulated category. 
   With a different language, this last result is known, see \cite{Kontsevich},  \cite{K1}{(7.6)-(7.7)}, \cite{hsb}(7.4), \cite{LH}(7.2), and \cite{SF}(3.29).
   The notion of $(b,\nu)$-algebra corresponds to the closure under shifts $\widehat{\cal A}$ of a given $A_\infty$-category ${\cal A}$ with strict units introduced in \cite{K1}, where it is denoted by $\hueca{Z}{\cal A}$. It provides an intrinsic  formulation which permits to make a more detailed description of this $A_\infty$-category and to exhibit some  nice arithmetical features.
   
   Here, we give a  detailed description of the triangular structure in ${\cal H}(\hat{Z})$, which involves a more explicit study of the higher operations of $\hat{Z}$ related to the actions forming part of the structure of $\hat{Z}$, and some special sequences in the precategory ${\cal Z}(\hat{Z})$ of cocycles with respect to the first higher operation $\hat{b}_1^{tw}$ of the $b$-category $\tw(\hat{Z})$. 
    By a precategory we mean  an algebraic structure ${\cal C}$ which satisfies all the requirements of a category, except for the associativity of the composition. We call these sequences \emph{special conflations} and we show that they provide the precategory ${\cal Z}(\hat{Z})$ with a structure which is similar to an exact structure in an exact special Frobenius category, see \cite{BMJ} and \cite{dgb}(8.6), and that they similarly induce a triangulated structure in ${\cal H}(\hat{Z})$.
   
   We believe that the elementary approach and degree of detail with which we work within ${\cal Z}(\hat{Z})$ to study its exact structure induces a familiarity with the internal environment of this structure and the cohomology category ${\cal H}(\hat{Z})$, making their study available to a broader audience.

   Our motivation for this study was to have a deeper understanding of the following theorem of Keller and Lef\`evre-Hasegawa, see  {\cite{LH}}\S7, which plays an essential role in our  argumentation in \cite{hsb},  which follows closely  that of \cite{KKO}. There, $\Delta$ denotes the direct sum of a finite set $\{\Delta_1,\ldots,\Delta_n\}$ of non-isomorphic indecomposable $\Lambda$-modules, where $\Lambda$ is a finite-dimensional $k$-algebra with unit.
    
  \begin{theorem}\label{T: Keller-Lefevre} If the Yoneda   $A_\infty$-algebra $A$ associated to the $\Lambda$-module
  $\Delta$ is strictly unital, then there is an equivalence of categories
  ${\cal F}(\Delta){\simeq} H^0(\tw(A))$.
  \end{theorem}
  
  We want an  explicit description of this equivalence, which permits us to keep track of the image of the exact sequences of  ${\cal F}(\Delta)$ in $H^0(\tw(A))$. Space problems do not allow us to include such description in this paper, but we will come back to this aim in a forthcoming paper, where we will  also study the relation between ${\cal H}(\hat{Z})$ and the triangulated category of twisted modules of the graded bocs associated to a finite section of $\hat{Z}$, see \cite{dgb}.

  \section{$b$-algebras}\label{S: Differential graded S-bocses and b-algebras}
    
 Throughout this article, we assume that $k$ is a fixed ground field. 
 Unless we specify it otherwise, the terms category and functor mean $k$-category and $k$-functor respectively.    
 
 We  first recall some basic well known notions. We consider $k$-algebras (possibly without unit) but with enough idempotents, as in \cite{S} and \cite{W}, in the following  sense. Moreover, we consider only unitary modules and bimodules over these type of algebras.

  \begin{definition}\label{D: anillos con suficientes idemptes}
  An \emph{(associative) $k$-algebra $A$}  is a vector space over the ground 
field $k$,
  endowed with a product (a binary operation) such that:
  \begin{enumerate}
   \item $(ab)c=(ab)c, (a+b)c=ac+bc, \hbox{ and } a(b+c)=ab+ac, \hbox{ for } 
a,b,c\in A$
   \item $a(\lambda b)=\lambda (ab)=(\lambda a)b$, for all $a,b\in A$ and $\lambda\in k$.
  \end{enumerate}
  A \emph{morphism of $k$-algebras} $\phi:A\rightmap{}A'$ is a $k$-linear map which preserves the product. 
  
  The $k$-algebra $A$ \emph{has enough idempotents} iff it is equipped with  a family 
  $\{e_i\}_{i\in {\cal P}}$ of pairwise primitive orthogonal  idempotents of $A$ such that 
  $$\bigoplus_{i\in {\cal P}}e_iA=A=\bigoplus_{i\in {\cal P}}Ae_i.$$
  Following \cite{S}, we call $\{e_i\}_{i\in {\cal P}}$ \emph{the distinguished family of idempotents of $A$.}

  The category $\Mod\g A$ is the category of \emph{unitary right $A$-modules}, 
that is the right $A$-modules $M$ such that $M=\bigoplus_{i\in {\cal P}}Me_i$. 
  The category of left $A$-modules $A\g \Mod$ is defined similarly. 
  \end{definition}
  
  \begin{remark}\label{R: centralidad del campo k} Let $A$ be a $k$-algebra with enough idempotents $\{e_i\}_{i\in {\cal P}}$. Then, each linear map $k\rightmap{}e_iAe_i$ determined by $1\mapsto e_i$ is a morphism of rings with unit, so if $M$ is a right $A$-module, each $Me_i$ is a right $e_iAe_i$-module, and it inherits a natural structure of right $k$-module. 
  
  A right unitary $A$-module admits, by definition, an abelian group decomposition $M=\bigoplus_{i\in {\cal P}}Me_i$. So we can equip $M$ naturally with the vector space structure given by the vector space structure of the $Me_i$'s, so we get a vector space decomposition $M=\bigoplus_{i\in {\cal P}}Me_i$. We proceed similarly with the left unitary $A$-modules. 
  
  So, any unitary $A$-$A$-bimodule $M=\bigoplus_{i,j\in {\cal P}}e_jMe_i$ is naturally a $k$-$k$-bimodule. Throughout this paper, our definition of unitary $A$-$A$-bimodule includes tacitely that the action of $k$, on every unitary $A$-$A$-bimodule $M$, is central: that is such that $\lambda m=m\lambda$, for all $m\in M$ and $\lambda\in k$.  
   
  Notice that the $k$-algebra $A$ itself is a unitary $A$-$A$-bimodule, where the left $k$-module structure of the unitary $A$-$A$-bimodule $A$ coincides with the original $k$-vector space structure of $A$ and the right $k$-module structure of the unitary $A$-$A$-bimodule $A$ coincides with the one defined by $a\lambda:=\lambda a$, for $\lambda\in k$ and $a\in A$, where $(\lambda,a)\mapsto \lambda a$ is the action of the original $k$-vector space underlying the $k$-algebra  $A$. So the action of $k$ on $A$ is indeed central. 
  \end{remark}

  \begin{definition} A \emph{graded $k$-algebra $A$ with enough idempotents} is 
  a graded $k$-algebra  $A=\bigoplus_{q\in \hueca{Z}}A_q$ equipped with a distinguished family of orthogonal idempotents $\{e_i\}_{i\in {\cal P}}$, which are all homogeneous of degree 0. 

   For such graded $k$-algebra with enough idempotents, the category $\GMod\g A$ is \emph{the category of graded (unitary)  right $A$-modules}, 
that is the graded right $A$-modules $M$ such that $M=\bigoplus_{i\in {\cal P}}Me_i$. The morphism spaces of $\GMod\g A$ are defined by 
 $$\Hom_{\GMod\g A}(M,N):=\bigoplus_{d\in \hueca{Z}}\Hom^d_{\GMod\g A}(M,N),$$
where $\Hom^d_{\GMod\g A}(M,N)$ denotes the vector space of homogeneous morphisms $f:M\rightmap{}N$ of graded right $A$-modules of degree $d$.

  The category of graded left  $A$-modules $A\g\GMod$ and the category of graded $A$-$A$-bimodules are defined similarly.

  Given a graded vector space $M$, the degree of any homogeneous element $m\in M$ will be denoted by $\vert m\vert$. Given homogeneous morphisms of graded  $A\g A$-bimodules $f:M\rightmap{}M'$ and $g:N\rightmap{}N'$, we consider the associated tensor product morphism $f\otimes g:M\otimes_A N\rightmap{}M'\otimes_AN'$ which is defined, following the Koszul sign convention, by the  formula 
  $(f\otimes g)(m\otimes n)=(-1)^{\vert g\vert\vert m\vert}f(m)\otimes g(n)$, for any homogeneous elements $m\in M$ and $n\in N$.
  \end{definition}

\begin{definition}\label{D:  morphism of algebras}
Let $A$ and $A'$ be graded $k$-algebras, with enough idempotents
$\{e_i\}_{i\in {\cal P}}$ and $\{e'_j\}_{j\in {\cal P}'}$, respectively. 
Then, a \emph{morphism of graded algebras
(with enough idempotents as above)} is a 
$k$-linear homogeneous map $\psi:A\rightmap{}A'$ such that $\psi$ preserves the product  
 and 
  $\psi(\{e_i \mid i\in {\cal P}\})\subseteq\{e'_j\mid j\in {\cal P}'\}$. 
\end{definition}

\begin{definition}\label{D: elementary k-algebra}
A $k$-algebra $S$ with a distinguished family of orthogonal idempotents  is called 
\emph{elementary} if $e_iSe_i=ke_i$, 
for all $i\in {\cal P}$, and $e_jSe_i=0$ for all
$i,j\in {\cal P}$, with $i\not= j$. We will consider such $k$-algebras as graded $k$-algebras concentrated in degree 0. 
\end{definition}

 We need to recall some notions and statements from \cite{hsb}, se also \cite{K1}, \cite{LH}, and \cite{KKO}, but we adapt those to the context of graded $k$-algebras with enough idempotents. We prefered the language of graded $k$-algebras with enough idempotents to the equivalent one of small  graded $k$-categories because the notation is simpler and many statements (and their proofs) on graded $k$-algebras with unit have similar formulations.

  We are interested in the following type of structures. 
   
   \begin{definition}\label{D: b-category} Let $S$ be an elementary 
   $k$-algebra with enough idempotents $\{e_i\}_{i\in {\cal P}}$. 
 A \emph{$b$-algebra $Z$ over $S$}   is a graded unitary $S$-$S$-bimodule  $Z=\bigoplus_{q\in \hueca{Z}}{Z_q}$ equipped  with a family $\{b_n:Z^{\otimes n}\rightmap{}Z\}_{n\in \hueca{N}}$  of homogeneous morphisms of $S$-$S$-bimodules of degree $\vert b_n\vert=1$, for all $n\in \hueca{N}$. 
  It is required that, for each $n\in \hueca{N}$, the maps of the family  satisfy the following relation 
   $$S_n^b:\sum_{\scriptsize\begin{matrix}r+s+t=n\\ s\geq 1;r,t\geq 0\end{matrix}} 
   b_{r+1+t}(id^{\otimes r}\otimes b_s\otimes id^{\otimes t})=0.$$
 \end{definition}
 
 Thus, a $b$-algebra is simply the bar construction of an $A_\infty$-algebra. We consider also $b$-categories as defined in \cite{hsb}(6.8), which again are simply the bar construction of  $A_\infty$-categories, see \cite{hsb}\S6.

We keep the notations introduced before, so $S=(S,\{e_i\}_{i\in {\cal P}})$ is an elementary $k$-algebra with enough idempotents and $(Z,\{b_n\}_{n\in \hueca{N}})$ is a $b$-algebra over $S$. We will denote by $\FSMod\g S$ the category of right (unitary) $S$-modules \emph{with finite support}, that is the unitary right $S$-modules $X$ such that $Xe_i=0$, for almost all $i\in {\cal P}$.

\begin{definition}\label{R: recordatorio de ad(Z)} Let $Z=(Z,\{b_n\}_{n\in \hueca{N}})$ be a $b$-algebra over $S$.
  As in \cite{hsb}(6.5), a
   $b$-category $\ad(Z)$ is defined by the following. The objects of $\ad(Z)$ are the right (unitary) $S$-modules with finite support,  the spaces of morphisms are given by 
  $$\ad(Z)(X,Y):=\bigoplus_{i,j\in {\cal P}}\Hom_k(Xe_i,Ye_j)\otimes_ke_jZe_i,$$
  with the canonical grading of the tensor product where $\Hom_k(Xe_i,Xe_j)$ is considered as a graded vector space concentrated in degree 0. The morphisms $b^{ad}_n$ are defined, for $n\in \hueca{N}$ and a sequence of right $S$-modules $X_0,X_1,\ldots,X_n$, on typical generators by
  $$\begin{matrix}\ad(Z)(X_{n-1},X_n)\otimes_k\cdots \otimes_k\ad(Z)(X_1,X_2)\otimes_k\ad(Z)(X_0,X_1)\hbox{\,}\rightmap{b^{ad}_n}\hbox{\,}\ad(Z)(X_0,X_n)\end{matrix}$$ 
  $$\begin{matrix}(f_n\otimes a_n)\otimes\cdots\otimes(f_2\otimes a_2)\otimes (f_1\otimes a_1)&\longmapsto&
  f_n\cdots f_2f_1\otimes b_n(a_n\otimes\cdots\otimes a_1).\end{matrix}$$
  
   In the preceding recipe, since for given support-finite unitary right $S$-modules $X$ and $Y$, we have  the unitary $S$-$S$-bimodule   
   $$\Hom_k(X,Y)=\bigoplus_{i,j\in {\cal P}}e_i\Hom_k(X,Y)e_j=\bigoplus_{i,j\in{\cal P}}\Hom_k(Xe_i,Ye_j),$$
  we identify the elements of $\Hom_k(Xe_i,Ye_j)$ with the corresponding elements in $\Hom_k(X,Y)$, so the composition $f_n\cdots f_2f_1$ makes sense. 
\end{definition}

 \begin{remark}\label{N: la base de B} 
  A non-zero element $a$ in an $S$-$S$-bimodule $Z$  is called
  \emph{directed}
  iff $a=e_jae_i$, for some $i,j\in {\cal P}$. In this case, we will write $v(a):=j$ and $u(a):=i$. 
  A subset $L$ of $Z$  is called \emph{directed}  iff each one of its elements is so. 
 
It is convenient to fix a \emph{directed basis} $\hueca{B}$ for the graded $S$-$S$-bimodule $Z=\bigoplus_{q\in \hueca{Z}}Z_q$. It is chosen as follows. For each $q\in \hueca{Z}$ and  $i,j\in {\cal P}$, we choose a $k$-basis $\hueca{B}_q(i,j)$ for the vector space $e_jZ_qe_i$; then, we consider the basis  $\hueca{B}_q=\bigcup_{i,j}\hueca{B}_q(i,j)$ of $Z_q$. 
Finally,  we can consider the $k$-basis $\hueca{B}=\bigcup_{q\in \hueca{Z} 
}\hueca{B}_q$ of $Z=\bigoplus_{q\in \hueca{Z}}Z_q$.  

The elements $f\in \ad(Z)(X,Y)$ are called  the morphisms of $\ad(Z)$ and we often say that $f:X\rightmap{}Y$ is a morphism in $\ad(Z)$ to  make explicit its domain and codomain. 
Any morphism $f\in \ad(Z)(X,Y)$, can be written uniquely as a sum $\sum_{a\in \hueca{B}}f_a\otimes a$. In the following, when we consider a morphism written as $f=\sum_af_a\otimes a$, we mean this description.
Moreover, such an $f$ is homogeneous of degree $d$ iff $f_a=0$ for all $a\in\hueca{B}$ with $\vert a\vert\not=d$.  
\end{remark}

\begin{definition}\label{D: strict morphisms mejor} A directed element $a$ of a $b$-algebra $Z$ is called \emph{strict} iff for any $n\not= 2$ and any sequence of  directed elements  
$a_1\in e_{u_1}Ze_{u_0},\ldots,a_n\in e_{u_n}Ze_{u_{n-1}}$, such that 
$a\in \{a_1,\ldots,a_n\}$, we have 
$b_n(a_n\otimes\cdots\otimes a_1)=0$. 

A morphism $f:X\rightmap{}Y$ of $\ad(Z)$ is called \emph{strict} iff it has the form $f=\sum_af_a\otimes a$, where each $a$, with $f_a\not=0$, is a strict element of $Z$. 
\end{definition}

\begin{definition}\label{D: unitary strict} We say that the $b$-algebra $Z=(Z,\{b_n\}_{n\in \hueca{N}})$, over the elementary algebra $S$ with distinguished idempotents $\{e_i\}_{i\in {\cal P}}$,  is \emph{unitary strict} iff for each $i\in {\cal P}$ there is a homogeneous element $\frak{e}_i\in Z$ with degree $\vert \frak{e}_i\vert=-1$ satisfying the following:
\begin{enumerate}
 \item $\frak{e}_i=e_i\frak{e}_ie_i$, for all $i\in {\cal P}$;
 \item $\frak{e}_i$ is a strict element of $Z$, for all $i\in {\cal P}$;  
 \item For each homogeneous element $a\in Z$, we have 
 $$b_2(\frak{e}_i\otimes a)=e_ia \hbox{ \ \ and \ \ } b_2(a\otimes \frak{e}_i)=(-1)^{\vert a\vert+1}ae_i.$$
\end{enumerate}
In this case, the elements of the family $\{\frak{e}_i\}_{i\in {\cal P}}$ are called \emph{the strict units of $Z$}. 
\end{definition}

When we are dealing with a unitary strict $b$-algebra $Z$, we always assume that the directed basis $\hueca{B}$ fixed in (\ref{N: la base de B}) contains the strict units of $Z$. 

\begin{notation}\label{D: hueca(I) y circ para ad(hat(B))} Assume that $Z=(Z,\{b_n\}_{n\in \hueca{N}})$ is a unitary strict $b$-algebra. Given $a_1,a_2\in Z$, we often write  $a_1\circ a_2:=b_2(a_1\otimes a_2)$. We have to be careful because here, for $a_1$ and $a_2$ homogeneous, we have $\vert a_1\circ a_2\vert=\vert a_1\vert+\vert a_2\vert+1$. 
With this notation, we have $\frak{e}_i\circ \frak{e}_i=\frak{e}_i$, $\frak{e}_i\circ a=e_ia$, and $a\circ \frak{e}_i=(-1)^{\vert a\vert+1}ae_i$, for all $i\in {\cal P}$ and all homogeneous $a\in Z$. 

Likewise, given morphisms $f\in \ad(Z)(X,Y)$ and $g \in \ad(Z)(Y,W)$, we will 
write  
$g\circ f:=b_2^{ad}(g\otimes f)\in \ad(Z)(X,W).$
Again, for $f$ and $g$ homogeneous, we have $\vert g\circ f\vert=\vert g\vert +\vert f\vert+1$.
\medskip

For each object $X$ of $\ad(Z)$, set
$\hueca{I}_X:=\sum_{u\in {\cal P}}id_{Xe_u}\otimes \frak{e}_u\in \ad(Z)(X,X)$. 
\end{notation}

\begin{lemma}\label{L: hueca(I)X} 
In the context of the last definition, 
 we see that the morphisms $\hueca{I}_X$ are strict morphisms in $\ad(Z)$. Moreover, for any homogeneous morphism $f:X\rightmap{}Y$ of $\ad(Z)$, we have 
 $\hueca{I}_Y\circ f=f$ and $f\circ \hueca{I}_X=(-1)^{\vert f\vert+1}f$.
\end{lemma}

\begin{proof} Let $f=\sum_af_a\otimes a$ be a homogeneous morphism of $\ad(Z)$, so $\vert a\vert=\vert f\vert$ for all index $a$. Then, from the properties of the strict units, we have 
$
  f\circ \hueca{I}_X=
  \sum_af_a\otimes b_2(a\otimes \frak{e}_{u(a)})
  =
  (-1)^{\vert f\vert +1}\sum_af_a\otimes a=(-1)^{\vert f\vert +1}f
  $
and, also,  $\hueca{I}_Y\circ f =
   \sum_af_a\otimes b_2(\frak{e}_{v(a)}\otimes a)=\sum_af_a\otimes a=f$. 
\end{proof}

\begin{remark}\label{R: bn ad strict}
Let  $X_0\rightmap{f_1}X_1,...,X_{n-1}\rightmap{f_n}X_n$ be a sequence of morphisms in $\ad(Z)$ with  $n\not=2$. Then, if 
 at least one of the morphisms $f_1,\ldots,f_n$ is strict, we have  
$b_n^{ad}(f_n\otimes\cdots\otimes f_1)=0.$ 
\end{remark}

\begin{remark} Until the end of this section, we assume that  $Z=(Z,\{b_n\}_{n\in \hueca{N}})$ is a unitary strict $b$-algebra with strict units $\{\frak{e}_u\}_{u\in {\cal P}}$. Then, there is an ``embedding functor'' $L:\FSMod\g S\rightmap{}\ad(Z)$ such that $L(f)=\sum_{u\in {\cal P}}f_u\otimes \frak{e}_u$, where $f_u:Xe_u\rightmap{}Ye_u$ denotes the restriction of the morphism $f:X\rightmap{}Y$. Here, the sum is finite because $X$ has finite support. The preceding phrase means that $L$ is a function on objects and on morphisms, which is the identity on objects and  maps  morphisms $f:X\rightmap{}Y$ onto  morphisms $L(f):X\rightmap{}Y$ in $\ad(Z)$ in such a way that $L(id_X)=\hueca{I}_X$ and it maps each  composition $gf$  of a pair of composable morphisms in $\FSMod\g S$ onto  $L(gf)=L(g)\circ L(f)$.

The morphisms in $\ad(Z)$ of the form $f=\sum_uf_u\otimes \frak{e}_u$, that is those in the image of $L$, play an important role in this work. We will call them \emph{special morphisms}. The image of $L$ is a category isomorphic to $\FSMod\g S$. 
\end{remark}

\begin{lemma}\label{L: morf especiales vs operaciones de ad(Z)}
Let $f_1:X_0\rightmap{}X_1,\ldots,f_n:X_{n-1}\rightmap{}X_n$ be homogeneous morphisms in $\ad(Z)$. Then, the following holds.
\begin{enumerate}
 \item For any special morphism  $g:X_n\rightmap{}U$, we have
 $$b_n^{ad}(g\circ f_n\otimes f_{n-1}\otimes \cdots\otimes f_1)=g\circ b_n^{ad}(f_n\otimes f_{n-1}\otimes \cdots \otimes f_1).$$
 \item For any special morphism $h:V\rightmap{}X_0$, we have 
 $$b_n^{ad}(f_n\otimes f_{n-1}\otimes \cdots\otimes f_1\circ h)=(-1)^{\vert f_n\vert+\cdots+\vert f_2\vert+1} b_n^{ad}(f_n\otimes f_{n-1}\otimes \cdots \otimes f_1)\circ h.$$
 \item If $n\geq 2$, $i\in [2,n]$, and $f_i=f'_i\circ h_i$, where $f'_i:U_i\rightmap{}X_i$ is homogeneus and $h_i:X_{i-1}\rightmap{}U_i$ is a special morphism, we have
 $$b_n^{ad}(f_n\otimes\cdots f'_i\circ h_i\otimes f_{i-1}\otimes\cdots\otimes f_1)=(-1)^{\vert f'_i\vert +1} b_n^{ad}(f_n\otimes\cdots\otimes f'_i\otimes h_i\circ f_{i-1}\cdots \otimes f_1)$$
\end{enumerate}
\end{lemma}

\begin{proof} We only prove (2), since the other verifications are similar. We write $h=\sum_uh_u\otimes \frak{e}_u$ and $f_i=\sum_{a_i}(f_i)_{a_i}\otimes a_i$, for all $i\in [1,n]$. The left term of the equation in (2) is 
$\sum_{u,a_1,\ldots,a_n}(f_n)_{a_n}\cdots (f_1)_{a_1}h_u\otimes b_n(a_n\otimes\cdots\otimes a_1\circ \frak{e}_u)$
  while the right one is
 $$(-1)^{\vert f_n\vert +\cdots+\vert f_2\vert+ 1}\sum_{u,a_1,\ldots,a_n}(f_n)_{a_n}\cdots (f_1)_{a_1}h_u\otimes b_n(a_n\otimes\cdots\otimes a_1)\circ \frak{e}_u.$$
 But 
 $
 b_n(a_n\otimes\cdots\otimes a_1\circ \frak{e}_u)
 =
 (-1)^{\vert f_1\vert +1}b_n(a_n\otimes\cdots\otimes a_1e_u)$ 
 and 
 $b_n(a_n\otimes\cdots\otimes a_1)\circ \frak{e}_u=
 (-1)^{\vert f_n\vert+\cdots+\vert f_1\vert}b_n(a_n\otimes\cdots\otimes a_1)e_u$, so (2) follows. 
\end{proof}

\begin{corollary}\label{L: compos entre tres morfs homogeneos en ad(Z)}
Let $f:X\rightmap{}Y$, $g:Y\rightmap{}U$, and $h:U\rightmap{}V$ be homogeneous morphisms in $\ad(Z)$. Then, we have:
\begin{enumerate}
 \item If $f$ or $g$ are  special, then $(h\circ g)\circ f=(-1)^{\vert h\vert+1}h\circ (g\circ f)$.
 \item If $h$ is special, then $h\circ (g\circ f)=(h\circ g)\circ f$.  
\end{enumerate}
\end{corollary}

\begin{remark}\label{R: componentes de morfismos homogeneos en ad(Z)} Given $X\in \FSMod\g S$, suppose that we have a direct sum decomposition  $X=\bigoplus_{i=1}^n X_i$ of modules. Then, we have the projections $\pi_{X_i}:X\rightmap{}X_i$ and the injections 
 $\sigma_{X_i}:X_i\rightmap{}X$ associated to this decomposition. If we write $p_{X_i}=L(\pi_{X_i})$ and $s_{X_i}=L(\sigma_{X_i})$, we get the standard relations  $p_{X_i}\circ s_{X_i}=\hueca{I}_{X_i}$, for all $i$, $p_{X_i}\circ s_{X_j}=0$, for all $i\not=j$, and
 $\hueca{I}_X=\sum_{i=1}^ns_{X_i}\circ p_{X_i}$.
 
 Given a homogeneous morphism $f:X=\bigoplus_{i=1}^nX_i\rightmap{}\bigoplus_{j=1}^mY_j=Y$ in $\ad(Z)$, we  define  the $(j,i)$-component of $f$ by   
 $$f_{j,i}:=(-1)^{\vert f\vert+1}p_{Y_j}\circ f\circ s_{X_i}, \hbox{ for all } i,j.$$ 
  Using (\ref{L: compos entre tres morfs homogeneos en ad(Z)}) and the preceding standard relations, we can recover the morphism $f$ from its matrix $M(f):=(f_{j,i})$ using the formula $$f=(-1)^{\vert f\vert+1}\sum_{i,j}s_{Y_j}\circ f_{j,i}\circ p_{X_i}.$$
 The sign in the definition of $f_{j,i}$ is convenient because of the following.

 If $g:\bigoplus_{j=1}^mY_j\rightmap{}\bigoplus_{t=1}^sW_t$ is another homogeneous morphism in $\ad(Z)$, we can verify, using again  (\ref{L: compos entre tres morfs homogeneos en ad(Z)}) and the preceding standard relations, that the component  $(g\circ f)_{t,i}$ of the composition $g\circ f$ coincides with the $(t,i)$-entry $\sum_jg_{t,j}\circ f_{j,i}$ of the corresponding matrix product $M(g)\circ M(f)=(g_{t,j})\circ(f_{j,i})$.
 
 Moreover, observe that if $f=\sum_af_a\otimes a$, then $f_{j,i}=\sum_a\pi_{Y_j}f_a\sigma_{X_i}\otimes a$. For $a\in \hueca{B}$, the linear map
$f_a:\bigoplus_{i=1}^nX_i\rightmap{}\bigoplus_{j=1}^m Y_j$ has a matrix of linear maps $[f_a]:=[(f_a)_{j,i}]$, where $(f_a)_{j,i}=\pi_{Y_j}f_a\sigma_{X_i}$. The preceding expression for $f_{j,i}$  implies that
$[f_a]=((f_{j,i})_a)$.

Finally, notice that if $n=1=m$, then $f_{1,1}=f$, so we write, as usual, $f$ instead of $M(f)$. In the following sections, for simplicity, when we say that certain morphism $f$ in $\ad(Z)$ has matrix form $f=(f_{j,i})$, we mean that $M(f)=(f_{j,i})$. Then, we work with these matrices using the matrix product formula mentioned before and with the more general formula given in the next remark.
\end{remark}

\begin{remark}[On finite direct sums in $\ad(Z)$]\label{R: descomposicon de morfismos en ad(hat(B))}
 Assume that we have $n$ composable morphisms $f_1:X_0\rightmap{}X_1$, $f_2:X_1\rightmap{}X_2,\ldots,f_n:X_{n-1}\rightmap{}X_n$  of $\ad(Z)$, with  
$X_s=\bigoplus_{i_s\in I_s}X_{s,i_s}$, for $s\in [0,n]$. 
We want to describe the matrix of the morphism $b_n^{ad}(f_n\otimes\cdots\otimes f_1)$ of $\ad(Z)$. 
Applying (\ref{L: morf especiales vs operaciones de ad(Z)})(1)\&(2), for $i\in I_0$ and $j\in I_n$, we get 
$$\begin{matrix}
b_n^{ad}(f_n\otimes\cdots\otimes f_1)_{j,i}
&=&
(-1)^{\vert f_n\vert+\cdots+\vert f_1\vert}p_{X_{n,j}}\circ b_n^{ad}(f_n\otimes\cdots\otimes f_1)\circ s_{X_{0,i}}\hfill\\ 
&=&
(-1)^{\vert f_n\vert+\cdots+\vert f_1\vert}b_n^{ad}(p_{X_{n,j}}\circ f_n\otimes\cdots\otimes f_1)\circ s_{X_{0,i}}\hfill\\
&=&
(-1)^{\vert f_1\vert+1}b_n^{ad}(p_{X_{n,j}}\circ f_n\otimes\cdots\otimes f_1\circ s_{X_{0,i}}).\hfill\\
\end{matrix}$$   
Applying (\ref{L: morf especiales vs operaciones de ad(Z)})(3),  we see that this coincides with the following three expessions
$$ \sum_{r_{n-1}}(-1)^{\vert f_1\vert+1}b_n^{ad}(p_{X_{n,j}}\circ f_n\otimes s_{X_{n-1,r_{n-1}}}\circ (p_{X_{n-1,r_{n-1}}}\circ f_{n-1})\otimes\cdots\otimes f_1\circ s_{X_{0,i}}), $$
$$\sum_{r_{n-1}}(-1)^{\vert f_1\vert+\vert f_n\vert}b_n^{ad}(p_{X_{n,j}}\circ f_n\circ  s_{X_{n-1,r_{n-1}}}\otimes (p_{X_{n-1,r_{n-1}}}\circ f_{n-1})\otimes\cdots\otimes f_1\circ s_{X_{0,i}}),$$
and 
$\sum_{r_{n-1}}(-1)^{\vert f_1\vert+1}b_n^{ad}( (f_n)_{j,r_{n-1}}\otimes (p_{X_{n-1,r_{n-1}}}\circ f_{n-1})\otimes\cdots\otimes f_1\circ s_{X_{0,i}}).$
Then, applying the last argument repeatedly, we finally get 
$$b_n^{ad}(f_n\otimes\cdots\otimes f_1)_{j,i}=\sum_{r_1,r_2,\ldots,r_{n-1}}b_n^{ad}( (f_n)_{j,r_{n-1}}\otimes   (f_{n-1})_{r_{n-1},r_{n-2}}\otimes \cdots\otimes 
(f_1)_{r_1,i}).$$
 We define 
  $M(f_n)\otimes \cdots \otimes M(f_1)$ as the $I_n\times I_0$ matrix with $(j,i)$-entry 
  $$(M(f_n)\otimes \cdots \otimes M(f_1))_{j,i}:=\sum_{r_1,\ldots,r_{n-1}} (f_n)_{j,r_{n-1}}\otimes   (f_{n-1})_{r_{n-1},r_{n-2}}\otimes\cdots\otimes  
(f_1)_{r_1,i}$$
 With this notation, the preceding calculations are summarized in the formula
  $$M(b_n^{ad}(f_n\otimes\cdots\otimes f_1))=b_n^{ad}(M(f_n)\otimes\cdots\otimes M(f_1)).$$
 In particular, given the morphisms 
$X=\bigoplus_{i\in I}X_i\rightmap{f}Y=\bigoplus_{j\in J}Y_j\rightmap{g}W=\bigoplus_{t\in T}W_t$
in $\ad(Z)$,  for each $i\in J$ and $t\in T$, 
we have 
$$
 (g\circ f)_{t,i}=[b^{ad}_2(g\otimes f)]_{t,i}=\sum_jb^{ad}_2(g_{t,j}\otimes f_{j,i})=\sum_jg_{t,j}\circ f_{j,i}.$$
\end{remark}

  \begin{definition}\label{R: record de tw(Z)}
 As in \cite{hsb}(6.1), we can consider the $b$-category $\tw(Z)$ described by the following. The objects of $\tw(Z)$ are  the pairs $\underline{X}=(X,\delta_X)$ where $X$ is a right $S$-module with finite support and 
 $\delta_X\in \ad(Z)(X,X)_0$. Moreover:  
 \begin{enumerate}
  \item There is a finite filtration   
  $0=X_0\subseteq X_1\subseteq\cdots\subseteq X_{\ell(X)}=X$ of right $S$-modules
  such that if we express $\delta_X=\sum_{x\in \hueca{B}}f_x\otimes x$, where the maps $f_x\in \Hom_k(X,X)$ are uniquely determined, we have $f_x(X_r)\subseteq X_{r-1}$, for all $r\in [1,\ell(X)]$.
  \item We have $\sum_{s\geq 1} b_s^{ad}((\delta_X)^{\otimes s})=0$, where we notice that the preceding condition $\it{1}$ implies that $b_s^{ad}((\delta_X)^{\otimes s})=0$ for $s\geq \ell(X)$, so we are dealing with a finite sum.
 \end{enumerate}
 Given  $\underline{X},\underline{Y}\in\Ob(\tw(Z))$, we have the hom 
 graded $k$-vector space 
 $$\tw(Z)(\underline{X},\underline{Y})=\ad(Z)(X,Y)=\bigoplus_{i,j\in {\cal P}}\Hom_k(Xe_i,Ye_j)\otimes_k e_jZe_i.$$
 If $n\geq 1$ and $\underline{X}_0,\underline{X}_1,\ldots,\underline{X}_n\in \Ob(\tw(Z))$,  
 we have the following homogeneous linear map of degree 1  
 $$\begin{matrix}\tw(Z)(\underline{X}_{n-1},\underline{X}_n)\otimes_k\cdots \otimes_k\tw(Z)(\underline{X}_1,\underline{X}_2)\otimes_k\tw(Z)(\underline{X}_0,\underline{X}_1)\hbox{\,}\rightmap{b^{tw}_n}\hbox{\,}\tw(Z)(\underline{X}_0,\underline{X}_n)\end{matrix}$$ 
 which maps each homogeneous generator  $t_n\otimes\cdots\otimes t_2\otimes t_1$ on 
  $$\sum_{\scriptsize\begin{matrix}
  i_0,\ldots,i_n\geq 0\end{matrix}} 
   b^{ad}_{i_0+\cdots+i_n+n}(\delta_{X_n}^{\otimes i_n}\otimes t_n\otimes \delta_{X_{n-1}}^{\otimes i_{n-1}}\otimes t_{n-1}\otimes \cdots\otimes \delta_{X_1}^{\otimes i_1}\otimes t_1\otimes\delta_{X_0}^{\otimes i_0}),$$
   which is a finite sum. 
\end{definition}

\begin{remark}\label{R: notacion star}
Given  $f:(X,\delta_X)\rightmap{}(Y,\delta_Y)$ and $g:(Y,\delta_Y)\rightmap{}(W,\delta_W)$, two morphisms in $\tw(Z)$, we will use the notation: 
$g\star f=b_2^{tw}(g\otimes f)$. 

For each $\underline{X}=(X,\delta_X)$ and $\underline{Y}=(Y,\delta_Y)\in \tw(Z)$, we have the complex of vector spaces $\tw(Z)(\underline{X},\underline{Y})$ with differential  $b_1^{tw}$, so we can consider the graded vector spaces 
$${\cal K}(Z)(\underline{X},\underline{Y}):=\Ker b_1^{tw}\leq \tw(Z)(\underline{X},\underline{Y}) \hbox{ \ ; \ }
{\cal I}(Z)(\underline{X},\underline{Y}):=\Im b_1^{tw}\leq {\cal K}(Z)(\underline{X},\underline{Y}).$$

Then, we have the following.
\begin{enumerate}
 \item ${\cal K}(Z)$ is closed under the product $\star$. That is,   
 if we have $f\in {\cal K}(Z)(\underline{X},\underline{Y})$ and $g\in {\cal K}(Z)(\underline{Y},\underline{W})$, then $g\star f\in {\cal K}(Z)(\underline{X},\underline{W})$. 
 
 Indeed, if $b_1^{tw}(f)=0$ and $b_1^{tw}(g)=0$, with $g$ homogeneous, since $\tw(Z)$ is a $b$-category, we have 
  $$\begin{matrix}
    0&=& [b_1^{tw}b_2^{tw}+b^{tw}_2(id\otimes b_1^{tw})+b^{tw}_2(b_1^{tw}\otimes id)](g\otimes f)\hfill\\
    &=&
    b_1^{tw}(g\star f)+(-1)^{\vert g\vert}b^{tw}_2(g\otimes b_1^{tw}(f))+b^{tw}_2(b_1^{tw}(g)\otimes f)\hfill\\
    &=&
    b_1^{tw}(g\star f).\hfill\\
    \end{matrix}$$
 \item ${\cal I}(Z)$ is an ideal of ${\cal K}(Z)$. That is, 
 if we have $f\in {\cal I}(Z)(\underline{X},\underline{Y})$ and 
 $g\in {\cal K}(Z)(\underline{Y},\underline{W})$ (or $g\in {\cal K}(Z)(\underline{W},\underline{X})$), then $g\star f\in {\cal I}(Z)(\underline{X},\underline{W})$ (resp. $f\star g\in {\cal I}(Z)(\underline{W},\underline{Y})$). 
 
  Indeed, if  we have $h\in \tw(Z)(\underline{X},\underline{Y})$ such that $f=b_1^{tw}(h)$ and $g\in {\cal K}(Z)(\underline{Y},\underline{W})$ homogeneous, then, as before, we have:
  $$
  b^{tw}_1(b_2^{tw}(g\otimes h))+
  (-1)^{\vert g\vert}b^{tw}_2(g\otimes b_1^{tw}(h))+b^{tw}_2(b_1^{tw}(g)\otimes h)=0.$$
  Thus, we get 
  $g\star f=b^{tw}_2(g\otimes f)=(-1)^{\vert g\vert+1}b^{tw}_1(b_2^{tw}(g\otimes h))\in {\cal I}(\underline{X},\underline{Y})$. 
  
   Similarly, if   $g\in {\cal K}(Z)(\underline{W},\underline{X})$, 
   and $f$ is as above, we have 
     $f\star g\in {\cal I}(\underline{W},\underline{Y})$.
\end{enumerate}
A morphism $f\in {\cal K}(Z)(\underline{X},\underline{Y})$ 
 is called \emph{homologically trivial} iff its class modulo 
 ${\cal I}(Z)(\underline{X},\underline{Y})$ is zero.  
 \end{remark}

\begin{remark}\label{R: desarmado f en estricto + otro simplifica btw's}
 Very often, we can decompose a morphism $f\in \ad(Z)(X,Y)$ as $f=f^0+f^1$, where $f^0,f^1\in \ad(Z)(X,Y)$ and $f^0$ is a strict morphism. Assume this is the case for $n$ composable morphisms $h_1,h_2,\ldots,h_n$ in $\ad(Z)$ where $n>2$. Then, we have 
 $b_n^{ad}(h_n\otimes\cdots \otimes h_2\otimes h_1)=b_n^{ad}(h_n^1\otimes\cdots \otimes h_2^1\otimes h_1^1).$
 In particular, we have the following two situations.
 \begin{enumerate}
  \item If $f:(X,\delta_X)\rightmap{}(Y,\delta_Y)$ is a morphism in $\tw(Z)$ with $f=f^0+f^1$, $\delta_X=\delta_X^0+\delta_X^1$, and $\delta_Y=\delta_Y^0+\delta_Y^1$ as before, we have:
  $$b_1^{tw}(f)=f\circ \delta_X+\delta_Y\circ f+R(f), \hbox{ where }$$
  $R(f)=b_1^{ad}(f^1)+\sum_{\scriptsize\begin{matrix}i_0,i_1\geq 0\\ 
                          i_0+i_1\geq 2\end{matrix}} b^{ad}_{i_0+i_1+1}((\delta_Y^1)^{\otimes i_1}\otimes f^1\otimes (\delta_X^1)^{\otimes i_0})$. 
  \item If $f:(X,\delta_X)\rightmap{}(Y,\delta_Y)$ and $g:(Y,\delta_Y)\rightmap{}(W,\delta_W)$ are morphisms in $\tw(Z)$ with $f=f^0+f^1$, $g=g^0+g^1$, $\delta_X=\delta_X^0+\delta_X^1$,  $\delta_Y=\delta_Y^0+\delta_Y^1$, and $\delta_W=\delta_W^0+\delta_W^1$  as before, we have:
  $$g\star f=b_2^{tw}(g\otimes f)=g\circ f+R(g,f), \hbox{ where }$$
  \hbox{\hskip-.5cm$R(g,f)=\sum_{\scriptsize\begin{matrix}i_0,i_1,i_2\geq 0\\ 
                          i_0+i_1+i_2\geq 1\end{matrix}} b^{ad}_{i_0+i_1+i_2+2}((\delta_W^1)^{\otimes i_2}\otimes g^1\otimes  (\delta_Y^1)^{\otimes i_1}\otimes f^1\otimes (\delta_X^1)^{\otimes i_0})$.}                         
 \end{enumerate}
\end{remark}

In the following, a morphism $f:(X,\delta_X)\rightmap{}(Y,\delta_Y)$ in $\tw(Z)$ is called \emph{strict} iff $f:X\rightmap{}Y$ is a strict morphism of $\ad(Z)$.  

\begin{lemma}\label{R: cuando un estricto esta en Z(hat(B))}\label{R: star-compos=circ-compos para special morfisms} The following holds:
\begin{enumerate}
 \item For any strict morphism $f:(X,\delta_X)\rightmap{}(Y,\delta_Y)$ 
  in $\tw(Z)$,   we have  
$$
b^{tw}_1(f)=\delta_Y\circ f+f\circ \delta_X.$$
Thus, if $f:(X,\delta_X)\rightmap{}(Y,\delta_Y)$ is strict homogeneous morphism of $\tw(Z)$  with degree $-1$, we have: 
$f$ is a morphism in ${\cal Z}(Z)$  iff  $\delta_Y\circ f+f\circ \delta_X=0$. 

\item Given any morphisms $f:(X,\delta_X)\rightmap{}(Y,\delta_Y)$ and $g:(Y,\delta_Y)\rightmap{}(W,\delta_W)$ in $\tw(Z)$, such that at least one of them is strict, then we have $g\star f=g\circ f$.

\item A morphism $f:X=\bigoplus_{i\in I}X_i\rightmap{}\bigoplus_{j\in J}Y_j=Y$ in  $\ad(Z)$ is special (resp. strict) iff the component $f_{j,i}:X_i\rightmap{}Y_j$ is special (resp. strict) for all $i,j$.

\item Every special morphism $f$ of $\ad(Z)$ is strict. 
\end{enumerate}
\end{lemma}

\begin{proof} (1): If $f:(X,\delta_X)\rightmap{}(Y,\delta_Y)$ is strict, with the notation of (\ref{R: desarmado f en estricto + otro simplifica btw's})(1),  we have $R(f)=0$ and we obtain the wanted formula. 
 (2) follows from (\ref{R: desarmado f en estricto + otro simplifica btw's})(2). (3) follows from the formula describing how  $f$ is determined by the  components of its matrix, see (\ref{R: componentes de morfismos homogeneos en ad(Z)}).
\end{proof}

 \begin{lemma}\label{L: transfiriendo diferenciales}
 The following holds: 
 \begin{enumerate}
  \item If a special morphism $h=\sum_uh_u\otimes \frak{e}_u:X\rightmap{}Y$ in $\ad(Z)$ has a two sided inverse $h'$ in $\ad(Z)$ (i.e.  $h':Y\rightmap{}X$ is a morphism in $\ad(Z)$ with $h\circ h'=\hueca{I}_Y$ and $h'\circ h=\hueca{I}_X$) then $h'$ is also special. In this case, we call $h$ a \emph{special isomorphism}. Thus the inverse of a special morphism $h$ is uniquely determined and denoted by $h^{-1}$. Moreover, a morphism $h$ in $\ad(Z)$ is a special isomorphism  iff $h$ is \emph{locally invertible} (i.e.  $h_u:Xe_u\rightmap{}Ye_u$ is a linear isomorphism for all $u\in {\cal P}$).
  
  \item If $h:X\rightmap{}Y$ is a special isomorphism in $\ad(Z)$ and 
  $(X,\delta_X)$ is an object in ${\cal Z}(Z)$, then the pair $(Y,\delta_Y)$,  with $\delta_Y:=-h\circ \delta_X\circ h^{-1}$, is an object of ${\cal Z}(Z)$.  
 Moreover, $h:(X,\delta_X)\rightmap{}(Y,\delta_Y)$ is an isomorphism in ${\cal Z}(Z)$ with inverse $h^{-1}:(Y,\delta_Y)\rightmap{}(X,\delta_X)$. 
 \end{enumerate}
  \end{lemma}

  \begin{proof} First observe that if $g=\sum_ag_a\otimes a$ is a morphism in $\ad(Z)$ with all $a$ not strict units, and $f_1$, $f_2$ are special morphisms, then $f_1\circ g$ and $g\circ f_2$ admit expressions of the same type, that is without involving strict units. 

Assume that $h:X\rightmap{}Y$ is a special morphism with a two sided inverse  $h'$ in $\ad(Z)$. Write $h'=h'_1+h'_2$, where $h'_1$ is a special morphism and $h'_2$ admits an expression not involving strict units. Then, we have $h'_1\circ h+h'_2\circ h=\hueca{I}_X$.  So, we get $h'_1\circ h=\hueca{I}_X$ and $h'_2\circ h=0$. Similarly, we have $h\circ h'_1=\hueca{I}_Y$ and $h\circ h'_2=0$. Therefore, we have $h'_2=\hueca{I}_X\circ h'_2=(h'_1\circ h)\circ h'_2= h'_1\circ (h\circ h'_2)=0$. Therefore, $h'$ is a special morphism. 

Moreover, since $h$ and $h'$ are special, they are of the form $h=L(\underline{h})$ and $h'=L(\underline{h}')$. Since $L$ is a faithful functor, we get that $\underline{h}$ and $\underline{h}'$ are mutual inverses in $\FSMod\g S$. This implies that $h$ is locally invertible.  
So (1) holds. 

For $(2)$, we have   $\vert \delta_X\vert=0$, thus $\delta_X=\sum_a(\delta_X)_a\otimes a$, with $\vert a\vert=0$, for all $a$. Then, 
$\delta_Y=\sum_ah_{v(a)}(\delta_X)_ah_{u(a)}^{-1}\otimes a$, so $\vert \delta_Y\vert=0$. 

By assumption, there is a filtration $0=X_0\subseteq X_1\subseteq\cdots\subseteq X_r=X$ of submodules of $X$ such that $(\delta_X)_a(X_i)\subseteq X_{i-1}$, for all $i$. Consider the filtration 
$0=Y_0\subseteq Y_1\subseteq \cdots\subseteq Y_r=Y$ of the right $S$-module 
$Y$ defined by 
$$Y_i=\sum_{u\in {\cal P}}h_u(X_ie_u), \hbox{ for } i\in [1,r].$$
We have that each $h_u(Xe_u)\subseteq Ye_u$, thus $h_u(X_ie_u)=h_u(X_ie_u)e_u$ is an $S$-submodule of $Y$, and then so is $Y_i$. Now, let $y_i=h_u(x_ie_u)\in Y_i$ be a generator of $Y_i$, with $u\in {\cal P}$ and $x_i\in X_i$. Then,  we have 
$$(\delta_Y)_a(y_i)=(\delta_Y)_a(h_u(x_ie_u))=(\delta_Y)_a(h_u(x_ie_u)e_u).$$
The last expression is zero if $u(a)\not=u$.  If $u=u(a)$, we have 
$$\begin{matrix}
(\delta_Y)_a(y_i)&=&(\delta_Y)_a(h_{u(a)}(x_ie_{u(a)}))\hfill\\
&=&h_{v(a)}(\delta_X)_ah^{-1}_{u(a)}(h_{u(a)}(x_ie_{u(a)}))\hfill\\
&=&
h_{v(a)}(\delta_X)_a(x_ie_{u(a)})\in  h_{v(a)}(X_{i-1})\subseteq Y_{i-1},\hfill
  \end{matrix}$$

Moreover,  from (\ref{L: morf especiales vs operaciones de ad(Z)}), for $s\geq 0$, we have 
$$\begin{matrix}
b_s^{ad}(\delta_Y^{\otimes s})
&=&
(-1)^sb_s^{ad}(h\circ\delta_X\circ h^{-1}\otimes h\circ \delta_X\circ h^{-1}\otimes \cdots\otimes h\circ \delta_X\circ h^{-1})\hfill\\
&=&
(-1)^sh\circ b_s^{ad}(\delta_X\circ h^{-1}\otimes h\circ \delta_X\circ h^{-1}\otimes \cdots\otimes h\circ \delta_X\circ h^{-1})\hfill\\
&=&
(-1)^s(-1)^{s-1}h\circ b_s^{ad}(\delta_X\otimes \delta_X\otimes \cdots\otimes \delta_X\circ h^{-1})\hfill\\
&=&
(-1)^s(-1)^{s-1}(-1)h\circ b_s^{ad}(\delta_X\otimes \delta_X\otimes \cdots\otimes \delta_X)\circ h^{-1}\hfill\\
&=&
h\circ b_s^{ad}(\delta_X^{\otimes s})\circ h^{-1}\hfill\\
\end{matrix}$$
Hence, we obtain 
$\sum_sb_s^{ad}(\delta_Y^{\otimes s})=h\circ \sum_sb^{ad}_s(\delta_X^{\otimes s})\circ h^{-1}=0.$
 Then, we have that $(Y,\delta_Y)$ is an object of ${\cal Z}(Z)$.  

  Finally, we have that  $h:(X,\delta_X)\rightmap{}(Y,\delta_Y)$ is a  morphism in ${\cal Z}(Z)$ because it is strict and satisfies 
  $$\begin{matrix}
     \delta_Y\circ h+h\circ \delta_X
     &=&
     -(h\circ \delta_X\circ h^{-1})\circ h+h\circ \delta_X\hfill\\
     &=&
     -((h\circ \delta_X)\circ h^{-1})\circ h+h\circ \delta_X=-h\circ \delta_X+h\circ \delta_X=0.\hfill\\
    \end{matrix}$$
\end{proof}

 In the following, we say that a candidate ${\cal C}$ to be a category is a  \emph{precategory} iff all the requirements of a category are satisfied by ${\cal C}$, with the only exception of the associativity of the composition.

\begin{proposition}\label{P: la cat Z(B)} Given a  unitary strict $b$-algebra $Z=(Z,\{b_n\}_{n\in \hueca{N}})$, 
   the following elements determine a \emph{precategory} ${\cal Z}(Z)$. 
   Its objects are those of $\tw(Z)$. The morphisms $f:(X,\delta_X)\rightmap{}(Y,\delta_Y)$ of ${\cal Z}(Z)$ are the morphisms in $\ad(Z)(X,Y)$ with degree $\vert f\vert=-1$ such that $b_1^{tw}(f)=0$. 
   
   Given morphisms $f:(X,\delta_X)\rightmap{}(Y,\delta_Y)$ and $g:(Y,\delta_Y)\rightmap{}(W,\delta_W)$ in ${\cal Z}(Z)$, its composition is defined by
   $g\star f=b_2^{tw}(g\otimes f)$.

   The quotient precategory ${\cal H}(Z)$, obtained from ${\cal Z}(Z)$ as the quotient modulo the ideal ${\cal I}=b_1^{tw}[\tw(Z)(-,?)_{-2}]$ of ${\cal Z}(Z)$,  is a category. 
  \end{proposition} 
  
  \begin{proof} By (\ref{R: notacion star})(1), the composition $\star$ of the precategory ${\cal Z}(Z)$ is well defined. From (\ref{R: notacion star})(2), we get that ${\cal I}$ is indeed an ideal of 
   the precategory ${\cal Z}(Z)$, which implies that the composition in the quotient precategory ${\cal H}(Z)$ is well defined.

  Again, from the fact that $\tw(Z)$ is a $b$-category, we have
  $$\begin{matrix}
  0&=&
  b^{tw}_3(b^{tw}_1\otimes id^{\otimes 2})+b^{tw}_3(id\otimes b^{tw}_1\otimes id)+
   b^{tw}_3(id^{\otimes 2}\otimes b^{tw}_1)\hfill\\
      &&+b^{tw}_2(b^{tw}_2\otimes id)+b^{tw}_2(id\otimes b^{tw}_2)+b^{tw}_1b^{tw}_3.\hfill\\
     \end{matrix}$$
   From this equation we obtain that, modulo the ideal ${\cal I}$, we indeed have the associativity property for the composition in the quotient precategory ${\cal H}(Z)$. 
   \medskip
   
   In the following, we use that $Z$ is unitary strict.  For $\underline{X}\in {\cal Z}(Z)$, we consider the special morphism  
   $\hueca{I}_{X}=\sum_{j\in {\cal P}}id_{Xe_j}\otimes \frak{e}_j\in \Hom_{\tw(Z)}(\underline{X},\underline{X}).$ 
    
   We have that $\hueca{I}_{\underline{X}}$ belongs to ${\cal Z}(Z)$ because, from (\ref{R: cuando un estricto esta en Z(hat(B))})(1), we have 
   $ b_1^{tw}(\hueca{I}_{X})= \delta_X\circ \hueca{I}_X+\hueca{I}_X\circ \delta_X=-\delta_X+\delta_X=0$.

 Now, given $t\in \tw(Z)(\underline{X},\underline{Y})_{-1}=\ad(Z)(X,Y)_{-1}$, from (\ref{R: cuando un estricto esta en Z(hat(B))})(2), we have 
 $t\star \hueca{I}_X=t\circ \hueca{I}_X=t$ and $\hueca{I}_Y\star t=\hueca{I}_Y\circ t=t$. 
 \end{proof}

\section{Special and canonical conflations in ${\cal Z}(Z)$}\label{special confl in Z(Z)}

We keep the preceding terminology, where $Z$ is a $b$-algebra over the elementary  algebra $S$, with enough idempotents $\{e_u\}_{u\in {\cal P}}$,  and we assume that it is unitary strict with  strict units $\{\frak{e}_u\}_{u\in {\cal P}}$, as in (\ref{D: unitary strict}). We have the associated $b$-category  $\ad(Z)$ over $S$, as in (\ref{R: recordatorio de ad(Z)}), and   
 a fixed basis ${\hueca{B}}$ for the vector space $Z$ formed by homogeneous directed elements, and containing the strict units of $Z$. 
 
 Then, we have the $b$-category $\tw(Z)$ reminded in (\ref{R: record de tw(Z)}). Recall that,  given  two morphisms $f:(X,\delta_X)\rightmap{}(Y,\delta_Y)$ and $g:(Y,\delta_Y)\rightmap{}(W,\delta_W)$  in $\tw(Z)$, we  use the notation $g\star f=b_2^{tw}(g\otimes f)$. Then, we have the precategory ${\cal Z}(Z)$ with composition $\star$.  
 
 In this section, we introduce a special class of pairs of morphisms in ${\cal Z}(Z)$, which we call \emph{special conflations} because they have properties which are similar to those of conflations of exact structures on  additive categories.

\begin{lemma}\label{L: diferenciales de la suma directa en tw(hat(B))}
Let $E=X\oplus Y$ be a decomposition of a right $S$-module $E$ and $\delta_E:E\rightmap{}E$ a morphism of degree $0$ in $\ad(Z)$ with matrix form 
$$\delta_E=\begin{pmatrix}
            \delta_X&\gamma\\
            0&\delta_Y\\
           \end{pmatrix}$$
           associated to the given decomposition of $E$. Then, 
           $$(E,\delta_E)\in \tw(Z)\hbox{ \  iff \ } (X,\delta_X),(Y,\delta_Y)\in \tw(Z)\hbox{  and }b_1^{tw}(\gamma)=0.$$ 
\end{lemma} 

 \begin{proof}  Assume first that  $(E,\delta_E)$ belongs to $\tw(Z)$, and let us show that the pairs $(X,\delta_X)$ and $(Y,\delta_Y)$ belong to  $\tw(Z)$. Suppose that 
 $\delta_E=\sum_{a\in \hueca{B}}(\delta_E)_a\otimes a$. From (\ref{R: componentes de morfismos homogeneos en ad(Z)}), we have  that  $(\delta_E)_a$ has the matrix form  $(\delta_E)_a=\begin{pmatrix}
      (\delta_X)_a&\gamma_a\\ 0&(\delta_Y)_a\\                                                           
      \end{pmatrix}$, 
with $\delta_X=\sum_{a\in \hueca{B}}(\delta_X)_a\otimes a$,  $\delta_Y=\sum_{a\in \hueca{B}}(\delta_Y)_a\otimes a$, and 
$\gamma=\sum_{a\in \hueca{B}}\gamma_a\otimes a$. 
 
We have a right $S$-module filtration $0=E_0\subseteq \cdots\subseteq  E_{l-1}\subseteq E_l=E$ such that $(\delta_E)_a(E_i)\subseteq E_{i-1}$, for all $a\in \hueca{B}$ and $i\in [1,l]$. 

If we define $X_i:=X\cap E_i$, for all $i$, we obtain the filtration $0=X_0\subseteq \cdots \subseteq X_{l-1}\subseteq X_l=X$. Given $x\in X_i$ and $a\in \hueca{B}$, we have $(\delta_X)_a(x)=(\delta_E)_a(x)\in X\cap E_{i-1}=X_{i-1}$, for all $i$. So, $(\delta_X)_a(X_i)\subseteq X_{i-1}$, for all $a\in \hueca{B}$ and $i\in [1,l]$ . 

If we define $Y_i=\pi_2(E_i)$, where $\pi_2:E\rightmap{}Y$ is the second projection, we obtain the filtration $0=Y_0\subseteq \cdots\subseteq Y_{l-1}\subseteq Y_l=Y$. Given $y\in Y_i$, there is $x\in X$ with $x+y\in E_i$. So, for $a\in \hueca{B}$, we have 
$(\delta_E)_a(x+y)=(\delta_X)_a(x)+\gamma_a(y)+(\delta_Y)_a(y)\in E_{i-1}$. Here, $(\delta_X)_a(x)+\gamma_a(y)\in X$ and, therefore, 
$(\delta_Y)_a(y)=\pi_2((\delta_E)_a(x+y))\in \pi_2(E_{i-1})=Y_{i-1}$. Therefore, $(\delta_Y)_a(Y_i)\subseteq Y_{i-1}$, for all $a\in \hueca{B}$ and $i\in [1,l]$. 
 
 From (\ref{R: descomposicon de morfismos en ad(hat(B))}), we have
 $0=\sum_{s\geq 1}b^{ad}_s(\delta_E^{\otimes s})=\begin{pmatrix}
                                                 \sum_sb_s^{ad}(\delta_X^{\otimes s})&b_1^{tw}(\gamma)\\
                                                 0& \sum_sb_s^{ad}(\delta_Y^{\otimes s})
                                                \end{pmatrix},$
so $(X,\delta_X), (Y,\delta_Y)\in \tw(Z)$ and $b_1^{tw}(\gamma)=0$.

Now, assume that $(X,\delta_X), (Y,\delta_Y)\in \tw(Z)$, and $b_1^{tw}(\gamma)=0$, and look at the canonical descriptions  $\delta_X=\sum_a(\delta_X)_a\otimes a$, $ \delta_Y=\sum_a(\delta_Y)_a\otimes a$, and 
$\gamma=\sum_a\gamma_a\otimes a$. 
Then, we have $\delta_E=\sum_a(\delta_E)_a\otimes a$, where 
$$(\delta_E)_a=\begin{pmatrix}
                (\delta_X)_a&\gamma_a\\ 0& (\delta_Y)_a\\
               \end{pmatrix}.$$
 By assumption, we have  filtrations of right $S$-submodules 
$$0=X_0\subseteq X_1\subseteq \cdots \subseteq X_r=X
  \hbox{  \ and \  } 
  0=Y_0\subseteq Y_1\subseteq \cdots \subseteq Y_t=Y,
  $$
  such that $(\delta_X)_a(X_i)\subseteq X_{i-1}$ and $(\delta_Y)_a(Y_j)\subseteq Y_{j-1}$, for all $i$ and  $j$. Consider the filtration
  $$0=X_0\oplus 0\subseteq X_1\oplus 0\subseteq \cdots\subseteq X_r\oplus 0\subseteq X\oplus Y_1\subseteq \cdots\subseteq X\oplus Y_t=X\oplus Y,$$
  which clearly satisfies $(\delta_E)_a(X_i\oplus 0)\subseteq X_{i-1}\oplus 0$ and $(\delta_E)_a(X\oplus Y_j)\subseteq X\oplus Y_{j-1}$. Since $\sum_{s\geq 1}b_s^{ad}(\delta_X^{\otimes s})=0$, $\sum_{s\geq 1}b_s^{ad}(\delta_Y^{\otimes s})=0$, and  $b_1^{tw}(\gamma)=0$, we also have $\sum_{s\geq 1}b_s^{ad}(\delta_E^{\otimes s})=0$.
 \end{proof}

 \begin{definition}\label{D: special conflations in Z(hat(B))}
 A \emph{special conflation}  is a sequence of morphisms in ${\cal Z}(Z)$  
 $$(X,\delta_X)\rightmap{f}(E,\delta_E)\rightmap{g}(Y,\delta_Y)$$
 formed by special morphisms $f=\sum_{u\in {\cal P}}f_u\otimes \frak{e}_u$ and $g=\sum_{u\in {\cal P}}g_u\otimes \frak{e}_u$ such that the sequence of vector spaces 
 $$0\rightmap{}Xe_u\rightmap{f_u}Ee_u\rightmap{g_u}Ye_u\rightmap{}0$$
 is exact for all $u\in {\cal P}$.
 
 A \emph{special inflation (resp. special deflation)} $f:(X,\delta_X)\rightmap{}(E,\delta_E)$  (resp. $g:(E,\delta_E)\rightmap{}(Y,\delta_Y)$) in ${\cal Z}({Z})$ is a special morphism  for which there is a special conflation $(X,\delta_X)\rightmap{f}(E,\delta_E)\rightmap{g}(Y,\delta_Y)$ in ${\cal Z}({Z})$. 
 \end{definition}

  \begin{lemma}\label{L: previo a forma canonica de conflaciones}
 Assume that we have a pair of composable morphisms in $\tw(Z)$:
 $$(X,\delta_X)\rightmap{ \ f \ }(E,\delta_E)\rightmap{ \ g \ }(Y,\delta_Y),$$
 where 
 $E=E^1\oplus E^2$, $f=(\tilde{f},0)^t$, $g=(0,\tilde{g})$, where  $\tilde{f}:X\rightmap{}E^1$ and $\tilde{g}:E^2\rightmap{}Y$ are special isomorphisms in $\ad(Z)$. Then,  the morphisms $f$ and $g$ belong to ${\cal Z}(Z)$ iff the morphism $\delta_E$ has the triangular form 
 $$\delta_{E}=
 \begin{pmatrix}
 -\tilde{f}\circ\delta_X\circ \tilde{f}^{-1}&\gamma\\ 0&-\tilde{g}^{-1}\circ\delta_Y\circ \tilde{g}\\
 \end{pmatrix},$$ for some homogeneous morphism 
 $\gamma:E^2\rightmap{}E^1$ in $\ad(Z)$ of degree $0$. In this case, if we define $\delta_{E^1}:=-\tilde{f}\circ \delta_X\circ \tilde{f}^{-1}$ and $\delta_{E^2}:=-\tilde{g}^{-1}\circ \delta_Y\circ \tilde{g}^{-1}$, from (\ref{L: transfiriendo diferenciales})(2), we obtain objects $(E^1,\delta_{E^1})$ and $(E^2,\delta_{E^2})$ in ${\cal Z}(Z)$, and $\delta_E=\begin{pmatrix}
            \delta_{E^1}&\gamma\\ 0& \delta_{E^2}\\                                                                                                                                                                                                                                                                                                                                                                                                                                                                                  \end{pmatrix}$.
 \end{lemma}
 
 \begin{proof} Assume that $f$ and $g$  are morphisms in ${\cal Z}(Z)$. 
 We have  $$\delta_E=\begin{pmatrix}
                            \alpha_{1,1}&\alpha_{1,2}\\
                            \alpha_{2,1}&\alpha_{2,2}\\
                           \end{pmatrix}$$
                           where $\alpha_{i,j}:E^j\rightmap{}E^i$ 
are morphisms in $\ad(Z)$ with degree 0. 
Since $f = (\tilde{f} , 0)^t : (X, \delta_X) \rightmap{} (E, \delta_E )$  is a strict morphism of  ${\cal Z}(Z)$, from (\ref{R: star-compos=circ-compos para special morfisms}), we have  
$0 = b^{tw}_1(f)=\delta_E\circ f+f\circ \delta_X=\delta_E\circ (\tilde{f},0)^t+(\tilde{f},0)^t\circ \delta_X$.
Then, we have 
$$0=\begin{pmatrix}
     \alpha_{1,1}\circ \tilde{f}\\ \alpha_{2,1}\circ \tilde{f}\\
    \end{pmatrix} + 
    \begin{pmatrix}
     \tilde{f}\circ \delta_X\\ 0\\
    \end{pmatrix}.$$
Hence, we obtain $\alpha_{2,1}\circ \tilde{f}=0$ and  
 $\alpha_{1,1}\circ \tilde{f}=-\tilde{f}\circ \delta_X$. Therefore, we have 
 $0=(\alpha_{2,1}\circ \tilde{f})\circ \tilde{f}^{-1}=-\alpha_{2,1}\circ (\tilde{f}\circ \tilde{f}^{-1})=-\alpha_{2,1}\circ \hueca{I}_X=\alpha_{2,1}$; and, also,  
 $\alpha_{1,1}=(\alpha_{1,1}\circ \tilde{f})\circ \tilde{f}^{-1}=-(\tilde{f}\circ \delta_X)\circ \tilde{f}^{-1}$. 

Since $g=(0,\tilde{g}):(E,\delta_E)\rightmap{}(Y,\delta_Y)$ is also strict, we have $0=b_1^{tw}(g)=
\delta_Y\circ(0,\tilde{g})+(0,\tilde{g})\circ \delta_E$, and we obtain
$$0=\begin{pmatrix}
    0,\delta_Y\circ \tilde{g} 
    \end{pmatrix}+
    \begin{pmatrix}
     \tilde{g}\circ \alpha_{2,1},\tilde{g}\circ \alpha_{2,2}\\
    \end{pmatrix}
    =
\begin{pmatrix}
   0,\delta_Y\circ \tilde{g}+\tilde{g}\circ \alpha_{2,2}
  \end{pmatrix}.$$
Hence, we have  $\tilde{g}^{-1}\circ(\delta_Y\circ \tilde{g})=-\tilde{g}^{-1}\circ (\tilde{g}\circ \alpha_{2,2})=-(\tilde{g}^{-1}\circ \tilde{g})\circ \alpha_{2,2}=-\alpha_{2,2}$. From (\ref{L: transfiriendo diferenciales}), we can omit the parenthesis. 

Conversely, notice that if $\delta_E$ has the matrix form described in the statement of this lemma, we can reverse the preceding argument to obtain that $b_1^{tw}(f)=0$ and $b_1^{tw}(g)=0$. Thus, the morphisms $f$ and $g$ are morphisms
 of ${\cal Z}(Z)$.
 \end{proof}
  
  Special inflations and deflations can be characterized by the following.
  
  \begin{lemma}\label{L: caracterizacion de special infl y defl} 
 We have:
 \begin{enumerate}
  \item 
 A special morphism $f=\sum_{u\in {\cal P}}f_u\otimes \frak{e}_u:(X,\delta_X)\rightmap{}(E,\delta_E)$  in ${\cal Z}(Z)$ is a special inflation iff each $f_u:Xe_u\rightmap{}Ee_u$ is a linear monomorphism. 
 \item A special morphism $g=\sum_{u\in {\cal P}}g_u\otimes \frak{e}_u:(E,\delta_E)\rightmap{}(Y,\delta_Y)$  in ${\cal Z}(Z)$ is a special deflation  iff each $g_u:Ee_u\rightmap{}Ye_u$ is a linear epimorphism. 
 \end{enumerate}   
 \end{lemma}
 
 \begin{proof} We only prove (1), since the proof of (2) is similar.  
  Assume that $f_u:Xe_u\rightmap{}Ee_u$ is injective for each $u\in {\cal P}$. Then, there is a decomposition of vector spaces $Ee_u=E_u^1\oplus E_u^2$ such that $f_u=(\tilde{f}_u,0)^t$, where 
  $\tilde{f}_u:Xe_u\rightmap{}E^1_u$ is a linear isomorphism. Now, consider $E^1=\bigoplus_{u\in {\cal P}}E_u^1$ and $E^2=\bigoplus_{u\in {\cal P}}E^2_u$. Both, $E^1$ and $E^2$ have natural structures of right ${S}$-modules with $E^ie_u=E^i_u$, for $u\in {\cal P}$ and $i\in \{1,2\}$. Clearly, $\tilde{f}:=\sum_u\tilde{f}_u\otimes \frak{e}_u:X\rightmap{}E^1$ is a special isomorphism in $\ad(Z)$ with $f=(\tilde{f},0)^t:X\rightmap{}E$ in $\ad(Z)$.  From the first part of the proof of (\ref{L: previo a forma canonica de conflaciones}), we know that $\delta_E$ has triangular form
  $\delta_E=\begin{pmatrix} \alpha_{1,1}&\alpha_{1,2}\\ 0&\alpha_{2,2}                                                                                                                                                                                                                                                                                                                                                                                                                                                                                                                                                                                                                                                                                                                                                                                                                                                                                                                                                                                        \end{pmatrix}$. Then, from (\ref{L: diferenciales de la suma directa en tw(hat(B))}), we have $(E^2,\alpha_{2,2})\in {\cal Z}(Z)$. So we have the special conflation
  $$(X,\delta_X)\rightmap{ \ f \ }(E,\delta_E)\rightmap{ \ (0,\hueca{I}_{E^2}) \ }(E^2,\alpha_{2,2}).$$
 \end{proof}

 \begin{definition}\label{D: equiv sequences in Z(hat(B))} Consider the following relation in the class of composable pairs of morphisms in ${\cal Z}(Z)$. 
 Given the composable pairs in ${\cal Z}(Z)$
 $$\xi\hbox{\,}:\hbox{\,}(X,\delta_X)\rightmap{f}(E,\delta_E)\rightmap{g}(Y,\delta_Y) \hbox{ and }
 \xi'\hbox{\,}:\hbox{\,}(X,\delta_X)\rightmap{f'}(E',\delta_{E'})\rightmap{g'}(Y,\delta_Y),$$
 we write  $\xi\rightmap{h\simeq}\xi'$ whenever $h:(E,\delta_E)\rightmap{}(E',\delta_{E'})$ is an isomorphism in ${\cal Z}(Z)$ such that  the following diagram 
 commutes in ${\cal Z}(Z)$
 $$\begin{matrix}
   (X,\delta_X)&\rightmap{f}&(E,\delta_E)&\rightmap{g}&(Y,\delta_Y)\hfill\\
   \shortlmapdown{\hueca{I}_X}&&\shortrmapdown{h}&&\shortrmapdown{\hueca{I}_Y}\\
   (X,\delta_X)&\rightmap{f'}&(E',\delta_{E'})&\rightmap{g'}&(Y,\delta_Y).\hfill\\
   \end{matrix}$$
    
   We will simply write $\xi\rightmap{\simeq}\xi'$, when there is a morphism  $h$ such that $\xi\rightmap{h\simeq}\xi'$.
   \end{definition}

We will show below that ``$\rightmap{\simeq}$'' is an equivalence relation in the class of special conflations.  For this we need to look closer into the structure of special conflations and the following canonical representatives.  
  
\begin{definition} A \emph{canonical conflation}  
$\xi:(X,\delta_X)\rightmap{f}(E,\delta_E)\rightmap{g}(Y,\delta_Y)$ in ${\cal Z}(Z)$
is a special conflation such that $E=X\oplus Y$ as $S$-modules, $f=(\hueca{I}_X,0)^t$, $g=(0,\hueca{I}_Y)$, and $\delta_E=\begin{pmatrix}
              \delta_X& \gamma\\ 0&\delta_Y\\                                                                                                                                                                                                     \end{pmatrix}$, for some 
some homogeneous morphism 
 $\gamma:Y\rightmap{}X$ in $\ad(Z)$ of degree $0$.
 
 A \emph{canonical inflation (resp. canonical deflation)} $f:(X,\delta_X)\rightmap{}(E,\delta_E)$  (resp. $g:(E,\delta_E)\rightmap{}(Y,\delta_Y)$) in ${\cal Z}(Z)$ is a special morphism  for which there is a canonical conflation $(X,\delta_X)\rightmap{f}(E,\delta_E)\rightmap{g}(Y,\delta_Y)$ in ${\cal Z}(Z)$. 
\end{definition}

 \begin{lemma}\label{L: may assume special conflation direct sum middle term}
  For any special conflation  
 $$\xi:\hbox{ \ }(X,\delta_X)\rightmap{f}(E,\delta_E)\rightmap{g}(Y,\delta_Y)$$
 of ${\cal Z}(Z)$ there is  a canonical conflation   
 $\overline{\xi}:(X,\delta_X)\rightmap{ \ \overline{f} \ }(\overline{E},\delta_{\overline{E}})\rightmap{ \ \overline{g} \ }(Y,\delta_Y)$ of ${\cal Z}(Z)$
 and a special isomorphism $h:(E,\delta_E)\rightmap{}(\overline{E},\delta_{\overline{E}})$ in ${\cal Z}(Z)$ such that  $$\xi\rightmap{h\simeq}\overline{\xi}\hbox{ \ \  and \ \  } \overline{\xi}\rightmap{ \ h^{-1}\simeq \ }\xi.$$ 
 \end{lemma}
 
 \begin{proof} The given special conflation $\xi$ is formed by special morphisms $f=\sum_{u\in {\cal P}}f_u\otimes \frak{e}_u$ and $g=\sum_{u\in {\cal P}}g_u\otimes \frak{e}_u$ such that the sequence 
 $$0\rightmap{}Xe_u\rightmap{f_u}Ee_u\rightmap{g_u}Ye_u\rightmap{}0$$
 is exact for all $u\in {\cal P}$.
 There are commutative diagrams of linear maps with exact rows 
  $$\begin{matrix}
    0&\rightmap{}&Xe_u&\rightmap{f_u}&Ee_u&\rightmap{g_u}&Ye_u&\rightmap{}&0\\
    &&\shortlmapdown{Id_{Xe_u}}&&\shortrmapdown{h_u}&&\shortrmapdown{Id_{Ye_u}}\\
    0&\rightmap{}&Xe_u&\rightmap{(Id_{Xe_u},0)^t}&Xe_u\oplus Ye_u&\rightmap{(0,Id_{Ye_u})}&Ye_u&\rightmap{}&0,\\
   \end{matrix}$$
 where $h_u$ is a linear isomorphism. Consider the special isomorphism 
 $h:=\sum_uh_u\otimes \frak{e}_u\in \ad(Z)(E,\overline{E})$, where $\overline{E}= X\oplus Y$. From (\ref{L: transfiriendo diferenciales}), we know that $\delta_{\overline{E}}:=-h\circ \delta_E\circ h^{-1}\in \ad(Z)(\overline{E},\overline{E})_0$ is a morphism such that $(\overline{E},\delta_{\overline{E}})$ is an object of ${\cal Z}(Z)$. We have a diagram of special morphisms which commutes in $\ad(Z)$:
 $$\begin{matrix}
    (X,\delta_X)&\rightmap{f}&(E,\delta_E)&\rightmap{g}&(Y,\delta_Y)\\
    \shortrmapdown{\hueca{I}_X}&&\shortrmapdown{h}&&\shortrmapdown{\hueca{I}_Y}\\
    (X,\delta_X)&\rightmap{\overline{f}}&(\overline{E},\delta_{\overline{E}})&\rightmap{\overline{g}}&(Y,\delta_Y),\\
   \end{matrix}$$
   where $\overline{f}=(\hueca{I}_X , 0)^t$ and $\overline{g}=(0,\hueca{I}_Y)$.  
  From (\ref{L: transfiriendo diferenciales}), we know that the special morphism 
  $h:(E,\delta_E)\rightmap{}(\overline{E},\delta_{\overline{E}})$ is an isomorphism in  ${\cal Z}(Z)$ with inverse 
  $h^{-1}=\sum_uh_u^{-1}\otimes e_u$ in ${\cal Z}(Z)$.
   We also have that $\overline{f}$ and $\overline{g}$ are morphisms of ${\cal Z}(Z)$, because from (\ref{R: star-compos=circ-compos para special morfisms}) we have $\overline{f}=h\circ f=h\star f$ and $\overline{g}=g\circ h^{-1}=g\star h^{-1}$. 
  Moreover, the preceding diagram commutes in ${\cal Z}(Z)$.
  Then, from (\ref{L: previo a forma canonica de conflaciones}) applied to the conflation $\overline{\xi}$, we obtain that 
 $\delta_{\overline{E}}=
 \begin{pmatrix}
 \delta_X&\gamma\\ 0&\delta_Y\\
 \end{pmatrix}$, for some homogeneous morphism 
 $\gamma:Y\rightmap{}X$ in $\ad(Z)$ of degree $0$, as  wanted.  
 \end{proof}

 \begin{lemma}\label{R: caract de infl (defl) canonic, compos de infl (defl) canonic es infl (defl) canonic} 
 The following holds. 
 \begin{enumerate}
  \item A morphism $f:(X,\delta_X)\rightmap{}(E,\delta_E)$ in ${\cal Z}(Z)$ is a canonical inflation iff $E=X\oplus Y$ as right $S$-modules and $f=(\hueca{I}_X,0)^t$.
  \item A morphism $g:(E,\delta_E)\rightmap{}(Y,\delta_Y)$ in ${\cal Z}(Z)$ is a canonical deflation iff $E=X\oplus Y$ as right $S$-modules and $g=(0,\hueca{I}_Y)$.
  \item Composition of canonical inflations (resp. canonical deflations) is a canonical inflation (resp. a canonical deflation). 
\end{enumerate}
 \end{lemma}

\begin{proof} 
Indeed, in order to verify (1), take a  morphism $f: (X,\delta_X)\rightmap{}(E,\delta_E)$ in ${\cal Z}(Z)$ is such that $E=X \oplus Y$ and $f=(\hueca{I}_X,0)^t$. Then, as in the first part of the proof of (\ref{L: previo a forma canonica de conflaciones}), we get 
 $\delta_E=\begin{pmatrix}
            \delta_X&\gamma\\ 0&\delta_Y\\
           \end{pmatrix},$
for some homogeneous morphisms $\gamma:Y\rightmap{}X$ and $\delta_Y:Y\rightmap{}Y$ in $\ad(Z)$ with zero degree. Then, from (\ref{L: diferenciales de la suma directa en tw(hat(B))}), we know that $(Y,\delta_Y)$ belongs to ${\cal Z}(Z)$.  Finally, from  (\ref{L: previo a forma canonica de conflaciones}), we obtain the canonical conflation 
$$(X,\delta_X)\rightmap{ \ f \ }(E,\delta_E)\rightmap{ \ (0,\hueca{I}_Y) \ }(Y,\delta_Y).$$

(2) is verified similarly, and (3) follows from (1) and (2).  
\end{proof}

 \begin{lemma}\label{L: isos matriciales entre terminos de enmedio}
 Let $(X,\delta_X)$, $(Y,\delta_Y)$, $(X',\delta_{X'})$, $(Y',\delta_{Y'})$, $E=(X\oplus Y,\delta_E)$, and $E'=(X'\oplus Y',\delta_{E'})$ be objects in ${\cal Z}(Z)$,  with  
 $$\delta_{E}=
 \begin{pmatrix}
 \delta_X&\gamma\\ 0&\delta_Y\\
 \end{pmatrix} \hbox{ and \ }
 \delta_{E'}=
 \begin{pmatrix}
 \delta_{X'}&\gamma'\\ 0&\delta_{Y'}\\
 \end{pmatrix},$$
 for homogeneous morphisms 
 $\gamma:Y\rightmap{}X$ and $\gamma':Y'\rightmap{}X'$ in $\ad(Z)$ of degree $0$.
   Suppose that $h:E\rightmap{}E'$ is a morphism in $\ad(Z)$ with  matrix form $$h=\begin{pmatrix}
      h_{1,1}&s\\ 0&h_{2,2}\\
     \end{pmatrix}$$
   with degree $-1$, where  $h_{1,1}:X\rightmap{}X'$ and  $h_{2,2}:Y\rightmap{}Y'$  are special isomorphisms in ${\cal Z}(Z)$.
   Then, $h$  belongs to ${\cal Z}(Z)(E,E')$  iff $b_1^{tw}(s)=-\gamma'\circ h_{2,2}-h_{1,1}\circ\gamma$.
\end{lemma}     

\begin{proof} Write $h=h_0+h_1$, where 
 $h_0=\begin{pmatrix}
      h_{1,1}&0\\ 0&h_{2,2}\\
      \end{pmatrix}$ and 
 $h_1=\begin{pmatrix}
      0&s\\ 0&0\\
      \end{pmatrix}$. Then, $h_0$ is a special morphism and, from (\ref{R: cuando un estricto esta en Z(hat(B))})(1),  we have 
  $$\begin{matrix}
   b_1^{tw}(h_0)
    &=&
   
    \begin{pmatrix}
     \delta_{X'} &\gamma'\\
     0&\delta_{Y'}\\
     \end{pmatrix}\circ
      \begin{pmatrix}
      h_{1,1}&0\\ 0&h_{2,2}\\
     \end{pmatrix}
   +
    \begin{pmatrix}
      h_{1,1}&0\\ 0&h_{2,2}\\
     \end{pmatrix}\circ
   \begin{pmatrix}
     \delta_{X} &\gamma\\
     0&\delta_{Y}\\
     \end{pmatrix}\hfill\\
     &=&
     \begin{pmatrix}
     \delta_{X'}\circ h_{1,1}+h_{1,1}\circ \delta_X& \gamma'\circ h_{2,2}+h_{1,1}\circ\gamma\\
     0&\delta_{Y'}\circ h_{2,2}+h_{2,2}\circ \delta_Y
     \end{pmatrix}\hfill\\
     &=&
     \begin{pmatrix}
      b_1^{tw}(h_{1,1})&\gamma'\circ h_{2,2}+h_{1,1}\circ\gamma\\
      0&b_1^{tw}(h_{2,2})
     \end{pmatrix}=
     \begin{pmatrix}
      0&\gamma'\circ h_{2,2}+h_{1,1}\circ\gamma\\
      0&0
     \end{pmatrix}
     ,\hfill\\
  \end{matrix}
   $$
   while
   $b_1^{tw}(h_1)=\begin{pmatrix}
            0&b_1^{tw}(s)\\
            0&0\\
           \end{pmatrix}.$ It follows that $b_1^{tw}(h)=0$ if and only if $b_1^{tw}(s)=-\gamma'\circ h_{2,2}-h_{1,1}\circ\gamma$, as claimed.   
\end{proof}

 \begin{lemma}\label{L: b1tw(conjugado(s))=conjugado(b1tw(s))}
 Let $h:(X',\delta_{X'})\rightmap{}(X,\delta_{X})$ and $h':(Y',\delta_{Y'})\rightmap{}(Y,\delta_{Y})$  be special morphisms in ${\cal Z}(Z)$. Then, for any 
 morphism $s:Y\rightmap{}X'$ in $\ad(Z)$, with degree $\vert s\vert=-1$, we have 
$b_1^{tw}(h\circ s\circ h')=-h\circ b_1^{tw}(s)\circ h'$. 
 \end{lemma}
 
 \begin{proof} Set $\Delta:=b^{ad}_{i_0+i_1+1}(\delta_X^{\otimes i_1}\otimes h\circ s\circ h'\otimes \delta_{Y'}^{\otimes i_0})$. Since $h$ and $h'$ are special morphisms in ${\cal Z}(Z)$, from (\ref{R: cuando un estricto esta en Z(hat(B))})(1), we have   $\delta_X\circ h=-h\circ\delta_{X'}$ and $h'\circ \delta_{Y'}=-\delta_Y\circ h'$. From (\ref{L: morf especiales vs operaciones de ad(Z)}), we obtain the following  equalities
$$\begin{matrix}
    \Delta&=&b^{ad}_{i_0+i_1+1}(\delta_X^{\otimes i_1}\otimes h\circ s\circ h'\otimes \delta_{Y'}^{\otimes i_0})\hfill\\
    &=&
    -b^{ad}_{i_0+i_1+1}(\delta_X^{\otimes (i_1-1)}\otimes \delta_X\circ h\otimes s\circ h'\otimes \delta_{Y'}^{\otimes i_0})\hfill\\
    &=&
     b^{ad}_{i_0+i_1+1}(\delta_X^{\otimes (i_1-1)}\otimes h\circ \delta_{X'}\otimes s\circ h'\otimes \delta_{Y'}^{\otimes i_0})=\cdots\hfill\\
    &=&
    b^{ad}_{i_0+i_1+1}(h\circ \delta_{X'}\otimes \delta_{X'}^{\otimes (i_1-1)}\otimes s\circ h'\otimes \delta_{Y'}^{\otimes i_0})\hfill\\
    &=&
    b^{ad}_{i_0+i_1+1}(h\circ\delta_{X'}\otimes \delta^{\otimes (i_1-1)}_{X'}\otimes s\otimes h'\circ \delta_{Y'}\otimes \delta_{Y'}^{\otimes (i_0-1)})\hfill\\
    &=&
    -b^{ad}_{i_0+i_1+1}(h\circ\delta_{X'}\otimes \delta^{\otimes (i_1-1)}_{X'}\otimes s\otimes \delta_Y\circ h'\otimes \delta_{Y'}^{\otimes (i_0-1)})=\cdots\hfill\\
    &=&
    -b^{ad}_{i_0+i_1+1}(h\circ\delta_{X'}\otimes \delta^{\otimes (i_1-1)}_{X'}\otimes s\otimes \delta_Y^{\otimes (i_0-1)}\otimes   \delta_Y\circ h')\hfill\\
    &=&
    -h\circ b^{ad}_{i_0+i_1+1}(\delta^{\otimes i_1}_{X'}\otimes s\otimes \delta_Y^{\otimes (i_0-1)}\otimes   \delta_Y\circ h')\hfill\\
    &=&
     -h\circ b^{ad}_{i_0+i_1+1}(\delta^{\otimes i_1}_{X'}\otimes s\otimes \delta_Y^{\otimes i_0})\circ h'\hfill\\
   \end{matrix}$$
for all $i_0,i_1\geq 0$. 
The wanted formula follows from this.
 \end{proof}

  \begin{lemma}\label{L: inverso de matriz u triangular en Z(hat(B))}
 Let $(X,\delta_X)$, $(Y,\delta_Y)$, $(X',\delta_{X'})$, $(Y',\delta_{Y'})$, $E=(X\oplus Y,\delta_E)$, and $E'=(X'\oplus Y',\delta_{E'})$ be objects in ${\cal Z}(Z)$, and  $s:Y\rightmap{}X'$  a homogeneous morphism in $\ad(Z)$ with degree $-1$. Suppose that 
   $$\delta_{E}=
 \begin{pmatrix}
 \delta_X&\gamma\\ 0&\delta_Y\\
 \end{pmatrix} \hbox{ \ and \  }\delta_{E'}=
 \begin{pmatrix}
 \delta_{X'}&\gamma'\\   0&\delta_{Y'}\\
 \end{pmatrix},$$
 where $\gamma,\gamma':Y\rightmap{}X$ are homogeneous morphisms with degree $0$.
 Assume that 
 $h=\begin{pmatrix}
    h_{1,1}&s\\
    0&h_{2,2}\\
  \end{pmatrix}:(E,\delta_E)\rightmap{}(E',\delta_{E'})$ in ${\cal Z}(Z)$, where $h_{1,1}:X\rightmap{}X'$ and $h_{2,2}:Y\rightmap{}Y'$ are special isomorphisms in $\ad(Z)$. Then, the matrix  
 $$h'=\begin{pmatrix}
    h_{1,1}^{-1}&-h_{1,1}^{-1}\circ s\circ h_{2,2}^{-1}\\
    0&h_{2,2}^{-1}\\
  \end{pmatrix}:E'\rightmap{}E $$
 is a morphism 
 $h':(E',\delta_{E'})\rightmap{}(E,\delta_E)$ in ${\cal Z}(Z)$. The morphisms $h$ and  $h'$ are mutual inverses in ${\cal Z}(Z)$.  
 \end{lemma}
 
 \begin{proof} Since $h$ is a morphism of ${\cal Z}(Z)$, by (\ref{L: isos matriciales entre terminos de enmedio}), we have the equality  
 $b_1^{tw}(s)=-\gamma'\circ h_{2,2}-h_{1,1}\circ \gamma$;  
 in order to show that the matrix $h'$ is a morphism of ${\cal Z}(Z)$, we need to show that 
 $b_1^{tw}(h^{-1}_{1,1}\circ s\circ h_{2,2}^{-1})=\gamma\circ h_{2,2}^{-1}+h_{1,1}^{-1}\circ \gamma'$. 
 Indeed, this follows from (\ref{L: b1tw(conjugado(s))=conjugado(b1tw(s))}). 
 
 Now, we have to show that $h'\star h=\hueca{I}_E$ and 
 $h\star h'=\hueca{I}_{E'}$. We only show the first one, since the other one is similar.

 We can write $h=h_0+h_1$ and $h'=h'_0+h'_1$, where
$h_0=
      \begin{pmatrix}
      h_{1,1}&0\\ 0&h_{2,2}\\
      \end{pmatrix}$,  
$$ h'_0=
      \begin{pmatrix}
      h_{1,1}^{-1}&0\\ 0&h_{2,2}^{-1}\\
      \end{pmatrix}, 
 h_1=
      \begin{pmatrix}
       0&s\\ 0&0\\
      \end{pmatrix}, \hbox{ and }
 h'_1= \begin{pmatrix}
       0&-h_{1,1}^{-1}\circ s\circ h_{2,2}^{-1}\\ 0&0\\
      \end{pmatrix}. $$

 Then, we have $h'\star h=h'_0\star h_0+h'_0\star h_1+h'_1\star h_0+h'_1\star h_1$. 
 We have $h'_0\star h_0=h'_0\circ h_0=
 \begin{pmatrix}
      \hueca{I}_X&0\\ 0&\hueca{I}_Y\\
      \end{pmatrix}$, 
      $h'_0\star h_1=h'_0\circ h_1=
      \begin{pmatrix}
       0&h_{1,1}^{-1}\circ s\\ 0&0\\
      \end{pmatrix}$, 
      $h'_1\star h_0=h'_1\circ h_0=
      \begin{pmatrix}0&-h_{1,1}^{-1}\circ s\\0&0\\ \end{pmatrix}$ and, finally, 
       $$h'_1\star h_1=b_2^{tw}(h'_1\otimes h_1)=\sum_{i_0,i_1,i_2\geq 0}b^{ad}_{i_0+i_1+i_2+2}(\delta^{\otimes i_2}_E\otimes h'_1\otimes\delta^{\otimes i_1}_{E'}\otimes h_1\otimes \delta_E^{\otimes i_0})=0.$$

Then, we have 
$h'\star h=
\begin{pmatrix} \hueca{I}_X&0\\ 0&\hueca{I}_Y\\                                                                                                               \end{pmatrix}=\hueca{I}_E$. 
 \end{proof}

 \begin{proposition}\label{P: morfismo entre conflac especiales es iso} 
Assume that the following diagram commutes  in ${\cal Z}(Z)$:
  $$\begin{matrix}
     \xi&:&(X,\delta_X)& \rightmap{ f}&(E,\delta_E)&\rightmap{g}&(Y,\delta_Y)\\
    && \shortlmapdown{\hueca{I}_X}&& \shortlmapdown{h}&&\shortlmapdown{\hueca{I}_Y}\\
     \xi'&:&(X,\delta_X)& \rightmap{ f'}&(E',\delta_{E'})&\rightmap{g'}&(Y,\delta_Y),\\
     \end{matrix}$$
 with $\xi$ and $\xi'$ special conflations.
 Then, $h$ is an isomorphism of ${\cal Z}(Z)$ and $\xi\rightmap{h\simeq}\xi'$. Moreover, we also have $\xi'\rightmap{ \ h^{-1}\simeq \ }\xi$. 
\end{proposition}

 \begin{proof} By assumption, we have $f=\sum_uf_u\otimes \frak{e}_u$ and $g=\sum_ug_u\otimes \frak{e}_u$.  Moreover, for each $u\in {\cal P}$, we have the exact sequence of vector spaces 
 $$0\rightmap{}Xe_u\rightmap{f_u}Ee_u\rightmap{g_u}Ye_u\rightmap{}0.$$
 Then, we have vector space decompositions $Ee_u=E^1_u\oplus E_u^2$. Moreover, we have    $f_u=(\tilde{f}_u,0)^t$, and $g_u=(0,\tilde{g}_u)$, with $\tilde{f}_u:Xe_u\rightmap{}E^1_u$ and $\tilde{g}_u:E^2_u\rightmap{}Ye_u$  linear isomorphisms. 
 
 Then, we have the right $S$-module decomposition $E=E^1\oplus E^2$, where  $E^1=\bigoplus_uE^1_u$ and $E^2=\bigoplus_uE^2_u$. Moreover, we have $f=(\tilde{f},0)^t$ and $g=(0,\tilde{g})$, where $\tilde{f}=\sum_u\tilde{f}_u\otimes\frak{e}_u:X\rightmap{}E^1$ and $\tilde{g}=\sum_u\tilde{g}_u\otimes\frak{e}_u:E^2\rightmap{}Y$ are special isomorphisms in $\ad(Z)$. We have a similar description for $\xi'$. Suppose that the matrix form of $h$ in $\ad(Z)$ is
 $$h=\begin{pmatrix}
   h_{1,1}&h_{1,2}\\ h_{2,1}&h_{2,2}\\  
    \end{pmatrix}.$$
Then, from the commutativity of the diagram, we obtain 
$(h_{1,1}\circ \tilde{f},h_{2,1}\circ \tilde{f})^t
= (\tilde{f'}\circ\hueca{I}_X,0)^t$ 
 and 
$(0, \hueca{I}_Y\circ\tilde{g})=(\tilde{g}'\circ h_{2,1},\tilde{g}'\circ h_{2,2}).$
From (\ref{L: compos entre tres morfs homogeneos en ad(Z)}), we obtain  $h_{2,1}=0$, $h_{1,1}\circ \tilde{f}=\tilde{f'}$, and $\tilde{g}=\tilde{g'}\circ h_{2,2}$. Therefore, we get 
$h= \begin{pmatrix}
      h_{1,1}&s\\ 0&h_{2,2}\\
     \end{pmatrix},$
 with  $h_{1,1}=\tilde{f'}\circ\tilde{f}^{-1}$,  $h_{2,2}=\tilde{g}'^{-1}\circ \tilde{g}$,  and the morphism  $s:E^2\rightmap{}E^1$ in $\ad(Z)$ is homogeneous of degree $-1$. From (\ref{R: cuando un estricto esta en Z(hat(B))}), since $f,f',g,g'$ are special morphisms, so are the components
 $\tilde{f},\tilde{g},\tilde{f}',\tilde{g}'$.    Thus, $h_{1,1}:E^1\rightmap{}E'^1$ and $h_{2,2}:E^2\rightmap{}E'^2$ are special isomorphisms in $\ad(Z)$.
     
    From (\ref{L: previo a forma canonica de conflaciones}), we get that the morphisms $\delta_E$ and $\delta_{E'}$ have triangular matrix form. Then, we can apply  (\ref{L: inverso de matriz u triangular en Z(hat(B))}), and obtain that the morphism $h$ is an isomorphism in ${\cal Z}(Z)$ with  inverse   $h':E'\rightmap{}E$ given by the matrix
    $h'=\begin{pmatrix}
    h_{1,1}^{-1}&-h_{1,1}^{-1}\circ s\circ h_{2,2}^{-1}\\
    0&h_{2,2}^{-1}\\
  \end{pmatrix}.$
  The verification of the commutativity of the following diagram  in ${\cal Z}(Z)$:
 $$\begin{matrix}
     \xi&:&(X,\delta_X)& \rightmap{ f}&(E,\delta_E)&\rightmap{g}&(Y,\delta_Y)\\
    && \shortlmapup{\hueca{I}_X}&& \shortrmapup{h^{-1}}&&\shortrmapup{\hueca{I}_Y}\\
     \xi'&:&(X,\delta_X)& \rightmap{ f'}&(E',\delta_{E'})&\rightmap{g'}&(Y,\delta_Y),\\
     \end{matrix}$$
     is straightforward. Thus $\xi'\rightmap{ \ h^{-1}\simeq \ } \xi$. 
 \end{proof}

  \begin{proposition}\label{L: equiv de canonical conflation is equiv relation} 
The relation ``$\rightmap{\simeq}$'' is an equivalence 
relation in the class of all the special conflations.
\end{proposition}

 \begin{proof} The relation ``$\rightmap{\simeq}$'' is symmetric by (\ref{P: morfismo entre conflac especiales es iso}). 
 Let us show that this relation is transitive. Consider the following diagram in ${\cal Z}(Z)$: 
 $$\begin{matrix}
     \xi&:&(X,\delta_X)& \rightmap{ f}&E&\rightmap{g}&(Y,\delta_Y)\\
    && \shortlmapdown{\hueca{I}_X}&& \shortrmapdown{h}&&\shortrmapdown{\hueca{I}_Y}\\
     \chi&:&(X,\delta_X)& \rightmap{f'}&E'&\rightmap{g'}&(Y,\delta_Y)\\
      && \shortlmapdown{\hueca{I}_X}&& \shortrmapdown{h'}&&\shortrmapdown{\hueca{I}_Y}\\
 \zeta&:&(X,\delta_X)& \rightmap{f''}&E''&\rightmap{g''}&(Y,\delta_Y),\\
   \end{matrix}$$
where the rows $\xi,\chi,\zeta$ are special conflations, every internal square is commutative, and $h$ and $h'$ are isomorphisms of ${\cal Z}(Z)$. As in the proof of (\ref{P: morfismo entre conflac especiales es iso}), we have triangular matrix expressions
$$\delta_{E}=\begin{pmatrix}
            \delta_{E^1}&\gamma\\
            0&\delta_{E^2}\\
           \end{pmatrix}, \delta_{E'}=\begin{pmatrix}
            \delta_{E'^1}&\gamma'\\
            0&\delta_{E'^2}\\
           \end{pmatrix}, \hbox{ and \ } \delta_{E''}=\begin{pmatrix}
            \delta_{E''^1}&\gamma''\\
            0&\delta_{E''^2}\\
           \end{pmatrix}.$$
Moreover, the morphisms $h$ and $h'$ have the following matrix form:
$$h= \begin{pmatrix}
      h_{1,1}&s\\ 0&h_{2,2}\\
     \end{pmatrix} \hbox{ \ and \ } h'= \begin{pmatrix}
      h'_{1,1}&s'\\ 0&h'_{2,2}\\
     \end{pmatrix},$$
     where the diagonal morphisms are special isomorphisms in ${\cal Z}(Z)$, and $s$ and $s'$ are morphisms in $\ad(Z)$. 
In order to show that $\xi$ is equivalent to $\zeta$, we will see that $h'\star h$ is an isomorphism in ${\cal Z}(Z)$ with matrix form  
$$h'\star h=\begin{pmatrix}
       h'_{1,1}\circ h_{1,1}&h'_{s,s}\circ s+s'\circ h_{2,2}\\ 0&h'_{2,2}\circ h_{2,2}\\
     \end{pmatrix}.$$
 
We can write $h=h_0+h_1$ and $h'=h'_0+h'_1$, where
$$h_0=
      \begin{pmatrix}
      h_{1,1}&0\\ 0&h_{2,2}\\
      \end{pmatrix}, 
 h'_0=
      \begin{pmatrix}
      h'_{1,1}&0\\ 0&h'_{2,2}\\
      \end{pmatrix}, 
 h_1=
      \begin{pmatrix}
       0&s\\ 0&0\\
      \end{pmatrix}, \hbox{ and }
 h'_1= \begin{pmatrix}
       0&s'\\ 0&0\\
      \end{pmatrix}. $$

 Then, we have $h'\star h=h'_0\star h_0+h'_0\star h_1+h'_1\star h_0+h'_1\star h_1$. 
 We have $h'_0\star h_0=h'_0\circ h_0= \begin{pmatrix}
      h'_{1,1}\circ h_{1,1}&0\\ 0&h'_{2,2}\circ h_{2,2}\\
      \end{pmatrix}$, 
      $h'_0\star h_1=h'_0\circ h_1=
      \begin{pmatrix}
       0&h'_{1,1}\circ s\\ 0&0\\
      \end{pmatrix}$, 
      $h'_1\star h_0=h'_1\circ h_0=\begin{pmatrix}0&s'\circ h_{2,2}\\0&0\\ \end{pmatrix}$ 
      and, finally, 
       $$h'_1\star h_1=b_2^{tw}(h'_1\otimes h_1)=\sum_{i_0,i_1,i_2\geq 0}b^{ad}_{i_0+i_1+i_2+2}(\delta^{\otimes i_2}_{E''}\otimes h'_1\otimes\delta^{\otimes i_2}_{E'}\otimes h_1\otimes \delta_E^{\otimes i_0})=0.$$
       So we get the wanted triangular matrix form for $h'\star h$.  This implies that the squares in the following diagram commute in $\ad(Z)$ (and in ${\cal Z}(Z)$):
        $$\begin{matrix}
     \xi&:&(X,\delta_X)& \rightmap{f}&E&\rightmap{g}&(Y,\delta_Y)\\
    && \shortlmapdown{\hueca{I}_X}&& \shortrmapdown{h'\star h}&&\shortrmapdown{\hueca{I}_Y}\\
     \zeta&:&(X,\delta_X)& \rightmap{f''}&E''&\rightmap{g''}&(Y,\delta_Y).\\
     \end{matrix}$$
     Indeed, this commutativity follows from the description of $h_{1,1},h_{2,2},h'_{1,1},h'_{2,2}$ in terms of the components of $f,f',f'',g,g',g''$ given in the proof of (\ref{P: morfismo entre conflac especiales es iso}). 
Again, from (\ref{P: morfismo entre conflac especiales es iso}), we know that $h'\star h$ is an isomorphism in ${\cal Z}(Z)$. So, the relation ``$\rightmap{\simeq}$'' is transitive in the class of special conflations. 
 \end{proof}
 
  \begin{lemma}\label{L: conlfacion es par exacto}
 Every special conflation 
 $\xi:\hbox{ \ }(X,\delta_X)\rightmap{f}(E,\delta_E)\rightmap{g}(Y,\delta_Y)$
is an exact pair in ${\cal Z}(Z)$. That is $f=\Ker g$ and $g=\Coker f$ in ${\cal Z}(Z)$.
 \end{lemma}
 
 \begin{proof} Because of (\ref{L: may assume special conflation direct sum middle term}), we may assume that $\xi$ is a canonical conflation. Thus,
 $E=X\oplus Y$, $f=(\hueca{I}_X,0)^t$, $g=(0,\hueca{I}_Y)$, and 
 $\delta_{E}=
 \begin{pmatrix}
 \delta_X&\gamma\\ 0&\delta_Y\\
 \end{pmatrix}$, for some homogeneous morphism 
 $\gamma:Y\rightmap{}X$ in $\ad(Z)$ of degree $0$.
 
 Let $h:(W,\delta_W)\rightmap{}(X,\delta_X)$ be a morphism in ${\cal Z}(Z)$ such that $f\star h=0$. Since $f$ is special, we have $0=f\circ h=(\hueca{I}_X,0)^t\circ h=(\hueca{I}_X\circ h,0)^t$ and, hence,  $0=\hueca{I}_X\circ h=h$, so $f$ is a monomorphism in ${\cal Z}(Z)$. Similarly, $g$ is an epimorphism in ${\cal Z}(Z)$. 
 
 Assume now that $h=(h_1,h_2):(E,\delta_E)\rightmap{}(W,\delta_W)$ is a morphism in ${\cal Z}(Z)$ such that $h\star f=0$. Again, we have $0=h\circ f=(h_1,h_2) \circ(\hueca{I}_X,0)^t= h_1\circ\hueca{I}_X$. Then, $h_1=0$ and we have the morphism $h_2:Y\rightmap{}W$ such that 
 $h_2\circ g=h_2\circ (0,\hueca{I}_Y)=(0,h_2\circ \hueca{I}_Y)=(h_1,h_2)=h.$ 
 
 By assumption, 
 $0=b_1^{tw}(h)=\sum_{i_0,i_1\geq 0}b^{ad}_{i_0+i_1+1}(\delta_W^{\otimes i_1}\otimes h\otimes \delta_E^{\otimes i_0}).$ 
 Then,  we have 
 $0=b_1^{tw}(h)=(0, \sum_{i_0,i_1\geq 0}b^{ad}_{i_0+i_1+1}(\delta_W^{\otimes i_1}\otimes h_2\otimes\delta_Y^{\otimes i_0}))=(0,b^{tw}_1(h_2)),$
 and $h_2$ is a morphism in ${\cal Z}(Z)$. It satisfies $h=h_2\circ g=h_2\star g$. So, we have that $g$ is the cokernel of $f$ in ${\cal Z}(Z)$. 
 
 The fact that $f$ is the kernel of $g$ is proved dually. 
 \end{proof}

\begin{lemma}\label{L: pushouts para conflaciones en Z(hat(B))}
Let   
$\xi:\hbox{ \ }(X,\delta_X)\rightmap{f}(E,\delta_E)\rightmap{g}(Y,\delta_Y)$
be a canonical conflation in ${\cal Z}(Z)$ with 
$\delta_E=
\begin{pmatrix}\delta_X&\gamma\\ 0&\delta_Y
\end{pmatrix}$ 
and $h:(X,\delta_X)\rightmap{}(X_1,\delta_{X_1})$ any morphism in ${\cal Z}(Z)$. Then, we have the following commutative diagram in ${\cal Z}(Z)$
 $$\begin{matrix}
   (X,\delta_X)&\rightmap{f}&(E,\delta_E)&\rightmap{g}&(Y,\delta_Y)\hfill\\
   \shortlmapdown{h}&&\shortrmapdown{t}&&\shortrmapdown{\hueca{I}_Y}\\
   (X_1,\delta_{X_1})&\rightmap{f_1}&(E_1,\delta_{E_1})&\rightmap{g_1}&(Y,\delta_Y),\hfill\\
   \end{matrix}$$
where $t=\begin{pmatrix}
         h&0\\ 0&\hueca{I}_Y\\ 
         \end{pmatrix}$
and the second row is a canonical conflation with 
$$ \delta_{E_1}=
\begin{pmatrix}
\delta_{X_1}&\gamma_1\\ 0&\delta_Y\\
\end{pmatrix} \hbox{ and \ } \gamma_1=h\star\gamma.$$
\end{lemma}

\begin{proof} By assumption,  we have    
 $E=X\oplus Y$, $f=(\hueca{I}_X,0)^t$, $g=(0,\hueca{I}_Y)$, and 
  $\gamma:Y\rightmap{}X$ is a homogeneous morphism 
  in $\ad(Z)$ of degree $0$.
  From (\ref{L: diferenciales de la suma directa en tw(hat(B))}), we know that $\gamma$ satisfies $b_1^{tw}(\gamma)=0$. By assumption $h$ has degree $-1$ and satisfies $b_1^{tw}(h)=0$. Then, 
  by (\ref{R: notacion star}),  the composition 
 $$\gamma_1:=h\star \gamma=b_2^{tw}(h\otimes \gamma):(Y,\delta_Y)\rightmap{}(X_1,\delta_{X_1})$$
 satisfies $b_1^{tw}(\gamma_1)=0$ and is homogeneous of degree $0$.
  Therefore, 
 by (\ref{L: diferenciales de la suma directa en tw(hat(B))}), we have the following 
 object of ${\cal Z}(Z)$:
 $$(E_1,\delta_{E_1}),\hbox{ \ where \ } E_1=X_1\oplus Y \hbox{ \ and \ }
 \delta_{E_1}=\begin{pmatrix}
 \delta_{X_1}& \gamma_1\\ 0& \delta_Y
 \end{pmatrix}.$$
 From (\ref{L: previo a forma canonica de conflaciones}), we obtain the canonical conflation                                                              
 $$\xi_1:(X_1,\delta_{X_1})\rightmap{f_1}(E_1,\delta_{E_1})\rightmap{g}(Y,\delta_Y),$$ 
           where  $f_1=(\hueca{I}_{X_1},0)^t$ and $g_1=(0,\hueca{I}_Y)$. 
  Consider the homogeneous morphism of degree $-1$ in $\tw(Z)$ given by the matrix 
  $t=\begin{pmatrix}
       h&0\\ 0&\hueca{I}_Y
      \end{pmatrix}:(E,\delta_E)\rightmap{}(E_1,\delta_{E_1}).$

      In order to show that $b_1^{tw}(t)=0$,  
      consider the morphisms $t_0:=\begin{pmatrix}
       0&0\\ 0&\hueca{I}_Y
      \end{pmatrix}$ 
      and
      $t_1:=\begin{pmatrix}
       h&0\\ 0&0
      \end{pmatrix}$, from $(E,\delta_E)$ to $(E_1,\delta_{E_1})$ is $\tw(Z)$. So, we get 
       $b_1^{tw}(t)= b_1^{tw}(t_0)+ b_1^{tw}(t_1)$.
      Since $t_0$ is strict, we have 
      $b_1^{tw}(t_0)=\delta_{E_1}\circ t_0+t_0\circ \delta_E$. Hence, 
      $$b_1^{tw}(t_0)=\begin{pmatrix}
                        0&\gamma_1\circ \hueca{I}_Y\\ 0&\delta_Y\circ \hueca{I}_Y\\
                       \end{pmatrix}+
                       \begin{pmatrix}
                        0&0\\ 0&\hueca{I}_Y\circ \delta_Y\\
                       \end{pmatrix}
                       =
                       \begin{pmatrix}
                        0&-\gamma_1\\ 0&0\\
                       \end{pmatrix}.$$

   Moreover, we have 
   $$\begin{matrix}
   b_1^{tw}(t_1)&=& \sum_{i_0,i_1\geq 0}b^{ad}_{i_0+i_1+1}(\delta_{E_1}^{\otimes i_1}\otimes t_1\otimes \delta_E^{\otimes i_0})\hfill\\
   &\,&\hfill\\
   &=&
   \begin{pmatrix}
    b_1^{tw}(h)
    & b_2^{tw}(h\otimes\gamma) \\
    0&0\\
   \end{pmatrix}=
   \begin{pmatrix}
    0&h\star \gamma\\
    0&0
   \end{pmatrix}.\hfill\\
     \end{matrix}$$
Therefore, we get 
$b_1^{tw}(t)=
\begin{pmatrix}
0&h\star \gamma- \gamma_1\\
0&0\\
\end{pmatrix}=0$, 
so $t:(E,\delta_E)\rightmap{}(E_1,\delta_{E_1})$ is a morphism in 
${\cal Z}(Z)$, as we wanted to show. The diagram commutes, because it commutes with respect to $\circ$ and all the implicit compositions involve the composition with a strict morphism, see (\ref{R: star-compos=circ-compos para special morfisms}).  
\end{proof}

Similarly, we have the following  statement.

 \begin{lemma}\label{L: pullbacks para conflaciones en Z(hat(B))}
Let   
$\xi:\hbox{ \ }(X,\delta_X)\rightmap{f}(E,\delta_E)\rightmap{g}(Y,\delta_Y)$
be a canonical conflation in ${\cal Z}(Z)$ with 
$\delta_E=
\begin{pmatrix}
\delta_X&\gamma\\ 0&\delta_Y\\           
\end{pmatrix}$ 
and $h:(Y_1,\delta_{Y_1})\rightmap{}(Y,\delta_Y)$ any morphism in ${\cal Z}(Z)$. Then, we have the following commutative diagram in ${\cal Z}(Z)$
 $$\begin{matrix}
  (X,\delta_X)&\rightmap{f_1}&(E_1,\delta_{E_1})&\rightmap{g_1}&(Y_1,\delta_{Y_1}).\hfill\\
   \shortlmapdown{\hueca{I}_X}&&\shortrmapdown{t}&&\shortrmapdown{h}\\
    (X,\delta_X)&\rightmap{f}&(E,\delta_E)&\rightmap{g}&(Y,\delta_Y),\hfill\\
   \end{matrix}$$
where $t=\begin{pmatrix}
         \hueca{I}_X&0\\ 0&h\\
        \end{pmatrix}$ and the first row is a canonical conflation with 
$$\delta_{E_1}=
\begin{pmatrix}
\delta_X&\gamma_1\\ 0&\delta_{Y_1}\\
\end{pmatrix} \hbox{ and \ } \gamma_1=-\gamma\star h.$$
\end{lemma}

\begin{proof} Similar to the proof of (\ref{L: pushouts para conflaciones en Z(hat(B))}).                         
\end{proof}

\begin{lemma}\label{L: Z(hat(B)) is an additive precategory}
The precategory ${\cal Z}(Z)$ has zero object $0=(0,0)$ and finite biproducts described as follows:
Given any finite family $(X_1,\delta_{X_1}),\ldots,(X_n,\delta_{X_n})$ of objects in ${\cal Z}(Z)$, we have the object $(X,\delta_X)$
in ${\cal Z}(Z)$, with $X=\bigoplus_{i=1}^nX_i$ and $\delta_X:X\rightmap{}X$  is the morphism  in $\ad(X)$ with diagonal matrix form with components $\delta_{X_1},\ldots,\delta_{X_n}$. We have  special morphisms
$s_{X_j}:(X_j,\delta_{X_j})\rightmap{}(X,\delta_{X})$ and $p_{X_j}:(X,\delta_{X})\rightmap{}(X_j,\delta_{X_j})$ in ${\cal Z}(Z)$, defined by  the morphisms $s_{X_j}:X_j\rightmap{}X$ and $p_{X_j}:X\rightmap{}X_j$ in $\ad(Z)$ introduced in (\ref{R: componentes de morfismos homogeneos en ad(Z)}).
They satisfy the relations:
$p_{X_j}\star s_{X_j}=id_{(X_j,\delta_{X_j})}$,  for all $j$,
    $p_{X_j}\star s_{X_i}=0$,  for all $ i\not=j$, and
 $id_{(X,\delta_{X})}=\sum_{i=1}^n s_{X_i}\star p_{X_i}$.

From now on, we use the notation  $\bigoplus_{i=1}^n(X_i,\delta_{X_i}):=(X,\delta_X)$.
As in any additive category,  each morphism
$f:\bigoplus_{i=1}^n(X_i,\delta_{X_i})\rightmap{}\bigoplus_{j=1}^m(Y_j,\delta_{Y_j})$ in ${\cal Z}(Z)$ is determined by its matrix $M(f):=(f_{j,i})$, where $f_{j,i}=p_{Y_j}\star f\star s_{X_i}$, and $f$ can be recovered from its matrix with the formula  $f=\sum_{i,j} s_{Y_j}\star f_{j,i}\star p_{X_i}$.
As usual, we will identify each morphism $f:\bigoplus_{i=1}^n(X_i,\delta_{X_i})\rightmap{}\bigoplus_{j=1}^m(Y_j,\delta_{Y_j})$ of ${\cal Z}(Z)$ with its matrix $M(f)$. When we forget the second components of the objects in ${\cal Z}(Z)$, the matrix notation for  the morphism $f:X\rightmap{}Y$ of $\ad(Z)$ of (\ref{R: componentes de morfismos homogeneos en ad(Z)}), coincides with the one mentioned here.
\end{lemma}

\begin{proof} It is easy to see that indeed $(X,\delta_X)\in {\cal Z}(Z)$. Using (\ref{R: cuando un estricto esta en Z(hat(B))}), the remaining  verifications are straightforward, see (\ref{R: componentes de morfismos homogeneos en ad(Z)}) and (\ref{R: descomposicon de morfismos en ad(hat(B))}).
\end{proof}

\begin{proposition}\label{P: list propiedades conflacionarias especiales} The class of special conflations in the additive precategory ${\cal Z}(Z)$ has the following properties:
\begin{enumerate}
 \item If 
 $(X,\delta_X)\rightmap{f}(E,\delta_E)\rightmap{g}(Y,\delta_Y)$ is a special conflation, then $f$ is kernel of $g$ and $g$ is cokernel of $f$ in the precategory  ${\cal Z}(Z)$.
 \item Composition of special inflations is a special inflation and composition of special deflations is a special deflation. 
 \item For each special inflation $f:(X,\delta_X)\rightmap{}(E,\delta_E)$ and each morphism $h:(X,\delta_X)\rightmap{}(X',\delta_{X'})$ there are a special inflation $f':(X',\delta_{X'})\rightmap{}(E',\delta_{E'})$ and a morphism $h':(E,\delta_E)\rightmap{}(E',\delta_{E'})$ such that $h'\star f=f'\star h$.  
 \item For each special deflation $g:(E,\delta_E)\rightmap{}(Y,\delta_Y)$ and each morphism $h:(Y',\delta_{Y'})\rightmap{}(Y,\delta_Y)$ there are a special deflation $g':(E',\delta_{E'})\rightmap{}(Y',\delta_{Y'})$ and a morphism $h':(E',\delta_{E'})\rightmap{}(E,\delta_E)$ such that $g\star h'=h\star g'$.  
 \item Identity morphisms are special inflations and special deflations. Moreover, if $f$ and $g$ are special composable morphisms and $g\star f$ is a special inflation (resp. a special deflation), then $f$  is a special inflation (resp. $g$ is a special deflation).
\end{enumerate}
\end{proposition}

\begin{proof} The additivity of ${\cal Z}(Z)$ was remarked in (\ref{L: Z(hat(B)) is an additive precategory}); (1) is just (\ref{L: conlfacion es par exacto}); (2) and (5) follow from (\ref{L: caracterizacion de special infl y defl}); (3) follows from (\ref{L: may assume special conflation direct sum middle term}) and (\ref{L: pushouts para conflaciones en Z(hat(B))}); (4) follows from (\ref{L: may assume special conflation direct sum middle term}) and (\ref{L: pullbacks para conflaciones en Z(hat(B))}).
\end{proof}

\begin{remark}\label{R: comentario al summary for special conflations}
The last summary shows that, although ${\cal Z}(Z)$ is not a category, it is an additive precategory and the special conflations satisfy  properties which are similar to those of conflations of exact structures in additive categories. 
\end{remark}

We close this section with a couple of remarks on split special conflations.
We say that a special conflation $(X,\delta_X)\rightmap{f}(E,\delta_E)\rightmap{g}(Y,\delta_Y)$ in ${\cal Z}(Z)$ \emph{splits} iff there are morphisms $f':(E,\delta_E)\rightmap{}(X,\delta_X)$ and $g':(Y,\delta_Y)\rightmap{}(E,\delta_E)$ in ${\cal Z}(Z)$ such that $f'\star f=\hueca{I}_X$  and $g\star g'=\hueca{I}_Y$. This is the case of \emph{the trivial ones}
$$\xi_0:\hbox{ \ }(X,\delta_X)\rightmap{ \ (\hueca{I}_X,0)^t \ }(E,\delta_E)\rightmap{ \ (0,\hueca{I}_Y) \ }(Y,\delta_Y),$$
where $(E,\delta_E)=(X,\delta_X)\oplus(Y,\delta_Y)$.
Indeed, from (\ref{R: cuando un estricto esta en Z(hat(B))}), we get that the special morphisms $(\hueca{I}_X,0):(E,\delta_{E})\rightmap{}(X,\delta_X)$ and $(0,\hueca{I}_Y)^t:(Y,\delta_Y)\rightmap{}(E,\delta_{E})$ belong to ${\cal Z}(Z)$, and they clearly provide a splitting of the special conflation
$\xi_0$.

\begin{lemma}\label{L: xi equiv a trivial implica xi splits}
 For any special conflation $\xi:(X,\delta_X)\rightmap{f}(E,\delta_E)\rightmap{g}(Y,\delta_Y)$ in ${\cal Z}(Z)$ the following statements are equivalent:
 \begin{enumerate}
  \item The special conflation $\xi$ splits;
  \item There is a morphism $f':(E,\delta_E)\rightmap{}(X,\delta_X)$  with $f'\star f=\hueca{I}_X$;
  \item There is a morphism $g':(Y,\delta_Y)\rightmap{}(E,\delta_E)$  with $g\star g'=\hueca{I}_Y$;
  \item The special conflation $\xi$ is equivalent to a trivial one.
 \end{enumerate}
\end{lemma}

\begin{proof} If $\xi\simeq \xi_0$, where $\xi_0$ is a trivial conflation, we have a commutative diagram
$$\begin{matrix}
     \xi&:&(X,\delta_X)& \rightmap{ \ f \ }&(E,\delta_E)&\rightmap{ \ g \ }&(Y,\delta_Y)\\
    && \shortlmapdown{\hueca{I}_X}&& \shortlmapdown{h}&&\shortlmapdown{\hueca{I}_Y}\\
     \xi_0&:&(X,\delta_X)& \rightmap{ \ s \ }&(\overline{E},\delta_{\overline{E}})&\rightmap{ \ p \ }&(Y,\delta_Y) \\
     \end{matrix}$$
 in ${\cal Z}(Z)$. We know that there is $s':(\overline{E},\delta_{\overline{E}})\rightmap{}(X,\delta_X)$ such that $s'\star s=\hueca{I}_X$. Hence $f':=s'\star h$ satisfies
$f'\star f=f'\circ f=s'\circ h\circ f=s'\circ s=\hueca{I}_X$.
The relation ``$\simeq$'' is symmetric, so $\xi_0\simeq \xi$, and we have a commutative diagram in ${\cal Z}(Z)$
$$\begin{matrix}
     \xi_0&:&(X,\delta_X)& \rightmap{ \ s \ }&(\overline{E},\delta_{\overline{E}})&\rightmap{ \ p \ }&(Y,\delta_Y)\\
    && \shortlmapdown{\hueca{I}_X}&& \shortlmapdown{h'}&&\shortlmapdown{\hueca{I}_Y}\\
    \xi&:&(X,\delta_X)& \rightmap{ \ f \ }&(E,\delta_E)&\rightmap{ \ g \ }&(Y,\delta_Y).\\
     \end{matrix}$$
 Now, we consider a morphism $p':(Y,\delta_Y)\rightmap{}(\overline{E},\delta_{\overline{E}})$ in ${\cal Z}(Z)$ such that $p\star p'=\hueca{I}_Y$ and notice that $g':=h'\star p'$ satisfies $g\star g'=\hueca{I}_Y$. So $\xi$ splits, and 4 implies 1.

 Now, assuming $2$, we get the commutative diagram in ${\cal Z}(Z)$
$$\begin{matrix}
     \xi&:&(X,\delta_X)& \rightmap{ \ f \ }&(E,\delta_E)&\rightmap{ \ g \ }&(Y,\delta_Y)\\
    && \shortrmapdown{\hueca{I}_X}&& \shortrmapdown{(f',g)^t}&&\shortrmapdown{\hueca{I}_Y}\\
     \xi_0&:&(X,\delta_X)& \rightmap{ \  (\hueca{I}_X,0)^t \ }&(X,\delta_X)\oplus(Y,\delta_Y)&\rightmap{ \ (0,\hueca{I}_Y) \ }&(Y,\delta_Y),\\
     \end{matrix}$$
     where $(f',g)^t$ is an isomorphism by (\ref{P: morfismo entre conflac especiales es iso}). So, $\xi\simeq \xi_0$, and 2 implies 4.
 The proof of 3 implies 4 is similar.
\end{proof}

\begin{remark}\label{L: sobre conflaciones triviales}
  Notice that a canonical conflation $(X,\delta_X)\rightmap{f}(E,\delta_E)\rightmap{g}(Y,\delta_Y)$ splits iff there are
   morphisms of the form $f'=(\hueca{I}_X , s) : (E, \delta_E ) \rightmap{} (X, \delta_X )$ and $ g'=(r,\hueca{I}_Y)^t:
(Y, \delta_Y) \rightmap{}(E,\delta_E)$ in ${\cal Z}(Z)$ such that  $f'\star f = \hueca{I}_X$ and  $g\star g' = \hueca{I}_Y$.
 \end{remark}

\section{Conflations in ${\cal Z}(Z)$}\label{Confl in Z(Z)}

In his short section, we keep the notation of the preceding one and continue the study of special conflations in ${\cal Z}(Z)$. The following statements will be applied later in section \ref{Z(hatZ)}.

\begin{lemma}\label{L: isos triangulares y con esquina estricta}
 Assume that we have objects $(X,\delta_X)$, $(X,\delta'_X)$, $(Y,\delta_Y)$ 
 in ${\cal Z}(Z)$. Suppose that $\gamma:Y\rightmap{}X$ is a strict homogeneous morphism in $\ad(Z)$ with degree $0$ and $b_1^{tw}(\gamma)=0$. Then, we can consider the objects 
 $(E,\delta_E)$ and $(E',\delta_{E'})$ of ${\cal Z}(Z)$ such that $E=X\oplus Y=E'$ and 
  $$\delta_E=\begin{pmatrix}
             \delta_X&\gamma\\ 0&\delta_Y\\
             \end{pmatrix}
             \hbox {\ and \ } 
             \delta_{E'}=\begin{pmatrix}
             \delta'_X&\gamma\\ 0&\delta_Y\\
             \end{pmatrix}.$$
For any homogeneous morphism $\rho:X\rightmap{}Y$  in $\ad(Z)$ with degree $-1$ such that $\rho\circ \gamma=0$, $\gamma\circ \rho=\delta'_X-\delta_X$, and the morphisms 
$\rho:(X,\delta_X)\rightmap{}(Y,\delta_Y)\hbox{ and } \rho:(X,\delta'_X)\rightmap{}(Y,\delta_Y)$
belong to ${\cal Z}(Z)$, 
we have an isomorphism in  ${\cal Z}(Z)$: 
$$h=\begin{pmatrix}
     \hueca{I}_X&0\\ \rho&\hueca{I}_Y\\
    \end{pmatrix}:(E,\delta_E)\rightmap{}(E',\delta_{E'}).$$
    Its inverse is given by
 $h'=\begin{pmatrix}
     \hueca{I}_X&0\\ -\rho&\hueca{I}_Y\\
    \end{pmatrix}:(E',\delta_{E'})\rightmap{}(E,\delta_E).$   
\end{lemma}

\begin{proof} In order to show that $h:(E,\delta_E)\rightmap{}(E',\delta_{E'})$ is a morphism in ${\cal Z}(Z)$ using (\ref{R: desarmado f en estricto + otro simplifica btw's})(1), we define  
$h^0:= \begin{pmatrix}
        \hueca{I}_X&0\\
        0&\hueca{I}_Y\\
       \end{pmatrix}$, 
 $h^1:=\begin{pmatrix}
        0&0\\
        \rho&0\\
       \end{pmatrix}$, 
 $\delta^0_E=\begin{pmatrix}
             0&\gamma\\ 0&0\\
             \end{pmatrix}$,
 $\delta^1_E=\begin{pmatrix}
             \delta_X&0\\ 0&\delta_Y\\
             \end{pmatrix}$, 
 $\delta^0_{E'}=\begin{pmatrix}
             0&\gamma\\ 0&0\\
             \end{pmatrix}$,  
  and  
 $\delta^1_{E'}=\begin{pmatrix}
             \delta'_X&0\\ 0&\delta_Y\\
             \end{pmatrix}$.
Then, we obtain that $b_1^{tw}(h)=h\circ \delta_E+\delta_{E'}\circ h+ R(h)$, where $$R(h)=b_1^{ad}(h^1)+\sum_{\scriptsize\begin{matrix}i_0,i_1\geq 0\\ 
                          i_0+i_1\geq 2\end{matrix}} b^{ad}_{i_0+i_1+1}((\delta_{E'}^1)^{\otimes i_1}\otimes h^1\otimes (\delta_E^1)^{\otimes i_0}).$$ 
Thus, $R(h)=\begin{pmatrix}
             0&0\\ b_1^{ad}(\rho)+\sum_{\scriptsize\begin{matrix}i_0,i_1\geq 0\\ 
                          i_0+i_1\geq 2\end{matrix}} b^{ad}_{i_0+i_1+1}(\delta_Y^{\otimes i_1}\otimes \rho\otimes \delta_X^{\otimes i_0})&0\\
            \end{pmatrix}.$ Moreover, we have 

$$h\circ \delta_E+\delta_{E'}\circ h=
\begin{pmatrix}
 \hueca{I}_X\circ \delta_X+\delta'_X\circ \hueca{I}_X+\gamma\circ \rho&\hueca{I}_X\circ \gamma+\gamma\circ \hueca{I}_Y\\
 \rho\circ \delta_X+\delta_Y\circ \rho&\rho\circ\gamma+\hueca{I}_Y\circ \delta_Y+\delta_Y\circ \hueca{I}_Y\\
\end{pmatrix}.$$
Thus, $h\circ \delta_E+\delta_{E'}\circ h=
\begin{pmatrix}
 \delta_X-\delta'_X+\gamma\circ \rho&0\\
 \rho\circ \delta_X+\delta_Y\circ \rho&\rho\circ\gamma\\
\end{pmatrix}.$ It follows that $b_1^{tw}(h)=0$ iff 
$\rho\circ \gamma=0$, $\gamma\circ \rho=\delta'_X-\delta_X$, and $\rho:(X,\delta_X)\rightmap{}(Y,\delta_Y)$ is a morphism in ${\cal Z}(Z)$. 
 
By the symmetry of the assumptions of the lemma, we also have
that $h':(E',\delta_{E'})\rightmap{}(E,\delta_E)$ is a morphism in 
${\cal Z}(Z)$. 
\medskip

It remains to show that $h$ and $h'$ are mutual inverses in ${\cal Z}(Z)$. We only show that $h\star h'=id_{(E',\delta_{E'})}$, since the verification of the other equality $h'\star h=id_{(E,\delta_E)}$ is similar. In order to apply, (\ref{R: desarmado f en estricto + otro simplifica btw's})(2),    
we consider also the following  morphisms $(h')^0:=\begin{pmatrix}
                        \hueca{I}_X&0\\ 0&\hueca{I}_Y\\ 
                        \end{pmatrix}$   and   $(h')^1:=
                    \begin{pmatrix}
                        0&0\\ -\rho&0\\ 
                        \end{pmatrix}$ in $\ad(Z)$. 
  Then, we have $h\star h'=h\circ h'+R(h,h')$, where 
  $$R(h,h')=\sum_{\scriptsize\begin{matrix}i_0,i_1,i_2\geq 0\\ 
                          i_0+i_1+i_2\geq 1\end{matrix}} b^{ad}_{i_0+i_1+i_2+2}((\delta_{E'}^1)^{\otimes i_2}\otimes h^1\otimes  (\delta_E^1)^{\otimes i_1}\otimes (h')^1\otimes (\delta_{E'}^1)^{\otimes i_0}).$$
     Since every tensor factor $h^1\otimes  (\delta_E^1)^{\otimes i_1}\otimes (h')^1$ is zero, we obtain $R(h,h')=0$, so 
     $$h\star h'=h\circ h'=\begin{pmatrix}
                                        \hueca{I}_X\circ\hueca{I}_X&0\\
                                        \rho\circ \hueca{I}_X-\hueca{I}_Y\circ\rho&\hueca{I}_Y\circ\hueca{I}_Y\\
                                       \end{pmatrix}=\hueca{I}_{E'},$$ 
                                       as we wanted to show. 
\end{proof}

Similarly, we have the following.

\begin{lemma}\label{L: isos triangulares sup entre objetos con esquina estricta}
 Assume that we have objects $(X,\delta_X)$, $(Y,\delta_Y)$, $(Y,\delta'_Y)$ 
 in ${\cal Z}(Z)$. Suppose that $\gamma:Y\rightmap{}X$ is a strict homogeneous morphism in $\ad(Z)$ with degree $0$ and $b_1^{tw}(\gamma)=0$. Then, we can consider the objects 
 $(E,\delta_E)$ and $(E',\delta_{E'})$ of ${\cal Z}(Z)$ such that $E=X\oplus Y=E'$ and 
  $$\delta_E=\begin{pmatrix}
             \delta_X&\gamma\\ 0&\delta_Y\\
             \end{pmatrix}\hbox {\ and \ } 
             \delta_{E'}=\begin{pmatrix}
             \delta_X&\gamma\\ 0&\delta'_Y\\
             \end{pmatrix}.$$
For any homogeneous morphism $\rho:X\rightmap{}Y$  in $\ad(Z)$ with degree $-1$ such that $\gamma\circ \rho=0$, $\rho\circ \gamma=\delta'_Y-\delta_Y$, and the morphisms 
$\rho:(X,\delta_X)\rightmap{}(Y,\delta_Y)$ and $\rho:(X,\delta_X)\rightmap{}(Y,\delta'_Y)$
belong to ${\cal Z}(Z)$, 
we have an isomorphism in  ${\cal Z}(Z)$: 
$$h=\begin{pmatrix}
     \hueca{I}_X&0\\ \rho&\hueca{I}_Y\\
    \end{pmatrix}:(E,\delta_E)\rightmap{}(E',\delta_{E'}).$$
    Its inverse is given by
 $h'=\begin{pmatrix}
     \hueca{I}_X&0\\ -\rho&\hueca{I}_Y\\
    \end{pmatrix}:(E',\delta_{E'})\rightmap{}(E,\delta_E).$   
\end{lemma}

\begin{proof} This is similar to the proof of (\ref{L: isos triangulares y con esquina estricta}).  
\end{proof}

\begin{lemma}\label{L: caractizacion de pares en Z(hat(B)) con izquierda variada} 
Assume that $(X,\delta_X)$ is an object in $\tw(Z)$ and 
 consider a sequence of morphisms in $\tw(Z)$ of the form 
 $$(X,\delta'_X)\rightmap{f}(E,\delta_E)\rightmap{g}(Y,\delta_Y),$$
with $E=X\oplus Y$, 
$\delta_E=\begin{pmatrix}
           \delta_X&\gamma\\ 0&\delta_Y\\
           \end{pmatrix}$, $f=(\hueca{I}_X,-\rho)^t$,  $g=(\rho,\hueca{I}_Y)$, 
           where $\rho:X\rightmap{}Y$ is a morphism  in $\ad(Z)$ with degree $-1$ 
           and  $\gamma:Y\rightmap{}X$ is a strict morphism in $\ad(Z)$ with degree $0$.  Then, the  sequence lies in ${\cal Z}(Z)$ iff 
     \begin{enumerate}
     \item $\rho:(X,\delta'_X)\rightmap{}(Y,\delta_Y)$  is a morphism in ${\cal Z}(Z)$ and $\gamma\circ\rho=\delta'_X-\delta_{X}$. 
     \item $\rho:(X,\delta_X)\rightmap{}(Y,\delta_Y)$ is a morphism in ${\cal Z}(Z)$ and $ \rho\circ \gamma=0$. 
     \end{enumerate}
\end{lemma}

\begin{proof} In order to apply (\ref{R: desarmado f en estricto + otro simplifica btw's}) to the computation of $b_1^{tw}(f)$ and  $b_1^{tw}(g)$, we consider $f^0=(\hueca{I}_X,0)^t$, $f^1=(0,-\rho)^t$, $g^0=(0,\hueca{I}_Y)$, $g^1=(\rho,0)$, $\delta_E^0=\begin{pmatrix}                                                                                                                             0&\gamma\\ 0&0\\   \end{pmatrix}$ and 
$\delta_E^1=\begin{pmatrix}                                                                                                                             \delta_X&0\\ 0&\delta_Y\\   \end{pmatrix}$. 
Then, we get 
$b_1^{tw}(f)=f\circ \delta'_X+\delta_E\circ f+R(f)$,  where 
  $R(f)=b_1^{ad}(f^1)+\sum_{\scriptsize\begin{matrix}i_0,i_1\geq 0\\ 
                          i_0+i_1\geq 2\end{matrix}}b^{ad}_{i_0+i_1+1}((\delta_E^1)^{\otimes i_1}\otimes f^1\otimes (\delta'_X)^{\otimes i_0})$. 
We also have  
$$f\circ \delta'_X+\delta_E\circ f=\begin{pmatrix}
                                   \hueca{I}_X\circ \delta'_X+\delta_X\circ \hueca{I}_X-\gamma\circ \rho\\
                                   -\rho\circ \delta'_X-\delta_Y\circ \rho\\
                                  \end{pmatrix}$$
and
$$R(f)=\begin{pmatrix} 0\\
       -b_1^{ad}(\rho)-\sum_{\scriptsize\begin{matrix}i_0,i_1\geq 0\\ 
                          i_0+i_1\geq 2\end{matrix}} b^{ad}_{i_0+i_1+1}(\delta_Y^{\otimes i_1}\otimes\rho\otimes (\delta'_X)^{\otimes i_0}) \\
       \end{pmatrix}.$$
 Then, we obtain that $b_1^{tw}(f)=0$ if and only if $\gamma\circ \rho=\delta'_X-\delta_X$ and    $\rho\circ \delta'_X+\delta_Y\circ \rho+ b_1^{ad}(\rho)+ \sum_{\scriptsize\begin{matrix}i_0,i_1\geq 0\\ 
                          i_0+i_1\geq 2\end{matrix}} b^{ad}_{i_0+i_1+1}(\delta_Y^{\otimes i_1}\otimes\rho\otimes (\delta'_X)^{\otimes i_0})=0$.   That is iff $\gamma\circ \rho=\delta'_X-\delta_X$ and $\rho:(X,\delta'_X)\rightmap{}(Y,\delta_Y)$ is a morphism in ${\cal Z}(Z)$. 
                          
For the computation of $b_1^{tw}(g)$, we have $b_1^{tw}(g)=g\circ \delta_E+\delta_Y\circ g+R(g)$, where   $R(g)=b_1^{ad}(g^1)+\sum_{\scriptsize\begin{matrix}i_0,i_1\geq 0\\ 
                          i_0+i_1\geq 2\end{matrix}} b^{ad}_{i_0+i_1+1}(\delta_Y^{\otimes i_1}\otimes g^1\otimes (\delta_E^1)^{\otimes i_0})$. We also have 
$$g\circ \delta_E+\delta_Y\circ g=
\begin{pmatrix}
\rho\circ \delta_X+\delta_Y\circ \rho \hbox{ \,,}&
\hueca{I}_Y\circ \delta_Y+\delta_Y\circ \hueca{I}_Y+\rho\circ\gamma\\
                                  \end{pmatrix}$$
and
$$R(g)=\begin{pmatrix} 
       b_1^{ad}(\rho)+\sum_{\scriptsize\begin{matrix}i_0,i_1\geq 0\\ 
                          i_0+i_1\geq 2\end{matrix}} b^{ad}_{i_0+i_1+1}(\delta_Y^{\otimes i_1}\otimes\rho\otimes \delta_X^{\otimes i_0}){\, ,}&0 \\
       \end{pmatrix}.$$
 Since $\hueca{I}_Y\circ \delta_Y+\delta_Y\circ \hueca{I}_Y=0$,  we obtain  that $b_1^{tw}(g)=0$ if and only if $\rho\circ\gamma=0$ and 
 $\rho\circ \delta_X+\delta_Y\circ \rho+ b_1^{ad}(\rho)+ \sum_{\scriptsize\begin{matrix}i_0,i_1\geq 0\\ 
                          i_0+i_1\geq 2\end{matrix}} b^{ad}_{i_0+i_1+1}(\delta_Y^{\otimes i_1}\otimes\rho\otimes \delta_X^{\otimes i_0})=0$.   That is iff $\rho\circ\gamma=0$ and $\rho:(X,\delta_X)\rightmap{}(Y,\delta_Y)$ is a morphism in ${\cal Z}(Z)$.      
\end{proof}

Similarly, we have the following. 

\begin{lemma}\label{L: caractizacion de pares en Z(hat(B)) con derecha variada}
 Assume that $(Y,\delta_Y)$ is an object in $\tw(Z)$ and  consider a sequence of morphisms in $\tw(Z)$ of the form 
 $$(X,\delta_X)\rightmap{f}(E,\delta_E)\rightmap{g}(Y,\delta'_Y),$$
with $E=X\oplus Y$, 
$\delta_E=\begin{pmatrix}
           \delta_X&\gamma\\ 0&\delta_Y\\
           \end{pmatrix}$, $f=(\hueca{I}_X,-\rho)^t$,  $g=(\rho,\hueca{I}_Y)$, 
           where $\rho:X\rightmap{}Y$ is a morphism  in $\ad(Z)$ with degree $-1$ 
           and  $\gamma:Y\rightmap{}X$ is a strict morphism in $\ad(Z)$ with degree $0$. Then, the  sequence lies in ${\cal Z}(Z)$ iff
     \begin{enumerate}
      \item  $\rho:(X,\delta_X)\rightmap{}(Y,\delta_Y)$  is a morphism in ${\cal Z}(Z)$ and $\gamma\circ\rho=0$. 
      \item $\rho:(X,\delta_X)\rightmap{}(Y,\delta'_Y)$  is a morphism in ${\cal Z}(Z)$ and $\rho\circ\gamma=\delta'_Y-\delta_Y$. 
     \end{enumerate}
\end{lemma}

\begin{proof} This is similar to the proof of (\ref{L: caractizacion de pares en Z(hat(B)) con izquierda variada}).        
\end{proof}

\begin{definition}\label{D: def de conflation}
We will say that a composable pair of morphisms of ${\cal Z}(Z)$
$$\xi:\hbox{ \ }(X,\delta_X)\rightmap{f}(E,\delta_E)\rightmap{g}(Y,\delta_Y),$$
is a \emph{conflation in ${\cal Z}(Z)$} iff there is a finite sequence of pairs of composable morphisms $\xi_0,\xi_1\ldots,\xi_n$ in ${\cal Z}(Z)$ such that  
$$  \xi_0\rightmap{\simeq}\xi_1\leftmap{\simeq}\xi_2\rightmap{\simeq}\xi_3\leftmap{\simeq}\xi_4\rightmap{\simeq}\cdots \leftmap{\simeq}\xi_{n-1}\rightmap{\simeq}\xi_{n-1}\leftmap{\simeq}\xi_n,$$
where $\xi=\xi_0$ and $\xi_n$ is a canonical conflation. In this case, we say that \emph{$\xi$ transforms into the canonical conflation $\xi_n$}. 
\end{definition}

\begin{lemma}\label{L: variando deltas en inicios de conflaciones}
Assume that we have an object $(X,\delta_X)$ in ${\cal Z}(Z)$ and the following  sequence of morphisms in ${\cal Z}(Z)$:  
 $$\xi:\hbox{ \ }(X,\delta'_X)\rightmap{f}(E,\delta_E)\rightmap{g}(Y,\delta_Y),$$
where $E=X\oplus Y$, 
$\delta_E=\begin{pmatrix}
           \delta_X&\gamma\\ 0&\delta_Y\\
           \end{pmatrix}$, $f=(\hueca{I}_X,-\rho)^t$, and  $g=(\rho,\hueca{I}_Y)$,  with $\rho:X\rightmap{}Y$ and  $\gamma:Y\rightmap{}X$ morphisms in $\ad(Z)$ with degree $-1$  and $0$, respectively. 
           Then, if $\gamma:Y\rightmap{}X$ is strict and 
           $\delta'_X\circ \gamma+\gamma\circ \delta_Y=0$, we have an isomorphism $h:(E,\delta_E)\rightmap{}(E_1,\delta_{E_1})$ and a  commutative diagram in ${\cal Z}(Z)$:
 $$\begin{matrix}
   (X,\delta'_X)&\rightmap{f}&(E,\delta_E)&\rightmap{g}&(Y,\delta_Y)\hfill\\
   \shortlmapdown{\hueca{I}_X}&&\shortrmapdown{h}&&\shortrmapdown{\hueca{I}_Y}\\
   (X,\delta'_X)&\rightmap{f_1}&(E_1,\delta_{E_1})&\rightmap{g_1}&(Y,\delta_Y)\hfill\\
   \end{matrix}$$                              
 where $E_1=X\oplus Y$, $\delta_{E_1}=\begin{pmatrix}
                                \delta'_X&\gamma\\ 0&\delta_Y\\
                               \end{pmatrix}$, $f_1=(\hueca{I}_X,0)^t$, and  $g_1=(0,\hueca{I}_Y)$.
 Since the lower row in the diagram is a canonical conflation, the upper row is a conflation. 
\end{lemma}

\begin{proof} The morphism $\gamma:Y\rightmap{}X$ of $\ad(Z)$ gives rise to the following two different morphisms   
$\gamma':(Y,\delta_Y)\rightmap{}(X,\delta'_X)$   and $\gamma:(Y,\delta_Y)\rightmap{}(X,\delta_X)$ in $\tw(Z)$,
We agree to add a prime to the first one to distinguish them. Thus, from 
 (\ref{L: diferenciales de la suma directa en tw(hat(B))}), we already know that $b_1^{tw}(\gamma)=0$, and, by assumption and (\ref{R: cuando un estricto esta en Z(hat(B))}), we have
 $b_1^{tw}(\gamma')=\delta'_X\circ \gamma+\gamma\circ \delta_Y=0$.
    Then, we get that the pair $(E_1,\delta_{E_1})$ such that $E_1=X\oplus Y$ and $\delta_{E_1}=
  \begin{pmatrix}  \delta'_X&\gamma\\
  0&\delta_Y\end{pmatrix}$ is an object of ${\cal Z}(Z)$. 

From (\ref{L: caractizacion de pares en Z(hat(B)) con izquierda variada}), we know that we can apply  (\ref{L: isos triangulares y con esquina estricta}) to the morphism $\rho:X\rightmap{}Y$ of $\ad(Z)$, and we have an isomorphism in ${\cal Z}(Z)$ of the form 
$$h:=\begin{pmatrix}
                     \hueca{I}_X&0\\ \rho&\hueca{I}_Y 
                     \end{pmatrix}:(E,\delta_E)\rightmap{}(E_1,\delta_{E_1}).$$
It remains to show that the diagram in the statement of the lemma commutes. Since $g_1$ is strict, we have 
$g_1\star h=g_1\circ h
=\begin{pmatrix}\hueca{I}_Y\circ \rho,\hueca{I}_Y\circ\hueca{I}_Y\end{pmatrix}
   =
   (\rho,\hueca{I}_Y)=g$. 

Let us verify that $h\star f=f_1$. Since $h\circ f=f_1$, it will be enough to show that $h\star f=h\circ f$. For this we use (\ref{R: desarmado f en estricto + otro simplifica btw's})(2). So, consider $f^0=(\hueca{I}_X,0)^t$, $f^1=(0,-\rho)^t$, $h^0=\begin{pmatrix}\hueca{I}_X&0\\ 0&\hueca{I}_Y\end{pmatrix}$, 
$h^1=\begin{pmatrix}0&0\\ \rho&0\end{pmatrix}$, $\delta_E^0=\begin{pmatrix}0&\gamma\\ 0&0\\ \end{pmatrix}$,
$\delta_E^1=\begin{pmatrix}\delta_X&0\\ 0&\delta_Y\\ \end{pmatrix}$, $\delta_{E_1}^0=\begin{pmatrix}0&\gamma\\ 0&0\\ \end{pmatrix}$, and $\delta_{E_1}^0=\begin{pmatrix}\delta'_X&0\\ 0&\delta_Y\\ \end{pmatrix}$. Then, we have $h\star f=h\circ f+R(h,f)$, where 
$R(h,f)=\sum_{\scriptsize\begin{matrix}i_0,i_1,i_2\geq 0\\ 
                          i_0+i_1+i_2\geq 1\end{matrix}} b^{ad}_{i_0+i_1+i_2+2}((\delta_{E_1}^1)^{\otimes i_2}\otimes h^1\otimes  (\delta_E^1)^{\otimes i_1}\otimes f^1\otimes (\delta'_X)^{\otimes i_0})$.
Since each tensor factor $h^1\otimes  (\delta_E^1)^{\otimes i_1}\otimes f^1$ is zero, we obtain $R(h,f)=0$, so $h\star f=h\circ f$ as wanted.  
\end{proof}

Similarly, we have the following.

\begin{lemma}\label{L: variando deltas en finales de conflaciones}
Assume that we have an object $(Y,\delta_Y)$ in ${\cal Z}(Z)$ and 
the following  sequence of morphisms in ${\cal Z}(Z)$: 
 $$\xi:\hbox{ \ }(X,\delta_X)\rightmap{f}(E,\delta_E)\rightmap{g}(Y,\delta'_Y),$$
where $E=X\oplus Y$, 
$\delta_E=\begin{pmatrix}
           \delta_X&\gamma\\ 0&\delta_Y\\
           \end{pmatrix}$, 
           $f=(\hueca{I}_X,-\rho)^t$, and  $g=(\rho,\hueca{I}_Y)$, 
           with   $\rho:X\rightmap{}Y$ and  $\gamma:Y\rightmap{}X$ morphisms in $\ad(Z)$ with degree $-1$  and $0$, respectively. 
           Then, if $\gamma:Y\rightmap{}X$ is strict and 
           $\delta_X\circ\gamma+\gamma\circ\delta'_Y=0$, we have an isomorphism $h:(E,\delta_E)\rightmap{}(E_1,\delta_1)$ and a commutative diagram in ${\cal Z}(Z)$:
 $$\begin{matrix}
  (X,\delta_X)&\rightmap{f}&(E,\delta_E)&\rightmap{g}&(Y,\delta'_Y)\\
   \shortlmapdown{\hueca{I}_X}&&\shortrmapdown{h}&&\shortrmapdown{\hueca{I}_Y}\\
     (X,\delta_X)&\rightmap{f_1}&(E_1,\delta_{E_1})&\rightmap{g_1}&(Y,\delta'_Y)\\
   \end{matrix}$$                              
 where $E_1=X\oplus Y$, $\delta_{E_1}=\begin{pmatrix}
                                \delta_X&\gamma\\ 0&\delta'_Y\\
                               \end{pmatrix}$, $f_1=(\hueca{I}_X,0)^t$, and  $g_1=(0,\hueca{I}_Y)$.
 Since the lower row of the diagram is a canonical conflation, the upper row is a conflation. 
\end{lemma}

\begin{proof} Similar to the proof of (\ref{L: variando deltas en inicios de conflaciones}), now using (\ref{L: caractizacion de pares en Z(hat(B)) con derecha variada}) and (\ref{L: isos triangulares sup entre objetos con esquina estricta}).                    
\end{proof}

\section{$(b,\nu)$-algebras}\label{S: (b,nu)-algebras}

In the following, we examine a special type of $b$-algebras $\hat{Z}=(\hat{Z},\{\hat{b}_n\}_{n\in \hueca{N}})$, over elementary $k$-algebras with enough idempotents $\hat{S}=(\hat{S},\{e_u\}_{u\in \hat{\cal P}})$, which admit a special two-sided action of an automorphism $\nu$ of $\hat{S}$.

\begin{definition}\label{D: bimod with translation} 
We say that a (unitary) graded $\hat{S}$-$\hat{S}$-bimodule $\hat{Z}$, over an elementary $k$-algebra with enough idempotents $\hat{S}=(\hat{S},\{e_u\}_{u\in \hat{\cal P}})$,  \emph{admits a two-sided translation $\nu$} iff the following two conditions hold:
\begin{enumerate}
 \item $\nu:\hat{S}\rightmap{}\hat{S}$ is an infinite order automorphism of $k$-algebras with enough idempotents \emph{acting freely on} $\{e_u\mid u\in \hat{\cal P}\}$  (that is $\nu^t(e_u)\not=e_u$, for all $u\in \hat{\cal P}$ and $t\in \hueca{Z}\setminus\{0\}$). In particular, $\nu$ induces a permutation of $\hat{\cal P}$ such that $\nu(e_u)=e_{\nu(u)}$, for all $u\in \hat{\cal P}$.  
 \item The infinite cyclic group $\langle \nu\rangle$ acts by the left and by the right on $\hat{Z}$ in such a way that the left and right actions by $\nu$ on the graded vector space $\hat{Z}$ are homogeneous $k$-linear automorphisms with degree $-1$ and they satisfy the following  for any $a\in \hat{Z}$ and $u\in \hat{\cal P}$:
 \begin{enumerate}
 \item $(\nu a)\nu=\nu(a\nu)$;
 \item $\nu(e_ua)=\nu(e_u)\nu a$ \hbox{ \ and \ } $\nu(ae_u)=(\nu a)e_u$;
 \item $(ae_u)\nu=(a\nu)\nu^{-1}(e_u)$ \hbox{ \ and \ } $(e_ua)\nu=e_u(a\nu)$.
 \end{enumerate}
 Notice that  ${\rm (a)}$ is equivalent to $(\nu^sa)\nu^t=\nu^s(a\nu^t)$, for all $s,t\in \hueca{Z}$. 
 \end{enumerate}
\end{definition}

\begin{remark}\label{R: las rho's de hat(Z) con accion}\label{R: la rho y los idempotentes} Given a 
graded $\hat{S}\g \hat{S}$-bimodule $\hat{Z}$ with a two-sided translation $\nu$, as above, we choose
 a complete set of representatives ${\cal P}$ of the $\langle \nu\rangle$-orbits of $\hat{\cal P}$ and set $\hat{\cal P}_s:=\nu^s({\cal P})$, for all $s\in \hueca{Z}$. In the following, we keep the set ${\cal P}$ fixed.
  Since $\nu$ acts freely on $\hat{\cal P}$, we have that $\hat{\cal P}_s\cap\hat{\cal P}_t=\emptyset$, whenever $s\not=t$. 
 $$\hbox{ For }s,t\in \hueca{Z},\hbox{ define }
 \hat{Z}_{s,t}:=\bigoplus_{v\in \hat{\cal P}_s,u\in \hat{\cal P}_t}e_v\hat{Z}
 e_u, \hbox{ thus, we have }\hat{Z}=\bigoplus_{s,t\in \hueca{Z}}\hat{Z}_{s,t}.$$ 
 Notice that, whenever $a\in \hat{Z}_{s,t}$, we have $\nu a\in \hat{Z}_{s+1,t}$ and $a\nu\in \hat{Z}_{s,t-1}$. For $s,t\in \hueca{Z}$, we consider the linear homogeneous isomorphism of degree  $\vert \rho_{s,t}\vert=s-t$ 
 $$\rho_{s,t}:\hat{Z}_{s,t}\rightmap{}\hat{Z}_{0,0} \hbox{ \ defined by \ } \rho_{s,t}(a)=\nu^{-s}a\nu^t.$$ 
 For each $\underline{s}=(s_0,s_1,\ldots,s_n)\in \hueca{Z}^{n+1}$, we consider the linear map 
 $$\rho_{\underline{s}}:=\rho_{s_0,s_1}\otimes\cdots\otimes\rho_{s_{n-1},s_n}:\hat{Z}_{s_0,s_1}\otimes \hat{Z}_{s_1,s_2}\otimes\cdots\otimes \hat{Z}_{s_{n-1},s_n}\rightmap{}\hat{Z}_{0,0}^{\otimes n}.$$
 We will also write $\hat{Z}_{\underline{s}}:=\hat{Z}_{s_0,s_1}\otimes \hat{Z}_{s_1,s_2}\otimes\cdots\otimes \hat{Z}_{s_{n-1},s_n}$, thus $\rho_{\underline{s}}:\hat{Z}_{\underline{s}}\rightmap{}\hat{Z}_{0,0}^{\otimes n}$ is homogeneous with degree $\vert \rho_{\underline{s}}\vert=s_0-s_n$. 
 
 Notice that, 
 for $s,t\in \hueca{Z}$, $a\in \hat{Z}_{s,t}$, and $u,v\in \hat{\cal P}$, we have $$\rho_{s,t}(e_vae_u)=e_{\nu^{-s}(v)}\rho_{s,t}(a)e_{\nu^{-t}(u)}.$$
\end{remark}

\begin{definition}\label{D: (b,sigma)-algebra} 
A $b$-algebra  $\hat{Z}=(\hat{Z},\{\hat{b}_n\}_{n\in \hueca{N}})$ over an elementary $k$-algebra with enough idempotents $\hat{S}=(\hat{S},\{e_u\}_{u\in \hat{\cal P}})$ is a \emph{$(b,\nu)$-algebra} iff the unitary graded $\hat{S}$-$\hat{S}$-bimodule $\hat{Z}$  admits a two-sided translation $\nu$ and there is a set of representatives ${\cal P}$ of the $\langle \nu\rangle$-orbits in $\hat{\cal P}$ such that, with the  notations of (\ref{R: las rho's de hat(Z) con accion}), we have    
$${\hat{b}_n \vert}_{\hat{Z}_{\underline{s}} }=
(-1)^{s_0-s_n}
\rho^{-1}_{s_0,s_n}b_n\rho_{\underline{s}},$$
for all $n\in \hueca{N}$ and $\underline{s}=(s_0,\ldots,s_n)\in \hueca{Z}^{n+1}$, where $Z:=\hat{Z}_{0,0}$ and $b_n:Z^{\otimes n}\rightmap{}Z$ is the restriction of $\hat{b}_n$. 
\end{definition}

We will show in a moment, in the proof of (\ref{P: construccion de hat(Z) a partir de Z}), how the preceding notion relates to Keller's construction $\hueca{Z}{\cal A}$, for an $A_\infty$-category ${\cal A}$ with strict units. Before this, we give some elementary useful arithmetical properties of these $(b,\nu)$-algebras.
In the following paragraphs, unless we clearly indicate otherwise, we assume that $\hat{Z}=(\hat{Z},\{\hat{b}_n\}_{n\in \hueca{N}})$ is a $(b,\nu)$-algebra as in (\ref{D: (b,sigma)-algebra}).

\begin{lemma}\label{R: hat(b)n evaluado}  
Given  homogeneous elements 
$a_1\in \hat{Z}_{s_0,s_1}, \ldots,a_n\in \hat{Z}_{s_{n-1},s_n}$, we have 
$\hat{b}_n(a_1\otimes\cdots\otimes a_n)=(-1)^{z}\rho^{-1}_{s_0,s_n}b_n(\rho_{s_0,s_1}(a_1)\otimes\cdots\otimes\rho_{s_{n-1},s_n}(a_n)),$
where $z=s_0-s_n+\sum_{l=1}^{n-1}(s_l-s_n)\vert a_l\vert$. 
\end{lemma}

\begin{proof} It follows from the application of the Koszul sign convention. 
\end{proof}

\begin{proposition} For $n\in \hueca{N}$ and $\underline{s},\underline{t}\in \hueca{Z}^{n+1}$, 
we have 
$${\hat{b}_n\vert}_{\hat{Z}_{\underline{s}} }=
(-1)^{s_0-t_0-(s_n-t_n)}
\rho^{-1}_{s_0-t_0,s_n-t_n}
{\hat{b}_n \vert}_{ \hat{Z}_{\underline{t}}  }\rho^{-1}_{\underline{t}}
\rho_{\underline{s}}$$
\end{proposition}     

\begin{proof} Denote by $\Delta$ the expression on the right of the equality above. From the definition of $(b,\nu)$-algebra, we have   
${\hat{b}_n \vert}_{\hat{Z}_{\underline{t}} }=
(-1)^{t_0-t_n}
\rho^{-1}_{t_0,t_n}b_n\rho_{\underline{t}}$. So, we have 
$$\begin{matrix}
\Delta
&=&
(-1)^{s_0-t_0-(s_n-t_n)}
\rho^{-1}_{s_0-t_0,s_n-t_n}(-1)^{t_0-t_n}
\rho^{-1}_{t_0,t_n}b_n\rho_{\underline{t}}\rho^{-1}_{\underline{t}}
\rho_{\underline{s}}\hfill\\
&=&
(-1)^{s_0-s_n}
\rho^{-1}_{s_0,s_n}b_n\rho_{\underline{s}}= {\hat{b}_n\vert}_{\hat{Z}_{\underline{s}}}. \hfill\\
\end{matrix}
$$
\end{proof}

\begin{corollary}\label{L: nu por hat(b)n}\label{R: hat(b)n por conjugados de nu internos}\label{R: hat(b)n por nu a la -1}  
Given  homogeneous elements 
$a_1\in \hat{Z}_{s_0,s_1}, \ldots,a_n\in \hat{Z}_{s_{n-1},s_n}$,  the following equalities hold.
\begin{enumerate}
 \item[(1)]$\hbox{\hskip1cm}
 \hat{b}_n(\nu a_1\otimes\cdots\otimes a_n)=(-1)^{s_1-s_n+1}\nu\hat{b}_n( a_1\otimes\cdots\otimes a_n);$
 
\item[(2)] $\hbox{\hskip1cm}
 \hat{b}_n( a_1\otimes\cdots\otimes a_n\nu^{-1})=(-1)^{1+\sum_{l=1}^{n-1}\vert a_l\vert}\hat{b}_n( a_1\otimes\cdots\otimes a_n)\nu^{-1};$
 
\item[(3)]
For $n\geq 2$ and $l\in [1,n-1]$,  we have 
 that 
 $$ \hat{b}_n( a_1\otimes a_2\otimes \cdots\otimes a_{l-1}\otimes a_l\nu^{-1}\otimes \nu a_{l+1}\otimes a_{l+2}\otimes  \cdots\otimes a_{n-1}\otimes a_n)$$
 coincides with 
  $(-1)^{\vert a_l\vert+s_l-s_{l+1}+1}\hat{b}_n( a_1\otimes a_2\otimes\cdots\otimes a_l\otimes a_{l+1}\otimes \cdots\otimes a_n)$.
\end{enumerate}
\end{corollary}

\begin{proof}(1): Take $\underline{s}=(s_0,s_1,\ldots,s_n)\in \hueca{Z}^{n+1}$ and set $\underline{t}=(s_0+1,s_1,\ldots,s_n)$. Denote by $\nu_L:\hat{Z}_{s_0,s_1}\rightmap{}\hat{Z}_{t_0,t_1}$ the left multiplication by $\nu$. Then,
$$
\rho_{\underline{t}}(\nu_L\otimes id^{\otimes(n-1)})=(-1)^{t_1-t_n}(\rho_{t_0,t_1}\nu_L\otimes\rho_{t_1,t_2}\otimes\cdots\otimes \rho_{t_{n-1},t_n})\hfill\\
=
(-1)^{s_1-s_n}\rho_{\underline{s}}.$$
 Hence, we have $\rho_{\underline{t}}^{-1}\rho_{\underline{s}}
 =
 (-1)^{s_1-s_n}(\nu_L\otimes id^{\otimes(n-1)})$. 
 So, in this case, we get
 $${\hat{b}_n\vert}_{\hat{Z}_{\underline{s}} }=(-1)^{1+s_1-s_n}\rho^{-1}_{-1,0}{\hat{b}_n\vert}_{\hat{Z}_{\underline{t}} }(\nu_L\otimes id^{\otimes(n-1)}).$$
 Therefore, given a typical generator $a_1\otimes a_2\otimes\cdots\otimes a_n\in \hat{Z}_{\underline{s}}$, we obtain 
 $$\hat{b}_n(a_1\otimes a_2\otimes\cdots\otimes a_n)=
 (-1)^{1+s_1-s_n}\nu^{-1}\hat{b}_n(\nu a_1\otimes a_2\otimes\cdots\otimes a_n).$$

\noindent(2): Take $\underline{s}=(s_0,s_1,\ldots,s_n)\in \hueca{Z}^{n+1}$ and set $\underline{t}=(s_0,s_1,\ldots,s_n+1)$. Denote by $\nu^{-1}_R:\hat{Z}_{s_{n-1},s_n}\rightmap{}\hat{Z}_{t_{n-1},t_n}$ the right multiplication by $\nu^{-1}$. Then,
$$\begin{matrix}
\rho_{\underline{t}}( id^{\otimes(n-1)}\otimes \nu^{-1}_R)
&=&
\rho_{t_0,t_1}\otimes\rho_{t_1,t_2}\otimes\cdots\otimes \rho_{t_{n-1},t_n}\nu_R^{-1}
=\rho_{\underline{s}}.\hfill\\
 \end{matrix}$$
 Hence, we have $\rho_{\underline{t}}^{-1}\rho_{\underline{s}}=( id^{\otimes(n-1)}\otimes\nu^{-1}_R)$. So, in this case, we get
 $${\hat{b}_n\vert}_{\hat{Z}_{\underline{s}} }=-\rho^{-1}_{0,-1}{\hat{b}_n\vert}_{\hat{Z}_{\underline{t}} }( id^{\otimes(n-1)}\otimes\nu_R^{-1}).$$
 Therefore, given a typical generator $a_1\otimes a_2\otimes\cdots\otimes a_n\in \hat{Z}_{\underline{s}}$, with homogeneous tensor factors, we obtain 
 $$\hat{b}_n(a_1\otimes a_2\otimes\cdots\otimes a_n)=
 (-1)^{1+\sum_{l=1}^{n-1}\vert a_l\vert}\hat{b}_n(a_1\otimes a_2\otimes\cdots\otimes a_n\nu^{-1})\nu.$$

 \noindent(3): Take $\underline{s}=(s_0,s_1,\ldots,s_n)\in \hueca{Z}^{n+1}$ and set
$$\underline{t}=(s_0,\ldots,s_{l-1},s_l+1,s_{l+1},\ldots,s_n).$$
Denote by  $\nu^{-1}_R:\hat{Z}_{s_{l-1},s_l}\rightmap{}\hat{Z}_{t_{l-1},t_l}$ the right  multiplication by $\nu^{-1}$ and by $\nu_L:\hat{Z}_{s_l,s_{l+1}}\rightmap{}\hat{Z}_{t_l,t_{l+1}}$ the left multiplication by $\nu$. Then, 
$$
\rho_{\underline{t}}( id^{\otimes(l-1)}\otimes \nu^{-1}_R\otimes \nu_L\otimes id^{\otimes(n-l-1)})=(-1)^{t_l-t_{l+1}}\rho_{\underline{s}}=(-1)^{s_l+1-s_{l+1}}\rho_{\underline{s}}.$$

 Hence, we have $\rho_{\underline{t}}^{-1}\rho_{\underline{s}}=(-1)^{s_l+1-s_{l+1}}
 ( id^{\otimes(l-1)}\otimes\nu^{-1}_R\otimes \nu_L\otimes id^{\otimes(n-l-1)})$. So, in this case, we get
 $${\hat{b}_n\vert}_{\hat{Z}_{\underline{s}}}
 =
 (-1)^{s_l+1-s_{l+1}}{\hat{b}_n\vert}_{\hat{Z}_{\underline{t}} }
 (id^{\otimes(l-1)}\otimes \nu^{-1}_R\otimes\nu_L\otimes id^{\otimes(n-l-1)}).$$
 From this (3) follows.  
\end{proof}

From the last part of the preceding result, we have the following.

\begin{corollary}\label{C: transito de nu y nu-1 entre tensores de bn}
For $n\geq 2$ and $l\in [1,n-1]$, the following holds.
\begin{enumerate}
 \item Given homogeneous 
 $a_1\in \hat{Z}_{s_0,s_1},\ldots,
 a_l\in \hat{Z}_{s_{l-1},s_l},
 a_{l+1}\in \hat{Z}_{s_l-1,s_{l+1}},a_{l+2}\in \hat{Z}_{s_{l+1},s_{l+2}},\ldots,a_n\in\hat{Z}_{s_{n-1},s_n}$ we have that 
 $$ \hat{b}_n( a_1\otimes a_2\otimes \cdots\otimes a_{l-1}\otimes a_l\otimes \nu a_{l+1}\otimes a_{l+2}\otimes  \cdots\otimes a_{n-1}\otimes a_n)$$
 coincides with 
  $(-1)^{\vert a_l\vert+1+s_{l}-s_{l+1}}
  \hat{b}_n( a_1\otimes a_2\otimes\cdots\otimes a_l\nu\otimes a_{l+1}\otimes \cdots\otimes a_n)$.

  \item Given homogeneous 
  $a_1\in \hat{Z}_{s_0,s_1},\ldots,a_l\in \hat{Z}_{s_{l-1},s_l},
  a_{l+1}\in \hat{Z}_{s_{l}+1,s_{l+1}},
  a_{l+2}\in \hat{Z}_{s_{l+1},s_{l+2}}\ldots,a_n\in\hat{Z}_{s_{n-1},s_n}$, we have that 
 $$ \hat{b}_n( a_1\otimes a_2\otimes \cdots\otimes a_{l-1}\otimes a_l\nu^{-1}\otimes  a_{l+1}\otimes a_{l+2}\otimes  \cdots\otimes a_{n-1}\otimes a_n)$$
 coincides with 
  $(-1)^{\vert a_l\vert+1+s_l-s_{l+1}}\hat{b}_n( a_1\otimes a_2\otimes\cdots\otimes a_l\otimes \nu^{-1}a_{l+1}\otimes \cdots\otimes a_n)$.
\end{enumerate}
\end{corollary}

\begin{corollary}\label{C: circ vs nu y nu-1}
For $n=2$, from (\ref{L: nu por hat(b)n}) and (\ref{C: transito de nu y nu-1 entre tensores de bn}),  the following holds.
\begin{enumerate}
 \item For any homogeneous $a_1\in \hat{Z}_{s_0,s_1}$ and $a_2\in \hat{Z}_{s_1,s_2}$ we have 
 $$\begin{matrix}(\nu a_1)\circ a_2&=&(-1)^{s_1-s_2+1}\nu(a_1\circ a_2)\\ 
 a_1\circ (a_2\nu^{-1})&=&(-1)^{\vert a_1\vert +1}(a_1\circ a_2)\nu^{-1}\\ 
 (a_1\nu^{-1})\circ (\nu a_2)&=&(-1)^{s_1-s_2+\vert a_1\vert+ 1}a_1\circ a_2.\end{matrix}$$
 \item For any homogeneous $a_1\in \hat{Z}_{s_0,s_1}$ and $a_2\in \hat{Z}_{s_1-1,s_2}$ we have  $$a_1\circ (\nu a_2)=(-1)^{s_1-s_2+\vert a_1\vert +1}(a_1\nu)\circ a_2.$$ 
 \item For any homogeneous $a_1\in \hat{Z}_{s_0,s_1}$ and $a_2\in \hat{Z}_{s_1+1,s_2}$ we have 
 $$(a_1\nu^{-1})\circ a_2=(-1)^{s_1-s_2+\vert a_1\vert+1}a_1\circ (\nu^{-1}a_2).$$
\end{enumerate}
\end{corollary}

\begin{definition}\label{D: a[1] y a[-1]}
For  any element $a\in \hat{Z}$, define
$$a[1]:=\nu a\nu^{-1} \hbox{ \ \ and \ \   }a[-1]:=\nu^{-1}a\nu.$$
Then, for $a\in \hat{Z}_{s,t}$, we have $a[1]\in \hat{Z}_{s+1,t+1}$ and $a[-1]\in \hat{Z}_{s-1,t-1}$. Moreover, we have that $a[1][-1]=a=a[-1][1]$. If $a$ is homogeneous, so are $a[1]$ and $a[-1]$, and  $\vert a[1]\vert=\vert a\vert=\vert a[-1]\vert$. 
\end{definition}

\begin{lemma}\label{R: hat(b)n de tensores de trasladados}
For any  homogeneous elements 
$a_1\in \hat{Z}_{s_0,s_1}, \ldots,a_n\in \hat{Z}_{s_{n-1},s_n}$, we have 
 $ \hat{b}_n( a_1[1]\otimes a_2[1]\otimes\cdots\otimes a_n[1])=
 \hat{b}_n( a_1\otimes a_2\otimes \cdots\otimes a_n)[1].$
\end{lemma}
 
\begin{proof} Set $\Delta:=\hat{b}_n( a_1[1]\otimes a_2[1]\otimes\cdots\otimes a_n[1])$. Then,
from the preceding lemmas, we have 
 $$\begin{matrix}
  \Delta
  &=&
  (-1)^{z_0}\nu\hat{b}_n( a_1\nu^{-1}\otimes a_2[1]\otimes\cdots\otimes a_n[1])\hfill\\
  &=&
  (-1)^{z_1}\nu\hat{b}_n( a_1\otimes a_2\nu^{-1}\otimes a_3[1]\otimes\cdots\otimes a_n[1])\hfill\\
 &=&
  (-1)^{z_2}\nu\hat{b}_n( a_1\otimes a_2\otimes a_3\nu^{-1}\otimes a_4[1]\cdots\otimes a_n[1])\hfill\\ 
  &\vdots&
  \vdots&\hfill\\
  &=&
  (-1)^{z_{n-1}}\nu\hat{b}_n( a_1\otimes a_2\otimes \cdots\otimes a_n\nu^{-1})\hfill\\ 
  &=&
 (-1)^{z_n}\hat{b}_n( a_1\otimes a_2\otimes \cdots\otimes a_n)[1]\hfill\\   
   \end{matrix}$$
 where, modulo 2, we have  
 $$\begin{matrix}
   z_0
   &=&
   (s_1+1)-(s_n+1)+1=s_1-s_n+1\hfill\\
   z_1
   &=&
   z_0+\vert a_1\vert+s_1-(s_2+1)+1\hfill\\
   z_2
   &=&
    z_1+\vert a_2\vert+s_2-(s_3+1)+1\hfill\\
   \vdots&\vdots&\hbox{\hskip1.5cm}\vdots\hfill\\
   z_{n-1}&=&z_{n-2}+\vert a_{n-1}\vert +s_{n-1}-(s_n+1)+1=1+\sum_{l=1}^{n-1}\vert a_l\vert\hfill\\
   z_n&=&z_{n-1}+1+\sum_{l=1}^{n-1}\vert a_l\vert=0.\hfill\\
   \end{matrix}$$
\end{proof}

\begin{remark}\label{R: restriccion de (b,sigma)-algebra es b-algebra} Let $\hat{Z}$ be a $(b,\nu)$-algebra, as in (\ref{D: (b,sigma)-algebra}). Consider the $k$-subalgebra with enough idempotents $S:=\bigoplus_{u\in {\cal P}}ke_u$ of $\hat{S}$,  the $S$-$S$-bimodule $Z:=\hat{Z}_{0,0}$, and the restrictions $b_n:Z^{\otimes n}\rightmap{}Z$ of the morphisms $\hat{b}_n:\hat{Z}^{\otimes n}\rightmap{}\hat{Z}$. Then, we obtain a $b$-algebra $Z=(Z,\{b_n\}_{n\in \hueca{N}})$ over the elementary $k$-algebra  with enough idempotents $S=(S,\{e_u\}_{u\in {\cal P}})$. There, we are identifying the tensor products over $S$ implicit in $Z^{\otimes n}$ with the tensor products over $\hat{S}$ implicit in $\hat{Z}_{0,0}^{\otimes n}$.  
We call the $b$-algebra $Z=(Z,\{b_n\}_{n\in \hueca{N}})$ \emph{a section}  of the $(b,\nu)$-algebra $(\hat{Z},\{\hat{b}_n\}_{n\in \hueca{N}})$. 
\end{remark}

\begin{lemma}\label{L: hat(Z) nu-bimod y Z(P) b-alg --> hat(Z) b-nu-alg con seccion Z}
 Let $\hat{Z}$ be a 
graded $\hat{S}\g \hat{S}$-bimodule with a two-sided translation $\nu$, as in (\ref{D: bimod with translation}), and take  
 any complete set of representatives ${\cal P}$ of the $\langle \nu\rangle$-orbits of $\hat{\cal P}$ as in (\ref{R: la rho y los idempotentes}).  
Consider the $k$-subalgebra with enough idempotents $S:=\bigoplus_{u\in {\cal P}}ke_u$ of $\hat{S}$ and the $S$-$S$-bimodule $Z:=\hat{Z}_{0,0}$.
Furthermore, suppose that we have a $b$-algebra $(Z,\{b_n\}_{n\in \hueca{N}})$ over $S$. Then, there is a $(b,\nu)$-algebra  $\hat{Z}=(\hat{Z},\{\hat{b}_n\}_{n\in \hueca{N}})$, over  $\hat{S}$,  with section $(Z,\{b_n\}_{n\in \hueca{N}})$. 
\end{lemma}

 \begin{proof} For $n\in \hueca{N}$, we have $\hat{Z}^{\otimes n}=\bigoplus_{\underline{w}\in \hueca{Z}^{(n+1)}}\hat{Z}_{\underline{w}}$. So we can consider the linear maps  $\hat{b}_n:\hat{Z}^{\otimes n}\rightmap{}\hat{Z}$ such that 
 $${\hat{b}_n \vert}_{\hat{Z}_{\underline{w}} }=
(-1)^{w_0-w_n}
\rho^{-1}_{w_0,w_n}b_n\rho_{\underline{w}},
\hbox{ \ for all \ }\underline{w}=(w_0,\ldots,w_n)\in \hueca{Z}^{n+1}.$$ 
It is readily seen that each $\hat{b}_n$ is a homogeneous morphism of $\hat{S}\g\hat{S}$-bimodules with degree 1. We have to show that 
 $\hat{S}_n:=\sum_{\scriptsize\begin{matrix}r+s+t=n\\ s\geq 1;r,t\geq 0\end{matrix}} 
   \hat{b}_{r+1+t}(id^{\otimes r}\otimes \hat{b}_s\otimes id^{\otimes t})=0$,   
for each $n\in \hueca{N}$.    It is enough to show that $\hat{S}_n\vert_{\hat{Z}_{\underline{w}}}=0$, for all $\underline{w}\in \hueca{Z}^{(n+1)}$. Given integers $r,t\geq 0$ and $s\geq 1$ such that $r+s+t=n$, we have $\hat{Z}_{\underline{w}}=\hat{Z}_{\underline{w}^s}\otimes\hat{Z}_{\underline{w}^r}\otimes\hat{Z}_{\underline{w}^t}$, where $\underline{w}^r=(w_0,\ldots,w_r)$,
    $\underline{w}^s=(w_r,\ldots,w_{r+s})$, and $\underline{w}^t=(w_{r+s},\ldots,w_n)$. 
    Thus, we have ${\hat{b}_s \vert}_{\hat{Z}_{\underline{w}^s} }=
(-1)^{w_r-w_{r+s}}
\rho^{-1}_{w_r,w_{r+s}}b_s\rho_{\underline{w}^s}$. 
   Now, consider a typical summand  $\Delta:=\hat{b}_{r+1+t}(id^{\otimes r}\otimes \hat{b}_s\otimes id^{\otimes t})\vert_{\hat{Z}_{\underline{w}}}$  of $\hat{S}_n\vert_{\hat{Z}_{\underline{w}}}$. We obtain 
   $$\begin{matrix}
    \Delta&=&
    (-1)^{w_0-w_n}\rho_{w_0,w_n}^{-1}b_{r+1+t}(\rho_{\underline{w}^r}\otimes\rho_{w_r,w_{r+s}}\otimes\rho_{\underline{w}^t})(id^{\otimes r}\otimes \hat{b}_s\vert_{\hat{Z}_{\underline{w}^s}}\otimes id^{\otimes t})\hfill\\
    &=&
     (-1)^{w_0-w_n+\vert\rho_{\underline{w}^t}\vert}\rho_{w_0,w_n}^{-1}b_{r+1+t}(\rho_{\underline{w}^r}\otimes\rho_{w_r,w_{r+s}}(\hat{b}_s\vert_{\hat{Z}_{\underline{w}^s}})\otimes\rho_{\underline{w}^t}) \hfill\\
     &=&
     (-1)^{w_0-w_n+\vert\rho_{\underline{w}^t}\vert+w_r-w_{r+s}}\rho_{w_0,w_n}^{-1}b_{r+1+t}(\rho_{\underline{w}^r}\otimes b_s\rho_{\underline{w}^s}\otimes\rho_{\underline{w}^t}) \hfill\\
     &=&
     (-1)^z\rho_{w_0,w_n}^{-1}b_{r+1+t}(id^{\otimes r}\otimes b_s\otimes id^{\otimes t})(\rho_{\underline{w}^r}\otimes \rho_{\underline{w}^s}\otimes\rho_{\underline{w}^t}) \hfill\\
     &=&
     \rho_{w_0,w_n}^{-1}b_{r+1+t}(id^{\otimes r}\otimes b_s\otimes id^{\otimes t})\rho_{\underline{w}}\hfill\\
   \end{matrix}$$
where 
$
z=w_0-w_n+\vert\rho_{\underline{w}^t}\vert+w_r-w_{r+s}+\vert \rho_{\underline{w}^r}\vert$ is zero modulo 2. So, adding up,  we obtain  
$\hat{S}_n\vert_{\hat{Z}_{\underline{w}}}=\rho_{w_0,w_n}^{-1}\sum_{\scriptsize\begin{matrix}r+s+t=n\\ s\geq 1;r,t\geq 0\end{matrix}} b_{r+1+t}(id^{\otimes r}\otimes b_s\otimes id^{\otimes t})\rho_{\underline{w}}=0$. 
 \end{proof}

\begin{proposition}\label{P: construccion de hat(Z) a partir de Z} 
Given a $b$-algebra $Z=(Z,\{b_n\}_{n\in \hueca{N}})$ over the elementary $k$-algebra  with enough idempotents $S=(S,\{e_i\}_{i\in {\cal P}})$, we can associate naturally  a $(b,\nu)$-algebra  $\hat{Z}=(\hat{Z},\{\hat{b}_n\}_{n\in \hueca{N}})$, over an elementary $k$-algebra with enough idempotents $\hat{S}=(\hat{S},\{e_u\}_{u\in \hat{\cal P}})$,  with section $Z=(Z,\{b_n\}_{n\in \hueca{N}})$. 
\end{proposition}

\begin{proof} For each $0\not=t\in \hueca{Z}$, 
 fix a copy $S[t]$ of the $k$-algebra $S$, and set $S[0]:=S$.  Then, consider the $k$-algebra without unit $\hat{S}:=\coprod_{t\in \hueca{Z}}S[t]\subset \prod_{t\in \hueca{Z}}S[t]$, with product induced by the product of the $k$-algebra $\prod_{t\in \hueca{Z}}S[t]$.
 
 For $(t,i)\in \hat{\cal P}:=\hueca{Z}\times {\cal P}$, define
 $e_{(t,i)}:=\sigma_t(e_i)\in \hat{S}$,  
 where $\sigma_t:S\rightmap{}S[t]$ is a fixed isomorphism of $k$-algebras.  Then, $\{e_u\mid u\in \hat{\cal P}\}$ is a family of primitive orthogonal  idempotents of $\hat{S}$ such that $\hat{S}= \bigoplus_{u,v\in \hat{\cal P}}e_v\hat{S}e_u$, and we can consider the automorphism  $\nu:\hat{S}\rightmap{}\hat{S}$ of $k$-algebras with enough idempotents, defined by $\nu(e_{(t,i)})=e_{(t+1,i)}$, for all $(t,i)\in \hat{\cal P}$, which acts freely on $\{e_u\mid u\in \hat{\cal P}\}$. 
 
 Then, consider for each $s,t\in \hueca{Z}$, a copy $\hat{Z}_{s,t}$ of the graded $S$-$S$-bimodule $Z[s-t]$.   We fix, for each $s,t\in \hueca{Z}$, an isomorphism $\phi_{s,t}:\hat{Z}_{s,t}\rightmap{}Z[s-t]$ of graded $S$-$S$-bimodules. We agree that $\hat{Z}_{0,0}=Z$ and $\phi_{0,0}$ is the identity map. 

Notice that any graded $S$-$S$-bimodule  $W$ is a graded $S[s]$-$S[t]$-bimodule by restriction via the isomorphism $\sigma_s^{-1}:S[s]\rightmap{}S$ on the left and $\sigma_t^{-1}:S[t]\rightmap{}S$ on the right. If we denote by $\pi_l:\hat{S}\rightmap{}S[l]$ the canonical projection on the $l$-factor of $\hat{S}$, for $l\in \hueca{Z}$, we can consider  the graded $\hat{S}$-$\hat{S}$-bimodule obtained from $W$ by restriction of scalars through $\pi_s$ on the left and $\pi_t$ on the right. This holds for the graded $S$-$S$-bimodules $\hat{Z}_{s,t}$ and $Z[s-t]$, and the given isomorphism $\phi_{s,t}:\hat{Z}_{s,t}\rightmap{}Z[s-t]$ of $S$-$S$-bimodules becomes an isomorphism of $\hat{S}$-$\hat{S}$-bimodules. Then, we consider the graded $\hat{S}\g\hat{S}$-bimodule 
$$\hat{Z}:=\bigoplus_{s,t\in \hueca{Z}}\hat{Z}_{s,t}.$$
Therefore, for $a\in \hat{Z}_{s,t}$ and $(s_1,j_1),(t_1,i_1)\in \hat{\cal P}$, we have 
$$\phi_{s,t}(e_{(s_1,j_1)}ae_{(t_1,i_1)})=
\begin{cases}
e_{j_1}\phi_{s,t}(a)e_{i_1}&\hbox{ if } s_1=s \hbox{ and }t_1=t\\
0&\hbox{ if } s_1\not=s \hbox{ or }t_1\not=t.\\
\end{cases}$$
For $l\in \hueca{Z}$,  denote by  $\tau(l):=Z[l]\rightmap{}Z$  the canonical homogeneous isomorphism induced by the identity map, so we have  $\vert \tau(l)\vert= l$. The following holds.

 For $s,t\in \hueca{Z}$, consider the homogeneous isomorphism of graded $S$-$S$-bimodules 
$$\rho_{s,t}:=\tau(s-t)\phi_{s,t}:\hat{Z}_{s,t}\rightmap{}Z,$$
which has degree $\vert\rho_{s,t}\vert=s-t$. 
Then, for $a\in \hat{Z}_{t,s}$ and $(s_1,j_1),(t_1,i_1)\in \hat{\cal P}$,  we have 
$$\rho_{s,t}(e_{(s_1,j_1)}ae_{(t_1,i_1)})=
\begin{cases}
e_{j_1}\rho_{s,t}(a)e_{i_1}&\hbox{ if } s_1=s \hbox{ and }t_1=t\\
0&\hbox{ if } s_1\not=s \hbox{ or }t_1\not=t.\\
\end{cases}$$

Now, let us specify the left and right actions of $\langle \nu\rangle$ on the $\hat{S}\g \hat{S}$-bimodule $\hat{Z}$. They are determined by the following formulas, if $a\in \hat{Z}_{s,t}$,  
$$\nu a:=\rho_{s+1,t}^{-1}\rho_{s,t}(a)\in \hat{Z}_{s+1,t} \hbox{ \ and \ } a\nu:=\rho_{s,t-1}^{-1}\rho_{s,t}(a)\in \hat{Z}_{s,t-1}.$$
So, indeed, $\nu$ acts on each side of the graded space $\hat{Z}$ as an homogeneous $k$-linear automorphism with degree $-1$. We readily see that $(\nu a)\nu=\nu(a\nu)$.

Now, we proceed to verify condition (2)(b) of (\ref{D: bimod with translation}). For the first part, let $u=(s_1,i_1)$ and $a\in \hat{Z}_{s,t}$, thus $\nu a\in \hat{Z}_{s+1,t}$. If $s\not=s_1$, we have 
$\nu(e_ua)=\rho_{s+1,t}^{-1}\rho_{s,t}(e_{(s_1,i_1)}a)=0
=e_{(s_1+1,i_1)}\nu a=\nu(e_u)\nu a$. 
While, if $s=s_1$, we have 
$\nu(e_ua)=\rho_{s+1,t}^{-1}\rho_{s,t}(e_{(s,i_1)}a)=
\rho_{s+1,t}^{-1}(e_{i_1}\rho_{s,t}(a))
=e_{i_1}\rho_{s+1,t}^{-1}(\rho_{s,t}(a))=e_{i_1}\nu a =e_{(s+1,i_1)}\nu a=\nu(e_u)\nu a$. 

Now assume $u=(t_1,i_1)$ and $a\in \hat{Z}_{s,t}$, thus $\nu a\in \hat{Z}_{s+1,t}$. If $t\not=t_1$, we have 
$\nu(ae_u)=\rho_{s+1,t}^{-1}\rho_{s,t}(ae_{(t_1,i_1)})=0
=(\nu a)e_{(t_1,i_1)}=(\nu a)e_u$. 
If $t=t_1$, we have 
$\nu(ae_u)=\rho_{s+1,t}^{-1}\rho_{s,t}(ae_{(t,i_1)})=
\rho_{s+1,t}^{-1}(\rho_{s,t}(a)e_{i_1})
=\rho_{s+1,t}^{-1}(\rho_{s,t}(a))e_{i_1}=(\nu a)e_{i_1}
=(\nu a)e_{(t,i_1)}=(\nu a )e_u$. 
The condition (\ref{D: bimod with translation})(2)(c) is verified similarly. 
Thus, the $\hat{S}$-$\hat{S}$ bimodule $\hat{Z}$ over the elementary $k$-algebra with enough idempotents $\hat{S}=(\hat{S},\{e_u\}_{u\in \hat{\cal P}})$  admits a two-sided translation $\nu$. 

Now, choosing the complete set of representatives ${\cal P}_0:=\{(0,i)\}_{i\in {\cal P}}$ of the $\langle \nu\rangle$-orbits in $\hat{\cal P}$, we have $Z=\hat{Z}_{0,0}=\bigoplus_{u,v\in \hat{\cal P}_0}e_u\hat{Z}e_v$. More generally, we have
$$\hat{Z}_{s,t}=\bigoplus_{i,j\in {\cal P}}e_{(s,j)}\hat{Z}e_{(t,i)}=
\bigoplus_{v\in \nu^s({\cal P}_0),u\in \nu^t({\cal P}_0)}e_v\hat{Z}e_u,$$
as in (\ref{R: las rho's de hat(Z) con accion}). Let us identify the $k$-algebra with idempotents $S=\bigoplus_{i\in {\cal P}}ke_i$ with the  $k$-subalgebra with idempotents $\bigoplus_{u\in {\cal P}_0}ke_u$ of $\hat{S}$, mapping each idempotent $e_i$ onto $e_{(0,i)}$. Then, 
we have have that $\rho_{s,t}(a)=\nu^{-s}a\nu^t$, for $a\in \hat{Z}_{s,t}$, thus we get the same maps $\rho_{s,t}$ considered in (\ref{R: las rho's de hat(Z) con accion}). 

 Then, from (\ref{L: hat(Z) nu-bimod y Z(P) b-alg --> hat(Z) b-nu-alg con seccion Z}), we obtain a $(b,\nu)$-algebra 
  $\hat{Z}=(\hat{Z},\{\hat{b}_n\}_{n\in \hueca{N}})$, over the elementary $k$-algebra with enough idempotents $\hat{S}=(\hat{S},\{e_u\}_{u\in \hat{\cal P}})$, with section  $(Z,\{b_n\}_{n\in \hueca{N}})$.    
\end{proof}

\begin{lemma}\label{L: unidades estrictas de hat(Z) son conjugadas}
 Assume that $\hat{Z}=(\hat{Z},\{\hat{b}_n\}_{n\in \hueca{N}})$ is a unitary strict $(b,\nu)$-algebra with strict units $\{\frak{e}_u\}_{u\in \hat{\cal P}}$, over an elementary $k$-algebra with enough idempotents $\hat{S}=(\hat{S},\{e_u\}_{u\in \hat{\cal P}})$. Then, the  strict units of $\hat{Z}$ 
   satisfy $\nu^{s}\frak{e}_u\nu^{-s}=\frak{e}_{\nu^s(u)}$, for all $u\in \hat{\cal P}$ and $s\in \hueca{Z}$.
\end{lemma}

\begin{proof} Choose a set of representatives ${\cal P}$ of $\hat{\cal P}$ as in (\ref{R: las rho's de hat(Z) con accion}). 
 For $u\in \hat{\cal P}$, we have  $u=\nu^s(v)$ for some $v\in {\cal P}$ and $s\in \hueca{Z}$, then $\frak{e}_u=e_u\frak{e}_ue_u\in \hat{Z}_{s,s}$. From the definition of $(b,\nu)$-algebra we obtain that $\nu\frak{e}_u\nu^{-1}\in e_{\nu(u)}\hat{Z}e_{\nu(u)}$. Hence, using that $\hat{Z}$ is unitary strict and (\ref{C: circ vs nu y nu-1}), we get 
 $\nu\frak{e}_u\nu^{-1}
=
\nu\frak{e}_u\nu^{-1}\circ \frak{e}_{\nu(u)}=\nu \frak{e}_u\circ \nu^{-1}\frak{e}_{\nu(u)}
=\nu(\frak{e}_u\circ \nu^{-1}\frak{e}_{\nu(u)})=\nu(\nu^{-1}\frak{e}_{\nu(u)})=\frak{e}_{\nu(u)}.$
It follows that $\nu^s\frak{e}_u\nu^{-s}=\frak{e}_{\nu^s(u)}$, for all integer $s$
\end{proof}

\begin{lemma}\label{L: hat(Z) es (b,sigma)-algebra estricta si Z lo es}
Assume that $\hat{Z}=(\hat{Z},\{\hat{b}_n\}_{n\in \hueca{N}})$ is a $(b,\nu)$-algebra, over an elementary $k$-algebra with enough idempotents $\hat{S}=(\hat{S},\{e_u\}_{u\in \hat{\cal P}})$, and consider its restriction 
$(Z,\{b_n\}_{n\in \hueca{N}})$, over an elementary $k$-algebra  with enough idempotents $S=(S,\{e_v\}_{v\in {\cal P}})$, as in (\ref{R: restriccion de (b,sigma)-algebra es b-algebra}). 

Then, if $\hat{Z}$ is a unitary
   strict $(b,\nu)$-algebra with strict units $\{\frak{e}_u\}_{u\in \hat{\cal P}}$,  the $b$-algebra $Z$ 
 is a unitary strict $b$-algebra with strict units $\{\frak{e}_u\}_{u\in {\cal P}}$.  

Conversely, if $Z$ is a unitary strict $b$-algebra with strict units $\{\frak{e}_v\}_{v\in {\cal P}}$.
Then, $\hat{Z}$ is naturally a unitary strict $(b,\nu)$-algebra with strict units $\{\frak{e}_u\}_{u\in \hat{\cal P}}$, where if $u\in \hat{\cal P}$, so $u=\nu^s(v)$, for some $v\in {\cal P}$, we have 
$\frak{e}_u=\frak{e}_{\nu^s(v)}:=\nu^s\frak{e}_v\nu^{-s}\in \hat{Z}_{s,s}.$

Moreover, the elements $\nu^s\frak{e}_u\nu^t\in \hat{Z}$ are strict
for all $u\in \hat{\cal P}$ and $s,t\in \hueca{Z}$. 
\end{lemma}

\begin{proof} For $v\in {\cal P}\subseteq \hat{\cal P}$, we have $\frak{e}_v=e_v\frak{e}_ve_v\in \bigoplus_{u,v\in {\cal P}}e_v\hat{Z}e_u=Z_{0,0}=Z$. So, the first claim of this lemma is clear because $b_n$ is the restriction of $\hat{b}_n$ to $Z^{\otimes n}$.

 For $u=\nu^s(v)$ with $v\in {\cal P}$, $\frak{e}_u=\rho_{s,s}^{-1}(\frak{e}_v)\in \hat{Z}_{s,s}$, and $a_2\in \hat{Z}_{s,t}$, we get 
 $$\begin{matrix}
 \hat{b}_2(\frak{e}_u\otimes a_2)
 &=&
 \rho_{s,t}^{-1}b_2(\rho_{s,s}\rho_{s,s}^{-1}(\frak{e}_v)\otimes\rho_{s,t}(a_2))\hfill\\
 &=&
 \rho_{s,t}^{-1}(e_v\rho_{s,t}(a_2))=e_{\nu^s(v)}a_2=e_ua_2.\hfill\\
 \end{matrix}$$
 Similarly, for $u=\nu^t(v)$ with $v\in {\cal P}$, $\frak{e}_u=\rho_{t,t}^{-1}(\frak{e}_v)\in \hat{Z}_{t,t}$, and $a_1\in \hat{Z}_{s,t}$, we have 
 $$\begin{matrix}\hat{b}_2(a_1\otimes \frak{e}_u)
 &=&
 (-1)^{s-t}\rho_{s,t}^{-1}b_2(\rho_{s,t}(a_1)\otimes\rho_{t,t}(\frak{e}_u))\hfill\\ &=&
 (-1)^{s-t}\rho_{s,t}^{-1}b_2(\rho_{s,t}(a_1)\otimes \frak{e}_v)=
 (-1)^{\vert a_1\vert+1}a_1e_{\nu^t(v)}\hfill\\
 &=&
 (-1)^{\vert a_1\vert+1}a_1e_u\hfill.\\
 \end{matrix}$$
 
 The fact that each element $\nu^s\frak{e}_v\nu^t$, with $v\in {\cal P}$, is strict follows  from the description of  $\hat{b}_n$ in terms of $b_n$ and the fact that $\frak{e}_v$ is strict in $Z$. In particular, the new elements $\frak{e}_u=\nu^s\frak{e}_v\nu^{-s}\in \hat{Z}$ are strict. So, $(\hat{Z}, \{\hat{b}_n\}_{n\in \hueca{N}})$ is a unitary strict $(b,\nu)$-algebra where all the elements $\nu^s\frak{e}_u\nu^t$, with $s,t\in \hueca{Z}$, are strict.  
\end{proof}

\begin{remark} Given a unitary strict  $b$-algebra  $Z=(Z,\{b_n\}_{n\in \hueca{N}})$, we are interested in a $(b,\nu)$-algebra  $\hat{Z}=(\hat{Z},\{\hat{b}_n\}_{n\in \hueca{N}})$,  which is unitary strict and has section $Z$, and their interaction. So from now on, until the end of this article,  we use freely the notations of   (\ref{D: (b,sigma)-algebra}) and the properties given in (\ref{L: hat(Z) es (b,sigma)-algebra estricta si Z lo es}).  
\end{remark}

\begin{lemma}\label{L: compos de a con trasladados de eu's} 
For each $a\in \hat{Z}_{s,t}$ and $u\in \hat{\cal P}$, we have:
\begin{enumerate}
 \item{\hskip.3cm}  $a\circ \frak{e}_u\nu^{-1}=(ae_u)\nu^{-1}$ \hbox{ \ and \  } $a\circ(\nu \frak{e}_u)=(a\nu)e_{u}$.
 \item{\hskip.3cm}  $\nu \frak{e}_u\circ a=(-1)^{{s-t}-1}\nu(e_ua)$ \hbox{ \ and \  } $(\frak{e}_u\nu^{-1})\circ a=(-1)^{s-t+1}e_u(\nu^{-1}a).$
\end{enumerate}
\end{lemma}

\begin{proof} 
It follows from the formulas in (\ref{C: circ vs nu y nu-1}).
\end{proof}

\begin{definition}\label{D:  tau(eu) y sigma(eu)}
For any $\frak{e}_u\in \hat{Z}_{s,s}$, 
 consider the following directed elements: 
$$
   \sigma(\frak{e}_u):=(-1)^s\nu \frak{e}_u\in \hat{Z}_{s+1,s}\hbox{ \ \ and \   \ }
   \tau(\frak{e}_u):=(-1)^s\frak{e}_u\nu^{-1}\in \hat{Z}_{s,s+1}.$$
Hence, we have 
$\frak{e}_u\in e_u\hat{Z}e_u, { \ } \sigma(\frak{e}_u)\in e_{\nu(u)}\hat{Z}e_u \hbox{  \ and \  }\tau(\frak{e}_u)\in e_{u}\hat{Z}e_{\nu(u)}$. These elements are homogeneous with  degrees 
$\vert \frak{e}_u\vert=-1$, $\vert \sigma(\frak{e}_u)\vert=-2$ and 
$\vert \tau(\frak{e}_u)\vert=0$. 
\end{definition}

With this notation, the preceding lemma (\ref{L: compos de a con trasladados de eu's}) yields the following. 

\begin{lemma}\label{L: aofrak y frakoa} 
For $a\in \hat{Z}_{s,t}$ and $u\in \hat{\cal P}$, we have  
$$\begin{matrix}
(1) &\hbox{ \ \ } 
 a\circ \sigma(\frak{e}_u)&=&(-1)^{t-1}(a\nu)e_u\hfill\\
 &\\
(2) &\hbox{ \ \ } 
  \sigma(\frak{e}_u)\circ a&=&(-1)^{t-1}\nu(e_ua)\\
 &\\
(3) &\hbox{ \ \ } 
  a\circ \tau(\frak{e}_u)&=&(-1)^t(ae_u)\nu^{-1} \\
  &\\
(4) &\hbox{ \ \ } 
 \tau(\frak{e}_u)\circ a &=&(-1)^te_u(\nu^{-1}a).\hfill \\
\end{matrix}$$
\end{lemma}

\begin{lemma}\label{L: tau(frak(eu)) o sigma(frak(eu)} For each $u\in \hat{\cal P}$, we have
$$
  \tau(\frak{e}_u)\circ\sigma(\frak{e}_u)=\frak{e}_u \hbox{ \ \ \ and  \ \ \  }
  \sigma(\frak{e}_u)\circ \tau(\frak{e}_u)=\frak{e}_{\nu(u)}. 
$$
\end{lemma}

\begin{proof} Assume that $\frak{e}_u\in \hat{Z}_{s,s}$, thus $\tau(\frak{e}_u)\in \hat{Z}_{s,s+1}$. From (\ref{L: aofrak y frakoa})(1), we have 
$\tau(\frak{e}_u)\circ\sigma(\frak{e}_u)=
(-1)^s(\tau(\frak{e}_u)\nu)e_u=\frak{e}_ue_u=\frak{e}_u$. 
From (\ref{L: aofrak y frakoa})(2) and (\ref{L: hat(Z) es (b,sigma)-algebra estricta si Z lo es})(1), we obtain 
$\sigma(\frak{e}_u)\circ \tau(\frak{e}_u)=
(-1)^s\nu(e_u\tau(\frak{e}_u))=\nu(e_u(\frak{e}_u\nu^{-1}))=e_{\nu(u)}\nu\frak{e}_u\nu^{-1}=e_{\nu(u)}\frak{e}_{\nu(u)}=\frak{e}_{\nu(u)}$. 
\end{proof}

\begin{remark}\label{R: sigma(fu)[-1] y tau(fu)[-1]}
For $u\in\hat{\cal P}$, we have 
$$\sigma(\frak{e}_u)[-1]=-\sigma(\frak{e}_{\nu^{-1}(u)}) \hbox{ \ and \ }  \tau(\frak{e}_u)[-1]=-\tau(\frak{e}_{\nu^{-1}(u)}).$$
Indeed, if $\frak{e}_u\in \hat{Z}_{s,s}$, hence 
$\frak{e}_{\nu^{-1}(u)}\in \hat{Z}_{s-1,s-1}$,  we have   
$
  \sigma(\frak{e}_u)[-1]=\nu^{-1}\sigma(\frak{e}_u)\nu
  =(-1)^s\frak{e}_u\nu=(-1)^s\nu \frak{e}_{\nu^{-1}(u)}=-\sigma(\frak{e}_{\nu^{-1}(u)}).
  $
We have used that $\frak{e}_u\nu=\nu \frak{e}_{\nu^{-1}(u)}$, which follows from (\ref{L: hat(Z) es (b,sigma)-algebra estricta si Z lo es}).
The other equality is verified similarly:
$\tau(\frak{e}_u)[-1]=\nu^{-1}\tau(\frak{e}_u)\nu=(-1)^s\nu^{-1}\frak{e}_u=(-1)^s\frak{e}_{\nu^{-1}(u)}\nu^{-1}=-\tau(\frak{e}_{\nu^{-1}(u)}).$
\end{remark}

\begin{lemma}\label{L: compos de tres en hat(B)}
For any element $a\in \hat{Z}$ and $u,v\in \hat{\cal P}$, we have:
$$\begin{matrix}
   (1)&\hbox{\hskip1cm}&\tau(\frak{e}_u)
   \circ(\sigma(\frak{e}_v)\circ a) =-e_ue_va\hfill\\
   \,\\
   (2)&\hbox{\hskip1cm}& \sigma(\frak{e}_u)\circ(\tau(\frak{e}_v)\circ a)=-e_{\nu(u)}e_{\nu(v)}a\hfill\\
   \end{matrix}\hbox{ \   }$$
 $$\begin{matrix}
    \hbox{\hskip.5cm}(3)&\hbox{\hskip1cm}& \sigma(\frak{e}_u)\circ(a\circ \tau(\frak{e}_v)) =e_{\nu(u)}a[1]e_{\nu(v)}\hfill\\
     \,\\
     \hbox{\hskip.5cm}(4)&\hbox{\hskip1cm}& (\sigma(\frak{e}_u)\circ  a)\circ \tau(\frak{e}_v)=-e_{\nu(u)}a[1]e_{\nu(v)}\hfill\\
    \,\\
     \hbox{\hskip.5cm}(5)&\hbox{\hskip1cm}&(a\circ\tau(\frak{e}_u))\circ \sigma(\frak{e}_v)=ae_ue_v.\hfill\\
   \end{matrix}\hbox{\ \ }$$  
\end{lemma}

\begin{proof} We may assume that $a\in \hat{Z}_{s,t}$. Then, from (\ref{L: aofrak y frakoa})(2)\&(4), we have 
 $$\begin{matrix}
 \tau(\frak{e}_u)\circ(\sigma(\frak{e}_v)\circ a)
   &=&
   (-1)^{t-1}\tau(\frak{e}_u)\circ \nu(e_va)\hfill\\
   &=&
   (-1)^{t-1}(-1)^te_u(\nu^{-1}\nu(e_va))=-e_ue_va.\\
   \end{matrix}$$
 From (\ref{L: aofrak y frakoa})(4)\&(2), we have 
 $$\begin{matrix}
   \sigma(\frak{e}_u)\circ(\tau(\frak{e}_v)\circ a)
   &=&
   (-1)^t\sigma(\frak{e}_u)\circ e_v(\nu^{-1}a)=
   (-1)^t(-1)^{t-1}\nu(e_ue_v(\nu^{-1}a))\hfill\\
   &=&
   -\nu(e_u)\nu(e_v)a=-e_{\nu(u)}e_{\nu(v)}a.\hfill\\
   \end{matrix}$$
  The verification of (3), (4), and (5) is similar, we use  (\ref{L: aofrak y frakoa})(3)\&(2), (\ref{L: aofrak y frakoa})(2)\&(3), and (\ref{L: aofrak y frakoa})(3)\&(1), 
respectively. 
\end{proof}

\begin{lemma}\label{L: bn(sigmala1...)=-sigmalbn(a1...)}
For any sequence $a_1\in \hat{Z}_{s_0,s_1},\ldots,a_n\in \hat{Z}_{s_{n-1},s_n}$ and $v\in \hat{\cal P}$, we have
$\hat{b}_n(\sigma(\frak{e}_v)\circ a_1\otimes a_2\otimes\cdots \otimes a_n)=
-\sigma(\frak{e}_v)\circ\hat{b}_n(a_1\otimes a_2\otimes\cdots \otimes a_n)$. 
\end{lemma}

\begin{proof} We may assume that the elements $a_1,\ldots,a_n$ are homogeneous. 
 From (\ref{L: aofrak y frakoa})(2) and (\ref{L: nu por hat(b)n}), we have
$$\begin{matrix}
  \hat{b}_n(\sigma(\frak{e}_v)\circ a_1\otimes a_2\otimes\cdots \otimes a_n)&=&(-1)^{s_1-1}  \hat{b}_n(\nu(e_v a_1)\otimes a_2\otimes\cdots \otimes a_n)\hfill\\
  &=&
  (-1)^{s_n}  \nu\hat{b}_n(e_v a_1\otimes a_2\otimes\cdots \otimes a_n).\hfill\\
  &=&
  -(-1)^{s_n-1}\nu[e_v\hat{b}_n(a_1\otimes a_2\otimes\cdots \otimes a_n)]\hfill\\
  &=&
  -\sigma(\frak{e}_v)\circ\hat{b}_n(a_1\otimes a_2\otimes\cdots \otimes a_n).\hfill\\
  \end{matrix}$$
\end{proof}

\begin{lemma}\label{L: bn(...al tau(fu) otimes sigma(fv) al+1...)=}
For any $n\geq 2$ and any sequence $a_1\in \hat{Z}_{s_0,s_1},\ldots,a_n\in \hat{Z}_{s_{n-1},s_n}$ of homogeneous elements and $l\in[1,n-1]$ and $u,v\in \hat{\cal P}$, we have that 
$$\hat{b}_n(a_1\otimes a_2\otimes\cdots\otimes a_{l-1}\otimes (a_l\circ\tau(\frak{e}_u))\otimes(\sigma(\frak{e}_v)\circ a_{l+1})\otimes a_{l+2}\otimes\cdots \otimes a_n)$$
coincides with
$(-1)^{\vert a_l\vert}\hat{b}_n(a_1\otimes a_2\otimes\cdots \otimes a_le_u\otimes e_va_{l+1}\otimes \cdots \otimes a_n)$.
\end{lemma}

\begin{proof} Denote by $\Delta$ the first expression in the statement of this lemma. From (\ref{L: aofrak y frakoa})(2-3) and (\ref{R: hat(b)n por conjugados de nu internos}) we have 
$$\begin{matrix}
  \Delta&=& 
  (-1)^{s_l}(-1)^{s_{l+1}-1}\hat{b}_n(a_1\otimes\cdots\otimes (a_le_u)\nu^{-1}\otimes \nu (e_va_{l+1})\otimes \cdots \otimes a_n)\hfill\\
  &=&
  (-1)^{\vert a_l\vert}\hat{b}_n(a_1\otimes a_2\otimes\cdots \otimes a_le_u\otimes e_va_{l+1}\otimes \cdots \otimes a_n).\hfill\\
  \end{matrix}$$
\end{proof}

\section{The $b$-category $\ad(\hat{Z})$}\label{S: la b-categoria ad(hat(Z))}

In this section we state some basic properties of the $b$-category $\ad(\hat{Z})$ associated to a unitary strict $(b,\nu)$-algebra $\hat{Z}$, see  (\ref{R: recordatorio de ad(Z)}).  We keep the notations of the last section. 
So, the objects of $\ad(\hat{Z})$ are the right support-finite $\hat{S}$-modules; given two such objects $X$ and $Y$, the corresponding space of morphisms is 
  $$\ad(\hat{Z})(X,Y):=\bigoplus_{u,v\in \hat{\cal P}} 
  \Hom_k(Xe_u,Ye_v)\otimes_ke_v\hat{Z}e_u,$$
  with the canonical grading of the tensor product.  The maps $\hat{b}^{ad}_n$ are defined, for $n\in \hueca{N}$ and a sequence of objects $X_0,X_1,\ldots,X_n$  
  in $\ad(\hat{Z})$, on generators by 
  $$\begin{matrix}\ad(\hat{Z})(X_{n-1},X_n)\otimes_k\cdots \otimes_k\ad(\hat{Z})(X_0,X_1)\hbox{\,}\rightmap{\hat{b}^{ad}_n}\hbox{\,}\ad(\hat{Z})(X_0,X_n)\\
   (f_n\otimes a_n)\otimes\cdots\otimes (f_1\otimes a_1)\hbox{\ }\longmapsto{\ } 
   f_n\cdots f_2f_1\otimes \hat{b}_n(a_n\otimes\cdots\otimes a_1). 
    \end{matrix}$$
  
\begin{remark}\label{R: la base hat(hueca(B)) de hat(Z)}
 We fix a \emph{directed basis} $\hat{\hueca{B}}$ for the graded vector space $\hat{Z}$, as follows.  
 First, we consider a directed basis $\hueca{B}$ of $Z=\hat{Z}_{0,0}$ as in (\ref{N: la base de B}). Then, we define 
 $\hat{\hueca B}_{s,t}:=\nu^s\hueca{B}\nu^{-t}$, for all $s,t\in\hueca{Z}$. 
 Thus, $\hat{\hueca{B}}_{s,t}$ is a directed basis of $\hat{Z}_{s,t}$, and we consider  the directed basis  $\hat{\hueca{B}}=\biguplus_{s,t\in \hueca{Z}}\hat{\hueca{B}}_{s,t}$ of $\hat{Z}=\bigoplus_{s,t}\hat{Z}_{s,t}$.  
 
 If $\hueca{B}_q$ is the subset of $\hueca{B}$ formed by its homogeneous basis elements of degree $q$, which span the homogeneous component $Z_q$ of $Z$ of degree $q$, then $\nu^s\hueca{B}_q\nu^{-t}$ spans  
 the homogeneous component of $\hat{Z}_{s,t}$ of degree $q+t-s$. Notice that $\hat{\hueca{B}}$ contains the strict units of $\hat{Z}$, see (\ref{L: unidades estrictas de hat(Z) son conjugadas}). 
\end{remark}

 \begin{definition}\label{D: X[1] en para X en ad(hat(B))}
For any object $X$ of $\ad(\hat{Z})$, we define $X[1]$ as the right $\hat{S}$-module obtained from $X$ by restriction of scalars through the automorphism $\nu^{-1}:\hat{S}\rightmap{}\hat{S}$. That is, by definition, the underlying group of $X[1]$ is the same  $X$ and each idempotent  $e_u$ acts on any element $x\in X[1]$ by the rule $x\cdot e_u:=xe_{\nu^{-1}(u)}$. In the following few lines, we keep using the notation $x\cdot s$
for the action of the element $s$ of $\hat{S}$ on an element $x$ in $X[1]$, while $xs$ denotes the action of the same $s$ on $x$ in $X$. 

We consider the linear map $\phi(X):X\rightmap{}X[1]$ given by the identity map. 
Then, we have $\phi(X)[xe_u]=xe_u=x\cdot e_{\nu(u)}$, for $x\in X$ and $u\in \hat{\cal P}$.  Therefore, we have the corresponding linear restriction $\phi(X)_u:Xe_u\rightmap{}X[1]e_{\nu(u)}$ of $\phi(X)$.
\end{definition}

\begin{remark}\label{R: X[-1] y psi(X)} With the preceding notation, we define the right $\hat{S}$-module $X[-1]$ as the right $\hat{S}$ module obtained from $X$ by restriction using the automorphism $\nu:\hat{S}\rightmap{}\hat{S}$, and we have the linear map $\psi(X):X\rightmap{}X[-1]$ given by the identity, which  induces linear  restrictions 
$\psi(X)_{u}:Xe_u\rightmap{}X[-1]e_{\nu^{-1}(u)}$. 
Clearly, we have the equality of $\hat{S}$-modules $X[1][-1]=X=X[-1][1]$. Moreover, we have that the following composition is the identity map 
$$Xe_u\rightmap{ \  \ \phi(X)_u \ \  }X[1]e_{\nu(u)}\rightmap{ \ \ \psi(X[1])_{\nu(u)} \ \ }X[1][-1]e_{\nu^{-1}\nu(u)}=  Xe_u.$$
Thus, we have $\phi(X)_u^{-1}=\psi(X[1])_{\nu(u)}$. 
 \end{remark}

\begin{definition}\label{D: trasladados en ad(hat(B))}
Given $f=\sum_{a}f_a\otimes a\in \ad(\hat{Z})(X,Y)$,   define 
$$f[1]:=\sum_{a}\phi(Y)_{v(a)}f_a\phi(X)^{-1}_{u(a)}\otimes a[1] $$
and
$$f[-1]:=\sum_{a}\psi(Y)_{v(a)}f_a\psi(X)^{-1}_{u(a)}\otimes a[-1] $$
so we have that each 
$\phi(Y)_{v(a)}f_a\phi(X)^{-1}_{u(a)}\in 
\Hom_k(X[1]e_{\nu(u(a))},Y[1]e_{\nu(v(a))})$, with $a[1]\in e_{\nu(v(a))}\hat{Z}e_{\nu(u(a))}$. So   $f[1]\in \ad(\hat{Z})(X[1],Y[1])$ and, similarly, $f[-1]\in \ad(\hat{Z})(X[-1],Y[-1])$.
Our choice of directed basis $\hat{\hueca{B}}$ of $\hat{Z}$ in (\ref{R: la base hat(hueca(B)) de hat(Z)}) guarantees that  the expressions of $f[1]$ and $f[-1]$ are given in terms of basis elements. If $f$ is homogeneous, so are $f[1]$ and $f[-1]$, with $\vert f[1]\vert=\vert f\vert=\vert f[-1]\vert$.
\end{definition}

\begin{remark}\label{R: f[1][-1]=f=f[-1][1]}
 For any $f\in \ad(\hat{Z})(X,Y)$, we have $f[1][-1]=f=f[-1][1]$. Indeed, if $f=\sum_{a}f_a\otimes a$, from (\ref{R: X[-1] y psi(X)}), we have 
 $$f[1][-1]=\sum_{a}\psi(Y[1])_{\nu(v(a))}\phi(Y)_{v(a)}f_a\phi(X)^{-1}_{u(a)}\psi(X[1])^{-1}_{\nu(u(a))}\otimes a[1][-1]=f$$
\end{remark}
 
\begin{lemma}\label{L: translacion automorfismo de la b-cat ad(hat(B))}
 Let $X_0\rightmap{f_1}X_1,...,X_{n-1}\rightmap{f_n}X_n$ be any sequence of morphisms in $\ad(\hat{Z})$, then  
 $\hat{b}^{ad}_n(f_n[1]\otimes \cdots \otimes f_2[1]\otimes f_1[1])=\hat{b}^{ad}_n(f_n\otimes \cdots \otimes f_2\otimes f_1)[1]$. 
\end{lemma}

\begin{proof} It is enough to show this equality for morphisms $f_1,\ldots,f_n$ of the form $f_i=h_i\otimes a_i$, where $a_i\in \hat{Z}$ are directed  and $h_i\in\Hom_k(X_{i-1}e_{u(a_i)},X_ie_{v(a_i)})$. 
  Notice that $\hat{b}^{ad}_n((h_n\otimes a_n)[1]\otimes \cdots\otimes (h_2\otimes a_2)[1]\otimes (h_1\otimes a_1)[1])$ 
coincides with 
  $[\phi(X_n)_{v(a_n)}h_n\cdots h_2h_1\phi(X_0)^{-1}_{u(a_1)}]\otimes \hat{b}_n(a_n[1]\otimes\cdots\otimes a_2[1]\otimes a_1[1])$. 
   From (\ref{R: hat(b)n de tensores de trasladados}),  
   the last expression coincides with  
  $\hat{b}^{ad}_n((h_n\otimes a_n) \cdots \otimes (h_2\otimes a_2)\otimes (h_1\otimes a_1))[1]$. 
\end{proof}

\begin{definition}\label{D:  sigmaX y tauX}
 Let $X$ be an object of $\ad(\hat{Z})$. We have 
$\phi(X)_{u}\otimes \sigma(\frak{e}_u)\in  \Hom_k(Xe_u,X[1]e_{\nu(u)})\otimes_ke_{\nu(u)}\hat{Z}e_u$, for each  $u\in \hat{\cal P}$. 
Define 
$$\sigma_X:=\sum_{u\in \hat{\cal P}}\phi(X)_{u}\otimes \sigma(\frak{e}_u)\in  \ad(\hat{Z})(X,X[1]),$$
the sum is finite because $Xe_{u}=0$, for almost all $u\in \hat{\cal P}$. 
The morphism $\sigma_X$ has degree $\vert\sigma_X\vert=\vert\sigma(\frak{e}_u)\vert=-2$.  

Similarly,   we have 
$\phi(X)^{-1}_u\otimes \tau(\frak{e}_u)\in  \Hom_k(X[1]e_{\nu(u)},Xe_u)\otimes_k
e_u\hat{Z}e_{\nu(u)}$, for each $u\in \hat{\cal P}$,
and we can define 
$$\tau_X:=\sum_{u\in \hat{\cal P}}\phi(X)_{u}^{-1}\otimes \tau(\frak{e}_u)\in  \ad(\hat{Z})(X[1],X).$$
The morphism $\tau_X$ has degree $\vert\tau_X\vert=\vert\tau(\frak{e}_u)\vert=0$.

Notice that the morphisms $\tau_X$ and $\sigma_X$ are strict, for any object $X$ of $\ad(\hat{Z})$. This follows from  (\ref{L: hat(Z) es (b,sigma)-algebra estricta si Z lo es}) and (\ref{D:  tau(eu) y sigma(eu)}).
\end{definition}

\begin{lemma}\label{L: tauX circ sigmaX}
 For any object $X$ of $\ad(\hat{Z})$, we have:
$$ \tau_X\circ \sigma_X=\hueca{I}_X \hbox{ \ \ and \ \ } 
\sigma_X\circ \tau_X=\hueca{I}_{X[1]}.$$
\end{lemma}

\begin{proof} From (\ref{L: tau(frak(eu)) o sigma(frak(eu)}), we obtain:
$\tau_X\circ \sigma_X=
   \sum_u\phi(X)^{-1}_u\phi(X)_u\otimes \tau(\frak{e}_u)\circ\sigma(\frak{e}_u)
   =
   \sum_uid_{Xe_u}\otimes \frak{e}_u=\hueca{I}_X$, 
and 
$ \sigma_X\circ \tau_X=
   \sum_u\phi(X)_u\phi(X)^{-1}_u\otimes \sigma(\frak{e}_u)\circ \tau(\frak{e}_u)
   =
   \sum_uid_{X[1]e_{\nu(u)}}\otimes \frak{e}_{\nu(u)}=\hueca{I}_{X[1]}$. 
\end{proof}

\begin{remark}\label{R: sigmaX[-1] y tauX[-1]}
For any object $X$ of $\ad(\hat{Z})$, we have 
$$\sigma_X[-1]=-\sigma_{X[-1]} \hbox{ \ and \ } \tau_X[-1]=-\tau_{X[-1]}.$$
Indeed, from (\ref{R: sigma(fu)[-1] y tau(fu)[-1]}), we have 
$$\begin{matrix}\sigma_X[-1]&=&
 [\sum_u\phi(X)_u\otimes \sigma(\frak{e}_u)][-1]\hfill\\
 &=&
 \sum_u\psi(X[1])_{\nu(u)}\phi(X)_u\psi(X)^{-1}_u\otimes \sigma(\frak{e}_u)[-1]\hfill\\
 &=&
 -\sum_u\psi(X)^{-1}_u\otimes \sigma(\frak{e}_{\nu^{-1}(u)})\hfill\\
 &=&
 -\sum_u\phi(X[-1])_{\nu^{-1}(u)}\otimes \sigma(\frak{e}_{\nu^{-1}(u)})
 =-\sigma_{X[-1]}.\hfill\\
  \end{matrix}$$
  The other equality is verified similarly. 
\end{remark}

\begin{lemma}\label{L: circ de tres sigmas y taus con morfismos}
The following holds:
\begin{enumerate}
 \item For any morphism $f:X\rightmap{}Y$ in $\ad(\hat{Z})$, we have  
 $$(f\circ \tau_X)\circ \sigma_X=f\hbox{ \ \ and \ \ } \tau_Y\circ(\sigma_Y\circ f)=-f.$$
 \item For any morphism $g:X\rightmap{}Y[1]$ in $\ad(\hat{Z})$, we have  
 $$\sigma_Y\circ(\tau_Y\circ g)=-g.$$
 \item For any morphism $f:X\rightmap{}Y$ in $\ad(\hat{Z})$, we have  
 $$\sigma_Y\circ( f\circ \tau_X)=f[1]\hbox{ \ \ and \ \ }
 (\sigma_Y\circ f)\circ \tau_X=-f[1].$$
\end{enumerate}
\end{lemma}

\begin{proof} From (\ref{L: compos de tres en hat(B)})(5), we have 
 $$\begin{matrix}
   (f\circ \tau_X)\circ \sigma_X
   &=&
   (\sum_{a,u}f_a\phi(X)^{-1}_u\otimes (a\circ \tau(\frak{e}_u)))\circ \sigma_X\hfill\\
   &=&
   \sum_{a,u,v}f_a\phi(X)^{-1}_u\phi(X)_v\otimes ((a\circ \tau(\frak{e}_u))\circ \sigma(\frak{e}_v))\hfill\\
   &=&
   \sum_{a,u,v}f_a\phi(X)^{-1}_u\phi(X)_v\otimes ae_ue_v\hfill\\
   &=&
   \sum_af_a\otimes a=f.\hfill\\
  \end{matrix}$$
From (\ref{L: compos de tres en hat(B)})(1), we have 
$$\begin{matrix}
   \tau_Y\circ(\sigma_Y\circ f)
   &=&
   \tau_Y\circ(\sum_{a,v}\phi(Y)_vf_a\otimes (\sigma(\frak{e}_v)\circ a))\hfill\\
   &=&
 \sum_{a,v,u}\phi(Y)^{-1}_u\phi(Y)_vf_a\otimes( \tau(\frak{e}_u)\circ(\sigma(\frak{e}_v)\circ a))\hfill\\
 &=&
 -\sum_{a,v,u}\phi(Y)^{-1}_u\phi(Y)_vf_a\otimes e_ue_va\hfill\\
 &=&
 -\sum_af_a\otimes a=-f.\hfill\\
  \end{matrix}$$
 The verification of (2) and (3) is similar, now using consecutively (\ref{L: compos de tres en hat(B)})(2), (\ref{L: compos de tres en hat(B)})(3), and (\ref{L: compos de tres en hat(B)})(4). 
\end{proof}

\begin{lemma}\label{L: bn(sigmaXncirc fn...)=sigmaXn(bn fn...)}
 Let  $X_0\rightmap{f_1}X_1,...,X_{n-1}\rightmap{f_n}X_n$ be any sequence of morphisms in $\ad(\hat{Z})$. Then,   we have 
$\hat{b}_n^{ad}(\sigma_{X_n}\circ f_n\otimes\cdots\otimes f_1)=-\sigma_{X_n}\circ\hat{b}_n^{ad}( f_n\otimes\cdots\otimes f_1)$.  
\end{lemma}

\begin{proof} It is enough to show this equality for morphisms $f_1,\ldots,f_n$ of the form $f_i=h_i\otimes a_i$, where $a_i\in \hat{Z}$ is directed and 
$h_i\in \Hom_k(X_{i-1}e_{u(a_i)},X_ie_{v(a_i)})$.

We have $\sigma_{X_n}\circ f_n=\phi(X_n)_{v(a_n)}h_n\otimes (\sigma(\frak{e}_{v(a_n)})\circ a_n)$. Then, from (\ref{L: bn(sigmala1...)=-sigmalbn(a1...)}), if we set $\Delta:=\hat{b}_n^{ad}(\sigma_{X_n}\circ f_n\otimes\cdots\otimes f_1)$, we have:
$$\begin{matrix}
  \Delta&=& \hat{b}_n^{ad}([\phi(X_n)_{v(a_n)}h_n\otimes (\sigma(\frak{e}_{v(a_n)})\circ a_n)]\otimes (h_{n-1}\otimes a_n)\otimes\cdots\otimes (h_1\otimes a_1))\hfill\\
  &=&
  \phi(X_n)_{v(a_n)}h_nh_{n-1}\cdots h_1\otimes \hat{b}_n[(\sigma(\frak{e}_{v(a_n)})\circ a_n)\otimes a_{n-1}\otimes \cdots\otimes a_1]\hfill\\
  &=&
   -\phi(X_n)_{v(a_n)}h_nh_{n-1}\cdots h_1\otimes (\sigma(\frak{e}_{v(a_n)})\circ\hat{b}_n[ a_n\otimes a_{n-1}\otimes \cdots\otimes a_1])\hfill\\
  &=&
   -\sigma_{X_n}\circ (h_nh_{n-1}\cdots h_1\otimes \hat{b}_n( a_n\otimes\cdots\otimes a_1))\hfill\\
  &=&
  -\sigma_{X_n}\circ\hat{b}_n^{ad}( f_n\otimes\cdots\otimes f_1).\hfill\\
  \end{matrix}$$
\end{proof}

\begin{lemma}\label{L: bn(...f_ntauX otimes sigmaXncirc fn...)=bn( fn...)}
 For $n\geq 2$, let $X_0\rightmap{f_1}X_1,...,X_{n-1}\rightmap{f_n}X_n$ be a sequence of homogeneous morphisms in $\ad(\hat{Z})$, take $l\in [1,n-1]$. Then,   
 $$\hat{b}^{ad}_n(f_n\otimes f_{n-1}\otimes\cdots\otimes (f_{l+1}\circ\tau_{X_l})
 \otimes (\sigma_{X_l}\circ f_l)\otimes\cdots \otimes f_2\otimes f_1)$$
coincides with
$(-1)^{\vert f_{l+1}\vert}\hat{b}^{ad}_n(f_n\otimes f_{n-1}\otimes\cdots \otimes f_{l+1}\otimes f_l\otimes \cdots \otimes f_2\otimes f_1)$. 
\end{lemma}

\begin{proof} It is enough to show this equality for homogeneous morphisms $f_1,\ldots,f_n$ of the form $f_i=h_i\otimes a_i$, where $a_i\in \hat{Z}$ is directed and 
$h_i\in \Hom_k(X_{i-1}e_{u(a_i)},X_ie_{v(a_i)})$.

We have 
$$\sigma_{X_l}\circ f_l=\phi(X_l)_{v(a_l)}h_l\otimes (\sigma(\frak{e}_{v(a_l)})\circ a_l)$$
and
$ f_{l+1}\circ\tau_{X_l}=h_{l+1}\phi(X_l)^{-1}_{u(a_{l+1})}\otimes (a_{l+1}\circ\tau(\frak{e}_{u(a_{l+1})}))$. 
Denote by $\Delta$ the first expression in the statement of this lemma. Then, from (\ref{L: bn(...al tau(fu) otimes sigma(fv) al+1...)=}), we have:  
$$\begin{matrix}
  \Delta&=&
  h_n\cdots h_1\otimes \hat{b}_n(a_n\otimes \cdots\otimes (a_{l+1}\circ\tau(\frak{e}_{u(a_{l+1})}))\otimes(\sigma(\frak{e}_{v(a_l)})\circ a_l)\otimes \cdots \otimes a_1)\hfill\\
  &=&
  (-1)^{\vert a_{l+1}\vert}h_n\cdots h_1\otimes \hat{b}_n(a_n\otimes \cdots\otimes a_{l+1}\otimes a_l\otimes \cdots \otimes a_1)\hfill\\
  &=&
  (-1)^{\vert f_{l+1}\vert}\hat{b}_n^{ad}((h_n\otimes a_n)\otimes\cdots\otimes (h_{l+1}\otimes a_{l+1})\otimes (h_l\otimes a_l)\otimes \cdots\otimes (h_1\otimes a_1) )\hfill\\
  &=&
  (-1)^{\vert f_{l+1}\vert}\hat{b}_n^{ad}(f_n\otimes\cdots\otimes f_{l+1}\otimes f_l\otimes \cdots\otimes f_1).\hfill\\
  \end{matrix}$$
\end{proof}

\begin{lemma}\label{L: como entra la sigma en medio de la hat(b)adn(... sigmaX...)}\label{L: reemplazo tensors delta por tensors de delta[1]}
 Let  $X_0\rightmap{f_1}X_1,...,X_{n-1}\rightmap{f_n}X_n$ be a sequence of homogeneous morphisms in $\ad(\hat{Z})$, take $l\in [1,n-1]$. Then,   
 $$\sigma_{X_n}\circ\hat{b}^{ad}_n(f_n\otimes f_{n-1}\otimes\cdots \otimes f_l\otimes\cdots\otimes f_2\otimes f_1)$$
coincides with
$$(-1)^{d_l}\hat{b}^{ad}_n(f_n[1]\otimes f_{n-1}[1]\otimes\cdots   \otimes f_{l+1}[1]
 \otimes (\sigma_{X_l}\circ f_l)\otimes f_{l-1}\otimes  \cdots \otimes f_2\otimes f_1),$$
 where $d_l=\vert f_{l+1}\vert+\cdots+\vert f_n\vert+1$. Equivalently, 
 $$\hat{b}^{ad}_n(f_n\otimes f_{n-1}\otimes\cdots \otimes f_l\otimes\cdots\otimes f_2\otimes f_1)$$
 coincides with 
 $$(-1)^{d_l+1}\tau_{X_n}\circ\hat{b}^{ad}_n(f_n[1]\otimes f_{n-1}[1]\otimes\cdots   \otimes f_{l+1}[1]
 \otimes (\sigma_{X_l}\circ f_l)\otimes f_{l-1}\otimes  \cdots \otimes f_2\otimes f_1).$$
\end{lemma}

\begin{proof} Recall that $\sigma_{X_{i+1}}\circ (f_i\circ \tau_{X_{i-1}})=f_i[1]$, for $i\in [1,n-1]$.  In the following, where we use repeatedly (\ref{L: bn(...f_ntauX otimes sigmaXncirc fn...)=bn( fn...)}), we set $s_i:=\vert f_{l+1}\vert+\cdots+\vert f_{l+i}\vert$, for $i\in [1,n-l]$, and
 $\Delta:=\hat{b}^{ad}_n(f_n\otimes f_{n-1}\otimes\cdots \otimes f_l\otimes\cdots\otimes f_2\otimes f_1)$. Then, 
$$\begin{matrix}
 \Delta&=&(-1)^{s_1}\hat{b}^{ad}_n(f_n\otimes f_{n-1}\otimes\cdots \otimes (f_{l+1}\circ \tau_{X_l})\otimes (\sigma_{X_l}\circ f_l)\otimes\cdots\otimes f_2\otimes f_1)
 \hfill\\
 &=&
 (-1)^{s_2}\hat{b}^{ad}_n(\cdots\otimes (f_{l+2}\circ \tau_{X_{l+1}}) \otimes (\sigma_{X_{l+1}}\circ (f_{l+1}\circ \tau_{X_l}))\otimes (\sigma_{X_l}\circ f_l)\otimes\cdots)
 \hfill\\
 &=&
 (-1)^{s_2}\hat{b}^{ad}_n(\cdots\otimes (f_{l+2}\circ \tau_{X_{l+1}}) \otimes 
 f_{l+1}[1]\otimes (\sigma_{X_l}\circ f_l)\otimes\cdots)
 \hfill\\
 &=&
 (-1)^{s_3}\hat{b}^{ad}_n(\cdots \otimes (f_{l+3}\circ \tau_{X_{l+2}})\otimes (\sigma_{X_{l+2}}\circ (f_{l+2}\circ \tau_{X_{l+1}}) )\otimes 
 f_{l+1}[1]\otimes (\sigma_{X_l}\circ f_l)\otimes\cdots)
 \hfill\\
 &=&
 (-1)^{s_3}\hat{b}^{ad}_n(\cdots \otimes (f_{l+3}\circ \tau_{X_{l+2}})\otimes f_{l+2}[1]\otimes 
 f_{l+1}[1]\otimes (\sigma_{X_l}\circ f_l)\otimes\cdots)
 \hfill\\
 &\cdots&\\
 &=&
 (-1)^{s_{n-l}}\hat{b}^{ad}_n((f_n\circ \tau_{X_{n-1}})\otimes f_{n-1}[1]\otimes\cdots \otimes f_{l+2}[1]\otimes 
 f_{l+1}[1]\otimes (\sigma_{X_l}\circ f_l)\otimes\cdots)
 \hfill\\
  \end{matrix}$$
Then, from (\ref{L: bn(sigmaXncirc fn...)=sigmaXn(bn fn...)}) and (\ref{L: circ de tres sigmas y taus con morfismos}), we get 
$$\begin{matrix}\sigma_{X_n}\circ \Delta&=&(-1)^{s_{n-l}+1}\hat{b}^{ad}_n(\sigma_{X_n}\circ(f_n\circ \tau_{X_{n-1}})\otimes f_{n-1}[1]\otimes\cdots \otimes  
 f_{l+1}[1]\otimes (\sigma_{X_l}\circ f_l)\otimes\cdots)\hfill\\
 &=&
 (-1)^{s_{n-l}+1}\hat{b}^{ad}_n(f_n[1]\otimes f_{n-1}[1]\otimes\cdots \otimes  
 f_{l+1}[1]\otimes (\sigma_{X_l}\circ f_l)\otimes\cdots).\hfill\\
 \end{matrix}$$
 The second part follows from the first one if we apply $\tau_{X_n}$ on both sides and use (\ref{L: circ de tres sigmas y taus con morfismos}). 
\end{proof}

\section{Conflations in ${\cal Z}(\hat{Z})$ and the functors $T$ and $J$}\label{Z(hatZ)}

Here, we keep the preceding terminology, where $\hat{Z}$ is a $(b,\nu)$-algebra over the elementary  algebra $\hat{S}$, with enough idempotents $\{e_u\}_{u\in \hat{\cal P}}$,  and we assume that it is unitary strict with  strict units $\{\frak{e}_u\}_{u\in \hat{\cal P}}$, as in (\ref{D: unitary strict}). We have the associated $b$-category  $\ad(\hat{Z})$ over $\hat{S}$, as in (\ref{R: recordatorio de ad(Z)}), and   
 a fixed basis $\hat{\hueca{B}}$ for the vector space $\hat{Z}$ formed by homogeneous directed elements, and containing the strict units of $\hat{Z}$. 
 
 Then, we have the $b$-category $\tw(\hat{Z})$ reminded in (\ref{R: record de tw(Z)}). Recall that,  given  two morphisms $f:(X,\delta_X)\rightmap{}(Y,\delta_Y)$ and $g:(Y,\delta_Y)\rightmap{}(W,\delta_W)$  in $\tw(\hat{Z})$, we use  the notation $g\star f=\hat{b}_2^{tw}(g\otimes f)$. Then, we have the precategory ${\cal Z}(\hat{Z})$ with composition $\star$ and we have at our disposal all the results on its conflations presented in sections \S\ref{special confl in Z(Z)} and \S\ref{Confl in Z(Z)}.

In the following, we investigate further the precategory ${\cal Z}(\hat{Z})$ and show that the analogy with exact categories remarked in  (\ref{R: comentario al summary for special conflations}), in this case, can be  extended to an analogy with special Frobenius categories, see \cite{H}, \cite{dgb}(8.6), \cite{DGCKeller}, and \cite{BMJ}\S3. We introduce a translation $T$ and a functor  $J$ on ${\cal Z}(\hat{Z})$, which associates  projective (and injective) objects relative to special conflations. The endofunctors $T$ and $J$ have similar properties to the corresponding functors on a special Frobenius category.

 \begin{lemma}\label{L: traslacion en Z(hat(B)) y en H(hat(B))}
  There is an autofunctor $T:{\cal Z}(\hat{Z})\rightmap{}{\cal Z}(\hat{Z})$, which induces an autofunctor $T:{\cal H}(\hat{Z})\rightmap{}{\cal H}(\hat{Z})$. Given $(X,\delta_X)\in {\cal Z}(\hat{Z})$, we have 
  $$T(X,\delta_X):=(X,\delta_X)[1]:=(X[1],\delta_X[1]),$$
  see (\ref{D: X[1] en para X en ad(hat(B))}) and (\ref{D: trasladados en ad(hat(B))}). Given $f\in {\cal Z}(\hat{Z})((X,\delta_X),(Y,\delta_Y))$,  by definition, 
  $$T(f):=f[1]:(X,\delta_X)[1]\rightmap{}(Y,\delta_Y)[1].$$
  Its inverse $T^{-1}$ is given by $(X,\delta_X)\mapsto (X[-1],\delta_X[-1])$ and $f\longmapsto f[-1]$. 
 \end{lemma}

 \begin{proof} Given $(X,\delta_X)\in {\cal Z}(\hat{Z})$, we know that $\delta_X\in \ad(\hat{Z})(X,X)_0$, and also $\delta_X[1]\in \ad(\hat{Z})(X[1],X[1])_0$. Moreover, from (\ref{L: translacion automorfismo de la b-cat ad(hat(B))}), we get  
 $\sum_s\hat{b}_s^{ad}((\delta_X[1])^{\otimes s})=\sum_s\hat{b}_s^{ad}(\delta_X^{\otimes s})[1]= 0.$
 
 We have $\delta_X=\sum_a(\delta_X)_a\otimes a$ and $\delta_X[1]=\sum_a\phi(X)_{v(a)}(\delta_X)_a\phi(X)^{-1}_{u(a)}\otimes a[1]$.
 Given the filtration of right $\hat{S}$-submodules 
 $0=X_0\subseteq X_1\subseteq \cdots \subseteq X_r=X$
 such that $(\delta_X)_a(X_i)\subseteq X_{i-1}$, for all $i$ and $a$, we have the 
 filtration of right $\hat{S}$-submodules 
 $0=X_0[1]\subseteq X_1[1]\subseteq \cdots \subseteq X_r[1]=X[1]$
 such that 
 $$\begin{matrix}
 ((\delta_X)[1])_{a[1]}(X_i[1])
 &=&
 \phi(X)_{v(a)}(\delta_X)_a\phi(X)^{-1}_{u(a)}(X_i[1])\hfill\\
  &\subseteq&
  \phi(X)_{v(a)}(\delta_X)_a(X_i)\hfill\\
  &\subseteq& \phi(X)_{v(a)}(X_{i-1})  \subseteq X_{i-1}[1].\hfill\\   
   \end{matrix}$$
 Then,  we have $(X[1],\delta_X[1])\in {\cal Z}(\hat{Z})$. If $f\in {\cal Z}(\hat{Z})((X,\delta_X),(Y,\delta_Y))$, then
 $$\begin{matrix}
   0& =&\sum_{i_0,i_1\geq 0}\hat{b}^{ad}_{i_0+i_1+1}(\delta_Y^{\otimes i_1}\otimes f\otimes \delta_X^{\otimes i_0})[1]\hfill\\
&=&
\sum_{i_0,i_1\geq 0}\hat{b}^{ad}_{i_0+i_1+1}((\delta_Y[1])^{\otimes i_1}\otimes f[1]\otimes (\delta_X[1])^{\otimes i_0}),\hfill\\
   \end{matrix}$$
and $f[1]\in {\cal Z}(\hat{Z})((X,\delta_X)[1],(Y,\delta_Y)[1])$. 
Whenever $f\in {\cal Z}(\hat{Z})((X,\delta_X),(Y,\delta_Y))$
and  $g\in {\cal Z}(\hat{Z})((Y,\delta_Y),(W,\delta_W))$, from (\ref{L: translacion automorfismo de la b-cat ad(hat(B))}), we have 
$T(g)\star T(f)=g[1]\star f[1]=\hat{b}_2^{tw}(g[1]\otimes f[1])=\hat{b}_2^{tw}(g\otimes f)[1]=(g\star f)[1]=T(g\star f).$
So $T$ preserves the composition of ${\cal Z}(\hat{Z})$. Moreover, we have $T(\hueca{I}_X)=T (\sum_{u\in \hat{\cal P}}id_{Xe_u}\otimes \frak{e}_u)
  =
  \sum_{u\in \hat{\cal P}}id_{X[1]e_{\nu(u)}}\otimes \frak{e}_{\nu(u)}
  =\hueca{I}_{X[1]}.$
So, $T$ preserves identities.  It is easy to see that the association $f\mapsto f[-1]$ determines an inverse for the functor $T$. 

In order to show that $T$ induces an autofunctor $T:{\cal H}(\hat{Z})\rightmap{}{\cal H}(\hat{Z})$, it is enough to show that $T({\cal I})={\cal I}$, where ${\cal I}=\hat{b}_1^{tw}[\tw(\hat{Z})(-,?)_{-2}]$. Indeed, if $f\in {\cal I}((X,\delta_X),(Y,\delta_Y))$, there is some   $h\in\tw(\hat{Z})((X,\delta_X),(Y,\delta_Y))_{-2}$  such that $f=\hat{b}_1^{tw}(h)$. Then, we have  
$
   f[1]=\sum_{i_0,i_1\geq 0}\hat{b}^{ad}_{i_0+i_1+1}(\delta_Y^{\otimes i_1}\otimes h\otimes \delta_X^{\otimes i_0})[1]
   =
   \sum_{i_0,i_1\geq 0}\hat{b}^{ad}_{i_0+i_1+1}((\delta_Y[1])^{\otimes i_1}\otimes h[1]\otimes (\delta_X[1])^{\otimes i_0}),
  $
so $f[1]=\hat{b}_1^{tw}(h[1])$. Similarly, we have  $f[-1]=\hat{b}_1^{tw}(h[-1])$. So $T({\cal I})={\cal I}$ and $T$ induces an autofunctor $T:{\cal H}(\hat{Z})\rightmap{}{\cal H}(\hat{Z})$ as we wanted to show. 
 \end{proof}

\begin{definition}\label{D: conflations in Z(hat(B))}
Given $(X,\delta_X),(Y,\delta_Y)\in {\cal Z}(\hat{Z})$, 
 we denote by 
 $$\Ext_{{\cal Z}(\hat{Z})}((Y,\delta_Y),(X,\delta_X))$$ 
 the collection of equivalence classes $[\xi]$ of special conflations  in ${\cal Z}(\hat{Z})$, for the equivalence relation ``$\rightmap{\simeq}$'', see (\ref{L: equiv de canonical conflation is equiv relation}), of the form 
 $$\xi:\hbox{ \ }(X,\delta_X)\rightmap{f}(E,\delta_E)\rightmap{g}(Y,\delta_Y).$$
 \end{definition}

 \begin{lemma}\label{L: cada gama Y-->X determina una confla X-->E-->Y }
  Every morphism $h\in \Hom_{{\cal Z}(\hat{Z})}((Y,\delta_Y),(X,\delta_X)[1])$ determines a canonical conflation $\xi_h$ in ${\cal Z}(\hat{Z})$ of the form
  $$\xi_h:\hbox{ \ } (X,\delta_X)\rightmap{f}(E,\delta_E)\rightmap{g}(Y,\delta_Y),$$
  where $E=X\oplus Y$, $\delta_E=\begin{pmatrix}
                                 \delta_X&-\tau_X\circ h\\ 0&\delta_Y 
                                 \end{pmatrix}$, $f=(\hueca{I}_X,0)^t$, and $g=(0,\hueca{I}_Y)$.
 \end{lemma}

  \begin{proof} We have $\gamma:=-\tau_X\circ h\in \ad(\hat{Z})(Y,X)$.  Since  $\vert h\vert=-1$, $\vert \tau_X\vert=0$, and $\vert \hat{b}_2^{tw}\vert=1$, we have  $\vert \gamma\vert=\vert \hat{b}_2^{tw}(\tau_X\circ h)\vert=0$. From (\ref{L: diferenciales de la suma directa en tw(hat(B))}), we will have that $(E,\delta_E)$ is an object of ${\cal Z}(\hat{Z})$ once we have verified that $\hat{b}_1^{tw}(\gamma)=0$. Indeed, from (\ref{L: reemplazo tensors delta por tensors de delta[1]}) and (\ref{L: bn(sigmaXncirc fn...)=sigmaXn(bn fn...)}), we have 
  $$\begin{matrix}
     \hat{b}_1^{tw}(\gamma)&=&\sum_{i_0,i_1\geq 0}\hat{b}^{ad}_{i_0+i_1+1}(\delta_X^{\otimes i_1}\otimes  \gamma\otimes \delta_Y^{\otimes i_0})\hfill\\
     &=&
     -\tau_X\circ \sum_{i_0,i_1\geq 0}\hat{b}^{ad}_{i_0+i_1+1}((\delta_X[1])^{\otimes i_1}\otimes (\sigma_X\circ \gamma)\otimes \delta_Y^{\otimes i_0})\hfill\\
     &=&
     -\tau_X\circ \sum_{i_0,i_1\geq 0}\hat{b}^{ad}_{i_0+i_1+1}((\delta_X[1])^{\otimes i_1}\otimes h\otimes \delta_Y^{\otimes i_0})=0,\hfill\\
    \end{matrix}$$
because $\sigma_X\circ \gamma=\sigma_X\circ (-\tau_X\circ h)=h$, according to (\ref{L: circ de tres sigmas y taus con morfismos}). 

From (\ref{L: diferenciales de la suma directa en tw(hat(B))}), we get that $(E,\delta_E)$ is an object of ${\cal Z}(\hat{Z})$.
Then, from (\ref{L: previo a forma canonica de conflaciones}), we know that $f$ and $g$ are morphisms in ${\cal Z}(\hat{Z})$. So, the composable pair $\xi_h$ is a special conflation of ${\cal Z}(\hat{Z})$.
  \end{proof}

 \begin{proposition}\label{P: biy :Hom(X,Y[1])<-->Ext(Y,X)}
  The map $$\Psi:\Hom_{{\cal Z}(\hat{Z})}((Y,\delta_Y),(X,\delta_X)[1])\rightmap{}\Ext_{{\cal Z}(\hat{Z})}((Y,\delta_Y),(X,\delta_X))$$
  such that $h\mapsto [\xi_h]$ is a surjection and induces a bijection 
  $$\Psi:\Hom_{{\cal H}(\hat{Z})}((Y,\delta_Y),(X,\delta_X)[1])\rightmap{}\Ext_{{\cal Z}(\hat{Z})}((Y,\delta_Y),(X,\delta_X)).$$
 \end{proposition}
 
\begin{proof} In order to show that $\Psi$ is surjective, we consider a canonical conflation 
 $$\xi:\hbox{ \ }(X,\delta_X)\rightmap{f}(E,\delta_E)\rightmap{g}(Y,\delta_Y).$$
 and let us find $h\in \Hom_{{\cal Z}(\hat{Z})}((Y,\delta_Y),(X,\delta_X)[1])$ with $[\xi_h]=[\xi]$. Recall that, as shown in  (\ref{L: may assume special conflation direct sum middle term}), any special conflation is equivalent to a canonical one. Since 
  $\xi$ is a canonical conflation, we have $E=X\oplus Y$,  $f=(\hueca{I}_X,0)^t$, $g=(0,\hueca{I}_Y)$, and   $$\delta_E=
 \begin{pmatrix}
 \delta_X&\gamma\\ 0&\delta_Y\\
 \end{pmatrix},$$ for some homogeneous morphism 
 $\gamma:Y\rightmap{}X$ in $\ad(\hat{Z})$ of degree $0$.
 
 From (\ref{L: diferenciales de la suma directa en tw(hat(B))}), we have 
 $0=\hat{b}^{tw}_1(\gamma)=\sum_{i_0,i_1\geq 0}\hat{b}^{ad}_{i_0+i_1+1}(\delta_X^{\otimes i_1}\otimes \gamma\otimes\delta_Y^{\otimes i_0})$. 
 Then, from (\ref{L: circ de tres sigmas y taus con morfismos}), (\ref{L: bn(sigmaXncirc fn...)=sigmaXn(bn fn...)}), and (\ref{L: reemplazo tensors delta por tensors de delta[1]}), we have 
 $$\begin{matrix}
   \hat{b}^{tw}_1(\sigma_X\circ\gamma)&=& \sum_{i_0,i_1\geq 0}\hat{b}^{ad}_{i_0+i_1+1}(\delta_X[1]^{\otimes i_1}\otimes\sigma_X\circ \gamma\otimes\delta_Y^{\otimes i_0})\hfill\\
   &=&
   -\sigma_X\circ(\tau_X\circ \sum_{i_0,i_1\geq 0}\hat{b}^{ad}_{i_0+i_1+1}(\delta_X[1]^{\otimes i_1}\otimes\sigma_X\circ \gamma\otimes\delta_Y^{\otimes i_0})\hfill\\
   &=&
   \sigma_X\circ \sum_{i_0,i_1\geq 0}\hat{b}^{ad}_{i_0+i_1+1}(\delta_X^{\otimes i_1}\otimes\gamma\otimes\delta_Y^{\otimes i_0})=0.\hfill\\
   \end{matrix}$$
 Then, we have that $h:=\sigma_X\circ \gamma\in \Hom_{{\cal Z}(\hat{Z})}((Y,\delta_Y),(X,\delta_X)[1])$ satisfies that $-\tau_X\circ h=\gamma$, thus $\Psi(h)=[\xi_h]=[\xi]$, and $\Psi$ is surjective.
 \medskip

 It remains to show that whenever $h,h_1\in \Hom_{{\cal Z}(\hat{Z})}((Y,\delta_Y),(X,\delta_X)[1])$ we have $\Psi(h)=\Psi(h_1)$ iff $h-h_1\in {\cal I}$. 
  
 Assume first that $[\xi_h]=[\xi_{h_1}]$.
 Set $\gamma:=-\tau_X\circ h$ and $\gamma_1:=-\tau_X\circ h_1$. Then, we have a commutative diagram in ${\cal Z}(\hat{Z})$
  $$\begin{matrix}
  \xi_h:\hbox{ \ } (X,\delta_X)&\rightmap{f}&(E,\delta_E)&\rightmap{g}&(Y,\delta_Y)\hfill\\
 \hbox{ \  \ \ }   \shortlmapdown{\hueca{I}_X}&&\shortrmapdown{t}&&\shortrmapdown{\hueca{I}_Y}\\
  \xi_{h_1}:\hbox{ \ }   (X,\delta_X)&\rightmap{f_1}&(E_1,\delta_{E_1})&\rightmap{g_1}&(Y,\delta_Y),\hfill\\
   \end{matrix}$$
   where $E,E_1$, $f,g$, and $f_1,g_1$ have the form described above and $t$ is an isomorphism of ${\cal Z}(\hat{Z})$. Thus $E=X\oplus Y=E_1$ as right $\hat{S}$-modules and
   $$\delta_{E}=
 \begin{pmatrix}
 \delta_X&\gamma\\ 0&\delta_Y\\
 \end{pmatrix} \hbox{ \ and \  }\delta_{E_1}=
 \begin{pmatrix}
 \delta_X&\gamma_1\\   0&\delta_Y\\
 \end{pmatrix}.$$
 
 By (\ref{R: star-compos=circ-compos para special morfisms}), the commutativity of the diagram in ${\cal Z}(\hat{Z})$ implies its commutativity in $\ad(\hat{Z})$ because $f,f_1,g,g_1$ are all special morphisms. Then the morphism  $t$ has the matrix form 
 $$t=\begin{pmatrix}
    \hueca{I}_X&s\\
    0&\hueca{I}_Y\\
  \end{pmatrix},$$ 
  where $s:Y\rightmap{}X$ is a homogeneous morphism in $\ad(\hat{Z})$ with degree $-1$.  From (\ref{L: isos matriciales entre terminos de enmedio}), we have  the equality $\gamma_1-\gamma=\sum_{i_0,i_1\geq 0}\hat{b}^{ad}_{i_0+i_1+1}(\delta_X^{\otimes i_1}\otimes s\otimes\delta_Y^{\otimes i_0})$. 
Notice that $\sigma_X\circ \gamma=-\sigma_X\circ (\tau_X\circ h)=h$ and, similarly, $\sigma_X\circ \gamma_1=h_1$. Then, from (\ref{L: como entra la sigma en medio de la hat(b)adn(... sigmaX...)}), we get
$$\begin{matrix}h-h_1&=&
   \sigma_X\circ \gamma-\sigma_X\circ \gamma_1\hfill\\
   &=&
   -\sum_{i_0,i_1\geq 0}\sigma_X\circ \hat{b}^{ad}_{i_0+i_1+1}(\delta_X^{\otimes i_1}\otimes s\otimes\delta_Y^{\otimes i_0})\hfill\\
   
   &=&
   \sum_{i_0,i_1\geq 0}\hat{b}^{ad}_{i_0+i_1+1} (\delta_X[1]^{\otimes i_1}\otimes \sigma_X\circ s\otimes\delta_Y^{\otimes i_0})\hfill\\
   &=&
   \hat{b}_1^{tw}(\sigma_X\circ s).\hfill\\
  \end{matrix}$$
Here, the composition $\sigma_X\circ s:Y\rightmap{}X[1]$ is a homogeneous 
 morphism in $\ad(\hat{Z})$ with degree $\vert\sigma_X\circ s\vert=-2$. Hence $h-h_1\in {\cal I}$ as we wanted to show.  

Conversely, if  $h-h_1\in {\cal I}$, we have 
$h-h_1=\hat{b}_1^{tw}(r)$, for some morphism $r:Y\rightmap{}X[1]$ in $\ad(\hat{Z})$ with degree $-2$. Then, the morphism $s:=-\tau_X\circ r: Y\rightmap{}X$ is homogeneous with degree $-1$, and we have 
$h-h_1=\hat{b}_1^{tw}(\sigma_X\circ s).$ 
Then, using again (\ref{L: isos matriciales entre terminos de enmedio}), we can reverse the above argument to show that  the morphism 
$t:=\begin{pmatrix}
    \hueca{I}_X&s\\
    0&\hueca{I}_Y\\
  \end{pmatrix}:(E,\delta_E)\rightmap{}(E_1,\delta_{E_1})$ 
belongs to ${\cal Z}(\hat{Z})$. It clearly makes the preceding diagram commutative with respect to $\circ$.
From (\ref{L: inverso de matriz u triangular en Z(hat(B))}), we know that $t$ is an isomorphism in ${\cal Z}(\hat{Z})$ and, hence, we get
$[\xi_h]=[\xi_{h_1}]$.
 \end{proof}
 
\begin{remark}\label{R: h homol trivial --> xih splits}
In (\ref{L: cada gama Y-->X determina una confla X-->E-->Y }), to the morphism $h=0$ corresponds the trivial conflation $\xi_0:(X,\delta_X)\rightmap{f}(X,\delta_X)\oplus(Y,\delta_Y)\rightmap{g}(Y,\delta_Y)$. From (\ref{L: xi equiv a trivial implica xi splits}), we get that, whenever $h$ is homologically trivial,
 the corresponding special conflation $\xi_h$ splits.
 \end{remark}

 For the construction of the endofunctor $J$, we will use the following maps.

\begin{lemma}\label{L: twb1(sigma)=0 y twb1(tau)=0}
Given any object $(X,\delta_X)$ in ${\cal Z}(\hat{Z})$, we have the following strict homogeneous morphisms of $\tw(\hat{Z})$ 
$$\sigma_X:(X,\delta_X)\rightmap{}(X[1],\delta_{X[1]}) \hbox{ \ \ and \ \ }
\tau_X:(X[1],\delta_{X[1]})\rightmap{}(X,\delta_X)$$
of degrees $-2$ and $0$, respectively, which satisfy $\hat{b}_1^{tw}(\sigma_X)=0$ and $\hat{b}_1^{tw}(\tau_X)=0$. 
\end{lemma}

\begin{proof} The strict morphism $\sigma_X:X\rightmap{}X[1]$ in $\ad(\hat{Z})$ has degree $-2$ and, from (\ref{R: cuando un estricto esta en Z(hat(B))}), it satisfies $\hat{b}_1^{tw}(\sigma_X)=\delta_{X[1]}\circ \sigma_X+\sigma_X\circ \delta_X$. Applying  (\ref{L: circ de tres sigmas y taus con morfismos}), we obtain  
$\delta_{X[1]}\circ \sigma_X=
-((\sigma_X\circ\delta_X)\circ \tau_X)\circ\sigma_X=
-\sigma_X\circ \delta_X.$
Hence, we get $\hat{b}_1^{tw}(\sigma_X)=0$. 

Similarly, the strict morphism $\tau_X:X[1]\rightmap{}X$ in $\ad(\hat{Z})$ has degree $0$ and, from (\ref{R: cuando un estricto esta en Z(hat(B))}), it satisfies $\hat{b}_1^{tw}(\tau_X)=\delta_X\circ \tau_X+\tau_X\circ \delta_{X[1]}$. Applying (\ref{L: circ de tres sigmas y taus con morfismos}), we get 
$\tau_X\circ \delta_{X[1]}=
\tau_X\circ (\sigma_X\circ(\delta_X\circ \tau_X))=
-\delta_X\circ \tau_X.$
Hence, we get $\hat{b}_1^{tw}(\tau_X)=0$.
\end{proof}

\begin{lemma}\label{L: introduccion de J(X,deltaX)}
Given any object $(X,\delta_X)$ in ${\cal Z}(\hat{Z})$, we consider the right $\hat{S}$-module $J(X):=X\oplus X[1]$ and the morphism 
$$\delta_{J(X)}:=\begin{pmatrix}
                 \delta_X&-\tau_X\\ 0&\delta_{X[1]}\\ 
                 \end{pmatrix}:J(X)\rightmap{}J(X) \hbox{ in } \ad(\hat{Z}).$$
 Then, the pair $J(X,\delta_X):=(J(X),\delta_{J(X)})$ is an object in ${\cal Z}(\hat{Z})$. It is homologically trivial in ${\cal H}(\hat{Z})$.
\end{lemma}

\begin{proof} From (\ref{L: twb1(sigma)=0 y twb1(tau)=0}), the homogeneous morphism $\tau_X:(X[1],\delta_{X[1]})\rightmap{}(X,\delta_X)$ in $\ad(\hat{Z})$ has degree $0$ and satisfies $\hat{b}_1^{tw}(\tau_X)=0$. 
From (\ref{L: diferenciales de la suma directa en tw(hat(B))}), we get that $(J(X),\delta_{J(X)})$ is an object in ${\cal Z}(\hat{Z})$. 

In order to show that $J(X,\delta_X)$ is homologically trivial, we have to exhibit   some $s\in \tw(\hat{Z})(J(X,\delta_X),J(X,\delta_X))_{-2}$ such that $\hat{b}_1^{tw}(s)=id_{J(X,\delta_X)}$. Consider the strict homogeneous morphism of degree $-2$ in $\tw(\hat{Z})$
$$s=\begin{pmatrix}
     0&0\\ -\sigma_X&0\\
    \end{pmatrix}:(J(X),\delta_{J(X)})\rightmap{}(J(X),\delta_{J(X)}).$$
From (\ref{R: cuando un estricto esta en Z(hat(B))}), we have that 
$\hat{b}_1^{tw}(s)=\delta_{J(X)}\circ s+s\circ \delta_{J(X)}$. From (\ref{L: tauX circ sigmaX}) and (\ref{L: twb1(sigma)=0 y twb1(tau)=0}), we get 
$$\begin{matrix}
   \hat{b}_1^{tw}(s)&=&\begin{pmatrix}
                        \tau_X\circ \sigma_X&0\\ -\delta_{X[1]}\circ \sigma_X&0\\
                       \end{pmatrix}
                       +
                       \begin{pmatrix}
                        0&0\\ -\sigma_X\circ \delta_X&\sigma_X\circ \tau_X\\
                       \end{pmatrix}\hfill\\
                       &\,&\\
                       &=&
                       \begin{pmatrix}
                        \hueca{I}_X&0\\ -\hat{b}_1^{tw}(\sigma_X)&\hueca{I}_{X[1]}\\
                       \end{pmatrix}
                       =
                       \begin{pmatrix}
                        \hueca{I}_X&0\\ 0&\hueca{I}_{X[1]}\\
                       \end{pmatrix}=id_{J(X,\delta_X)}.\hfill\\
 \end{matrix}$$
Hence $J(X,\delta_X)$ is homologically trivial. 
\end{proof}

\begin{lemma}\label{L: morf homol triv se fact por conflaciones esp}
 If $(X,\delta_X)\rightmap{f}(E,\delta_E)\rightmap{g}(Y,\delta_Y)$ is a special conflation, we have:
 \begin{enumerate}
  \item For any homologically trivial morphism $h:(X,\delta_X)\rightmap{}(X_1,\delta_{X_1})$ there is a morphism $h':(E,\delta_E)\rightmap{}(X_1,\delta_{X_1})$ in ${\cal Z}(\hat{Z})$ such that $h'\star f=h$. 
  \item For any homologically trivial morphism $h:(Y_1,\delta_{Y_1})\rightmap{}(Y,\delta_Y)$ there is a morphism $h':(Y_1,\delta_{Y_1})\rightmap{}(E,\delta_E)$ in ${\cal Z}(\hat{Z})$ such that $g\star h'=h$. 
 \end{enumerate}
\end{lemma}

\begin{proof} We only prove (1), because the proof of (2) is similar. 
 From (\ref{L: may assume special conflation direct sum middle term}), we see that it is enough to prove this for canonical conflations. From (\ref{L: pushouts para conflaciones en Z(hat(B))}), we have a commutative diagram 
 $$\begin{matrix}
    \xi&:&(X,\delta_X)&\rightmap{f}&(E,\delta_E)&\rightmap{g}&(Y,\delta_Y)\\
    &&\shortlmapdown{h}&&\shortrmapdown{t}&&\shortrmapdown{\hueca{I}_Y}\\
    \xi'&:&(X_1,\delta_{X_1})&\rightmap{f_1}&(E_1,\delta_{E_1})&\rightmap{g_1}&(Y,\delta_Y)\\
   \end{matrix}$$
where the second row is a canonical conflation. Notice that the morphism $h':=\sigma_{X_1}\star (h\star \gamma):(Y,\delta_Y)\rightmap{}(X_1,\delta_{X_1})[1]$ is homologically trivial. As a consequence of  (\ref{R: h homol trivial --> xih splits}), we obtain that  the canonical conflation $\xi_{h'}=\xi'$   splits. Here we have the $\hat{S}$-modules $E=X\oplus Y$ and $E_1=X_1\oplus Y$, with  differentials of the form
$\delta_E=\begin{pmatrix}   \delta_X&\gamma\\ 0&\delta_Y                                                                                                  \end{pmatrix}$ and
$\delta_{E_1}=\begin{pmatrix}   \delta_{X_1}&\gamma_1\\ 0&\delta_Y                                                                                                  \end{pmatrix}$. Moreover,
$t=\begin{pmatrix}     h&0\\ 0&\hueca{I}_Y\\                                                                                                                                                                                                                                                                                                                                                                                                                                                                                                   \end{pmatrix}$.

 Consider a left inverse
 $f'_1=(\hueca{I}_{X_1},s):
 (E_1,\delta_{E_1})=(X_1\oplus Y, \delta_{E_1})
 \rightmap{}(X_1,\delta_{X_1})$ for $f_1$ in ${\cal Z}(\hat{Z})$, see (\ref{L: sobre conflaciones triviales}).  Define $h':=f'_1\star t:(E,\delta_E)\rightmap{}(X_1,\delta_{X_1})$.
 From (\ref{R: desarmado f en estricto + otro simplifica btw's})(2), we have $f'_1\star t=f'_1\circ t+R$, where
 $$\begin{matrix}
   R=\sum_{\scriptsize\begin{matrix}i_0,i_1,i_2\geq 0\\
                          i_0+i_1+i_2\geq 1\end{matrix}}  \hat{b}^{ad}_{i_0+i_1+i_2}(\delta_{X_1}^{\otimes i_2}\otimes (0,s)\otimes \delta_{E_1}^{\otimes i_1}\otimes \begin{pmatrix}   h&0\\ 0&0\\                                                                                                          \end{pmatrix}\otimes \delta_{E}^{\otimes i_0})=0,
   \end{matrix}$$
because $(0,s)\otimes \begin{pmatrix}   \delta_{X_1}&\gamma_1\\ 0&\delta_Y                                                                                                  \end{pmatrix}^{\otimes i_1}\otimes \begin{pmatrix}   h&0\\ 0&0\\                                                                                                          \end{pmatrix}=0$. Thus, $h'=f'_1\star t=f'_1\circ t=(h,s)$. Finally, we get $h'\star f=h'\circ f=h$, as we wanted to show. 
\end{proof}

\begin{corollary}\label{C: J(X) proy e iny cra confl especiales} 
Any object in ${\cal Z}(\hat{Z})$ of the form $J(U,\delta_U)$ is projective and injective relative to special conflations.
\end{corollary}

\begin{proof} If the sequence $(X,\delta_X)\rightmap{f}(E,\delta_E)\rightmap{g}(Y,\delta_Y)$ is a special conflation and $h:(X,\delta_X)\rightmap{}J(U,\delta_U)$ is a morphism in ${\cal Z}(\hat{Z})$, we get from (\ref{L: introduccion de J(X,deltaX)}) that $h$ is homologically trivial. From (\ref{L: morf homol triv se fact por conflaciones esp}), we obtain that $h$ factors through $f$, and $J(U,\delta_U)$ is injective relative to special conflations. The statement on projectivity is proved similarly. 
\end{proof}

\begin{lemma}\label{L: el funtor J:Z(hat(Z))-->Z(hat(Z))}
There is a functor $J:{\cal Z}(\hat{Z})\rightmap{}{\cal Z}(\hat{Z})$ which maps each morphism $f:(X,\delta_X)\rightmap{}(X',\delta_{X'})$ of ${\cal Z}(\hat{Z})$ on the morphism 
$$J(f):=\begin{pmatrix}
         f&0\\ 0&f[1]\\
        \end{pmatrix}:(J(X),\delta_{J(X)})\rightmap{}(J(X'),\delta_{J(X')}).$$
\end{lemma}

\begin{proof} In order to show that $\hat{b}_1^{tw}(J(f))=0$ using (\ref{R: desarmado f en estricto + otro simplifica btw's})(1), we consider the morphisms  
$\delta_{J(X)}^0=\begin{pmatrix}
0&-\tau_X\\ 0&0\\                                                                         \end{pmatrix}$, $\delta_{J(X)}^1=\begin{pmatrix}
\delta_X&0\\ 0&\delta_{X[1]}\\                                                                                       \end{pmatrix}$, $\delta_{J(X')}^0=\begin{pmatrix}
0&-\tau_{X'}\\ 0&0\\                                                                         \end{pmatrix}$, and $\delta_{J(X')}^1=\begin{pmatrix}
\delta_{X'}&0\\ 0&\delta_{X'[1]}\\                                                                                       \end{pmatrix}$. Then, we have 
$$\hat{b}_1^{tw}(J(f))=\delta_{J(X')}\circ J(f)+J(f)\circ \delta_{J(X)}+R(J(f)),$$
where 
$R(J(f))=\hat{b}_1^{ad}(J(f))+\sum_{\scriptsize\begin{matrix}i_0,i_1\geq 0\\ 
                          i_0+i_1\geq 2\end{matrix}} \hat{b}^{ad}_{i_0+i_1+1}((\delta_{J(X')}^1)^{\otimes i_1}\otimes J(f)\otimes (\delta_{J(X)}^1)^{\otimes i_0})$, which is a $2\times 2$ diagonal matrix with diagonal terms
$$\begin{matrix}D_1=\hat{b}_1^{ad}(f)+\sum_{\scriptsize\begin{matrix}i_0,i_1\geq 0\\ 
                          i_0+i_1\geq 2\end{matrix}} \hat{b}^{ad}_{i_0+i_1+1}(\delta_{X'}^{\otimes i_1}\otimes f\otimes \delta_{X}^{\otimes i_0})\hbox{ \ and }\\ 
D_2=
    \hat{b}_1^{ad}(f[1])+\sum_{\scriptsize\begin{matrix}i_0,i_1\geq 0\\ 
                          i_0+i_1\geq 2\end{matrix}} \hat{b}^{ad}_{i_0+i_1+1}((\delta_{X'[1]})^{\otimes i_1}\otimes f[1]\otimes (\delta_{X[1]})^{\otimes i_0}).
                          \end{matrix}$$
Moreover, we have 
$$\delta_{J(X')}\circ J(f)+J(f)\circ \delta_{J(X)}=
\begin{pmatrix}
f\circ \delta_X+\delta_{X'}\circ f&-f\circ \tau_X-\tau_{X'}\circ f[1]\\
 0& f[1]\circ \delta_{X[1]}+\delta_{X'[1]}\circ f[1]
                                                   \end{pmatrix}.$$
From (\ref{L: circ de tres sigmas y taus con morfismos})(3), we have $f[1]=\sigma_{X'}\circ (f\circ \tau_X)$; so, by (\ref{L: circ de tres sigmas y taus con morfismos})(1), we have $$\tau_{X'}\circ f[1]=\tau_{X'}\circ (\sigma_{X'}\circ (f\circ \tau_X))=-(f\circ \tau_X).$$ 
Therefore, 
$\hat{b}_1^{tw}(J(f))=\begin{pmatrix}\hat{b}_1^{tw}(f)&0\\ 0&\hat{b}_1^{tw}(f[1])\end{pmatrix}=0$. 

In order to verify that $J$ preserves the composition $\star$, we consider another morphism $g:(X',\delta_{X'})\rightmap{}(X'',\delta_{X''})$ and set $\delta_{J(X'')}^0=\begin{pmatrix}
0&-\tau_{X''}\\ 0&0\\                                                                         \end{pmatrix}$, and $\delta_{J(X'')}^1=\begin{pmatrix}
\delta_{X''}&0\\ 0&\delta_{X''[1]}\\                                                                                       \end{pmatrix}$. We  use   (\ref{R: desarmado f en estricto + otro simplifica btw's})(2) to show that $J(g\star f)=J(g)\star J(f)$. 
We have $J(g)\star J(f)=J(g)\circ J(f)+R(J(g),J(f))$, where $R(J(g),J(f))$ denotes 
$$\sum_{\scriptsize\begin{matrix}i_0,i_1,i_2\geq 0\\ 
                          i_0+i_1+i_2\geq 1\end{matrix}} \hat{b}^{ad}_{i_0+i_1+i_2+2}((\delta_{J(X'')}^1)^{\otimes i_2}\otimes J(g)\otimes(\delta_{J(X')}^1)^{\otimes i_1}\otimes J(f)\otimes (\delta_{J(X)}^1)^{\otimes i_0}),$$
                          which is a $2\times 2$ diagonal matrix with diagonal terms
$$\begin{matrix}
 \hueca{D}_1=\sum_{\scriptsize\begin{matrix}i_0,i_1,i_2\geq 0\\ 
                          i_0+i_1+i_2\geq 1\end{matrix}} \hat{b}^{ad}_{i_0+i_1+i_2+2}(\delta_{X''}^{\otimes i_2}\otimes g\otimes\delta_{X'}^{\otimes i_1}\otimes f\otimes \delta_X^{\otimes i_0}) \hbox{ and }\\ 
\hueca{D}_2=  \sum_{\scriptsize\begin{matrix}i_0,i_1,i_2\geq 0\\ 
                          i_0+i_1+i_2\geq 1\end{matrix}} \hat{b}^{ad}_{i_0+i_1+i_2+2}(\delta_{X''[1]}^{\otimes i_2}\otimes g[1]\otimes\delta_{X'[1]}^{\otimes i_1}\otimes f[1]\otimes \delta_{X[1]}^{\otimes i_0}).                        
  \end{matrix}$$
We also have that $g\star f=g\circ f+R(g,f)$, where $R(g,f)=\hueca{D}_1$; and $(g\star f)[1]=g[1]\star f[1]=g[1]\circ f[1]+R(g[1],f[1])$, where $R(g[1],f[1])=\hueca{D}_2$. Therefore, 
$$J(g)\star J(f)=\begin{pmatrix}
                  g&0\\ 0&g[1]\\
                 \end{pmatrix}\circ 
                 \begin{pmatrix}
                  f&0\\ 0&f[1]\\
                 \end{pmatrix}+
                 \begin{pmatrix}
                  \hueca{D}_1&0\\ 0&\hueca{D}_2
                 \end{pmatrix}=
                 \begin{pmatrix}
                  g\star f&0\\ 0&(g\star f)[1]
                 \end{pmatrix}=J(g\star f).$$
\end{proof}

\begin{remark}\label{L: naturalidad de alfa y beta hacia J y desde J}
Given $(X,\delta_X)\in {\cal Z}(\hat{Z})$, since $\hat{b}_1^{tw}(-\tau_X)=0$, from (\ref{L: previo a forma canonica de conflaciones}),  we have in ${\cal Z}(\hat{Z})$  the canonical conflation
$$\xi(X,\delta_X):(X,\delta_X)\rightmap{\alpha_X}(J(X),\delta_{J(X)})\rightmap{\beta_X}(X,\delta_X)[1],$$
 with $\alpha_X=(\hueca{I}_X,0)^t$  and $\beta_X=(0,\hueca{I}_{X[1]})$.
 
Given any morphism $f:(X,\delta_X)\rightmap{}(X',\delta_{X'})$ in 
${\cal Z}(\hat{Z})$, we have the commutative diagram in ${\cal Z}(\hat{Z})$ 
$$\begin{matrix}
  \xi(X,\delta_X):&(X,\delta_X)&\rightmap{\alpha_X}&(J(X),\delta_{J(X)})&\rightmap{\beta_X}&(X,\delta_X)[1]\\
  &\shortlmapdown{f}&&\shortlmapdown{J(f)}&&\shortrmapdown{f[1]}\\
  \xi(X',\delta_{X'}):&(X',\delta_{X'})&\rightmap{\alpha_{X'}}&(J(X'),\delta_{J(X')})&\rightmap{\beta_{X'}}&(X',\delta_{X'})[1].\\
  \end{matrix}$$
\end{remark}

\begin{remark}\label{R: H(hat(Z)) y stable cat de Z(hat(Z))} Having in mind (\ref{L: xi equiv a trivial implica xi splits}),
notice that if $(X,\delta_X)$ is a projective (resp. injective) object of ${\cal Z}(\hat{Z})$ relative to special conflations, then the special  conflation given by (\ref{L: naturalidad de alfa y beta hacia J y desde J}) splits and, therefore,  $(X,\delta_X)$ is a direct summand of $J(X,\delta_X)$. From this and
(\ref{C: J(X) proy e iny cra confl especiales}), we get that  projective (resp. injective) objects in ${\cal Z}(\hat{Z})$ relative to special conflations are the direct summands of the  objects of the form $J(U,\delta_U)$. 

Again from (\ref{L: naturalidad de alfa y beta hacia J y desde J}), we see that ${\cal Z}(\hat{Z})$ has enough injectives and enough projectives, meaning that for any object $(X,\delta_X)$ we have a special deflation from a relative projective onto $(X,\delta_X)$, and a special inflation from $(X,\delta_X)$ into a relative injective. 

Moreover, from (\ref{L: morf homol triv se fact por conflaciones esp}), we get that a morphism $f:(X,\delta_X)\rightmap{}(Y,\delta_Y)$ in ${\cal Z}(\hat{Z})$ is homologically trivial iff it factors through a relative projective object in ${\cal Z}(\hat{Z})$. Thus, the cohomology category 
${\cal H}(\hat{Z})$ is \emph{the stable category of} ${\cal Z}(\hat{Z})$, that is the category obtained from ${\cal Z}(\hat{Z})$ by factoring out morphisms which factor through relative projectives. 
\end{remark}

\begin{lemma}\label{L: morfism de conflaciones con morfism medio en J(X)}
 Suppose that 
 $\xi:(X,\delta_X)\rightmap{f}(E,\delta_E)\rightmap{g}(Y,\delta_Y)$ is a canonical conflation in ${\cal Z}(\hat{Z})$ with 
 $E=X\oplus Y$ and 
 $\delta_E=\begin{pmatrix}
            \delta_X&\gamma\\ 0&\delta_Y\\
            \end{pmatrix}$, for some homogeneous  
$\gamma:Y\rightmap{}X$ in $\ad(\hat{Z})$ of degree $0$. Then, we have the following commutative diagram in ${\cal Z}(\hat{Z})$ 
    $$\begin{matrix}
   \xi:&(X,\delta_X)&\rightmap{f}&(E,\delta_E)&\rightmap{g}&\hbox{ \ }(Y,\delta_Y)\hfill\\
   &\shortrmapdown{\hueca{I}_X}&&\shortrmapdown{h_\xi}&&\shortrmapdown{h_\gamma}\\
   \xi(X,\delta_X):&(X,\delta_X)&\rightmap{\alpha_X}&(J(X),\delta_{J(X)})&\rightmap{\beta_X}&(X,\delta_X)[1],\hfill\\
   \end{matrix}$$ 
   where $h_\gamma=-\sigma_X\circ \gamma$ and $h_\xi=\begin{pmatrix}
             \hueca{I}_X&0\\ 0&h_\gamma\\
            \end{pmatrix}$.
\end{lemma}

\begin{proof}  From (\ref{L: twb1(sigma)=0 y twb1(tau)=0}) and (\ref{L: diferenciales de la suma directa en tw(hat(B))}), we have $\hat{b}_1^{tw}(\sigma_X)=0$ and  $\hat{b}_1^{tw}(\gamma)=0$. Therefore,   $\hat{b}_1^{tw}(h_\gamma)=-\hat{b}_1^{tw}(\sigma_X\circ \gamma)=0$, and the morphism $h_\gamma$  
belongs to  ${\cal Z}(\hat{Z})$.

From (\ref{L: pullbacks para conflaciones en Z(hat(B))}) applied to $\xi(X,\delta_X)$ and the morphism $h_\gamma$ of ${\cal Z}(\hat{Z})$
there is a commutative diagram in ${\cal Z}(\hat{Z})$ of the form  
  $$\begin{matrix}
 \hbox{ \ \ }& (X,\delta_X)&\rightmap{ (\hueca{I}_X,0)^t }&(E',\delta_{E'})&\rightmap{ (0,\hueca{I}_Y) }&\hbox{ \ \ }(Y,\delta_Y)\hfill\\
 &\hbox{ \  \ \ }  \shortlmapdown{\hueca{I}_X}&&\shortrmapdown{h_\xi}&&\shortrmapdown{h_\gamma}\\
 \xi(X,\delta_X):&\hbox{ \ }  (X,\delta_X)&\rightmap{ \ \alpha_X \ }&(J(X),\delta_{J(X)})&\rightmap{ \ \beta_X \ }&(X[1],\delta_{X[1]}),\hfill\\
   \end{matrix}$$
 where $E'=X\oplus Y$,  
 $\delta_{E'}=\begin{pmatrix}
                             \delta_X&\tau_X\star h_\gamma\\0&\delta_Y\\
                            \end{pmatrix}=
                            \begin{pmatrix}
                            \delta_X&\gamma\\ 0&\delta_Y\\ 
                            \end{pmatrix}=\delta_E$, 
and  $h_\xi=\begin{pmatrix}
    \hueca{I}_X&0\\ 0&h_\gamma\\
    \end{pmatrix}:(E,\delta_E)=(E',\delta_{E'})\rightmap{}(J(X),\delta_{J(X)})$. 
\end{proof}

\begin{proposition}\label{P: cada confl determina equiv de confl con J(X)}
 Suppose that 
 $\xi:(X,\delta_X)\rightmap{f}(E,\delta_E)\rightmap{g}(Y,\delta_Y)$ is a canonical conflation in  ${\cal Z}(\hat{Z})$, so $E=X\oplus Y$ and $\delta_E=
 \begin{pmatrix}
 \delta_X&\gamma\\ 0&\delta_Y\\                                                                                                                                                                                                                                                                                                            \end{pmatrix}$, for some homogeneous morphism 
 $\gamma:Y\rightmap{}X$ in $\ad(\hat{Z})$ of degree $0$.
  It determines the following pair of composable morphisms in ${\cal Z}(\hat{Z})$: 
 $$\eta:\hbox{ \ }(E,\delta_E)\rightmap{\alpha}J(X,\delta_X)\oplus(Y,\delta_Y)\rightmap{\beta}(X,\delta_X)[1],$$
 where $\alpha=(h_\xi,g)^t$ and $\beta=(\beta_X,-h_\gamma)$, with the notation of (\ref{L: morfism de conflaciones con morfism medio en J(X)}). 
   The composable pair of morphisms $\eta$ is a conflation, as in (\ref{D: def de conflation}).
\end{proposition}

\begin{proof}
 We will construct another composable pair $\overline{\eta}$ and isomorphisms $s$ and $s'$ in ${\cal Z}(\hat{Z})$ such that $\eta\rightmap{s\simeq}\overline{\eta}$ and $\overline{\eta}\rightmap{s'\simeq}\eta_1$, where $\eta_1$ is the canonical conflation 
$$\eta_1:\hbox{ \ }(E,\delta_E)\rightmap{\alpha_1}(E_1,\delta_{E_1})\rightmap{\beta_1}(X,\delta_X)[1],$$
where, 
$E_1=E\oplus X[1]$, $\delta_{E_1}=\begin{pmatrix}
                                \delta_E&\gamma_1\\ 0&\delta_{X[1]}\\
                               \end{pmatrix}$, $\gamma_1=(-\tau_X,0)^t$, $\alpha_1=(\hueca{I}_E,0)^t$, and  $\beta_1=(0,\hueca{I}_{X[1]})$.

 If we define $(E_\xi,\delta_{E_\xi}):=J(X,\delta_X)\oplus(Y,\delta_Y)$, we get
 $$E_\xi=J(X)\oplus Y\hbox{ \ \ and \ \ } \delta_{E_\xi}=
 \begin{pmatrix}
 \delta_{J(X)}&0\\ 0&\delta_Y\\
 \end{pmatrix}.$$
                      
 \noindent We have  the special isomorphism 
 $s=\begin{pmatrix}
     s_{1,1}&s_{1,2}\\ s_{2,1}&s_{2,2}
    \end{pmatrix}:J(X)\oplus Y\rightmap{}E\oplus X[1],$
    with $s_{1,1}=\begin{pmatrix}
                   \hueca{I}_X&0\\ 0&0\\
                  \end{pmatrix}$,
         $s_{2,1}=(
                   0,\hueca{I}_{X[1]})$,
         $s_{1,2}=\begin{pmatrix}
                   0\\ \hueca{I}_Y\\
                  \end{pmatrix}$, and $s_{2,2}=0$. Its special inverse is 
                  $r=s^{-1}=\begin{pmatrix}
                             r_{1,1}&r_{1,2}\\ r_{2,1}&r_{2,2}\\
                            \end{pmatrix}:E\oplus X[1]\rightmap{}J(X)\oplus Y$,
      with components $r_{1,1}=\begin{pmatrix}
                     \hueca{I}_X&0\\ 0&0\\
                    \end{pmatrix}$,
           $r_{2,1}=(0,\hueca{I}_Y)$,
           $r_{1,2}=\begin{pmatrix}
                     0\\ \hueca{I}_{X[1]}
                    \end{pmatrix}$, and $r_{2,2}=0$. 

     From (\ref{L: transfiriendo diferenciales})(2), we have the following object in ${\cal Z}(\hat{Z})$
     $$(\overline{E}_\xi,\delta_{\overline{E}_\xi}):=(E\oplus X[1],-s\circ\delta_{E_\xi}\circ s^{-1}).$$
  Therefore, we have   
  $\delta_{\overline{E}_\xi}=\begin{pmatrix}
                             \delta'_E&\overline{\gamma}\\
                             0&\delta_{X[1]}\\
                            \end{pmatrix}$,
with 
$\delta'_E=
                     \begin{pmatrix}
                     \delta_X&0\\ 0&\delta_Y\\
                     \end{pmatrix}$ 
 and 
                         $\overline{\gamma}=
                            \begin{pmatrix}
                             -\tau_X\\ 0\\
                            \end{pmatrix}:X[1]\rightmap{}X\oplus Y=E.$  Moreover, the morphism $s:(E_\xi,\delta_{E_\xi})\rightmap{}(\overline{E}_\xi,\delta_{\overline{E}_\xi})$ is an isomorphism in ${\cal Z}(\hat{Z})$. 
Consider the following  diagram in $\tw(\hat{Z})$:
$$\begin{matrix}
 \eta:\hbox{ \ }(E,\delta_E)&\rightmap{\alpha}&(E_\xi,\delta_{E_\xi})&\rightmap{\beta}&(X,\delta_X)[1]\hfill\\
  \hbox{ \ \ \ \ \ }  \shortlmapdown{\hueca{I}_E}&&\shortrmapdown{s}&&\shortrmapdown{\hueca{I}_{X[1]}}\\
 \overline{\eta}:\hbox{ \ }   (E,\delta_E)&\rightmap{\overline{\alpha}}&(\overline{E}_\xi,\delta_{\overline{E}_\xi})&\rightmap{\overline{\beta}}&(X,\delta_{X})[1]\hfill\\
    \end{matrix}$$
where $\overline{\alpha}=(\hueca{I}_E,-\rho)^t$ and $\overline{\beta}=(\rho,\hueca{I}_{X[1]})$, with $\rho:=(0,-h_\gamma):X\oplus Y\rightmap{}X[1]$.
The preceding diagram commutes in $\tw(\hat{Z})$ because:
$$s\star \alpha=s\circ \alpha=(s_{1,1}\circ h_\xi+s_{1,2}\circ 
g,s_{2,1}\circ h_\xi)^t=(\hueca{I}_E,-\rho)^t=\overline{\alpha}$$
 and 
$\beta\star s^{-1}=\beta\circ s^{-1}=(\beta_X\circ r_{1,1}+\sigma_X\circ \gamma\circ r_{2,1},\beta_X\circ r_{1,2})=(\rho,\hueca{I}_{X[1]})=\overline{\beta}$.

From (\ref{L: morfism de conflaciones con morfism medio en J(X)}), we know that the morphisms $h_\xi:(E,\delta_E)\rightmap{}J(X,\delta_X)$, 
$\beta_X:J(X,\delta_X)\rightmap{}(X,\delta_X)[1]$, 
and $h_\gamma:(Y,\delta_Y)\rightmap{}(X,\delta_X)[1]$ belong to ${\cal Z}(\hat{Z})$. By assumption, so does the morphism  
$g:(E,\delta_E)\rightmap{}(Y,\delta_Y)$. 
So, the components of the morphisms $\alpha$ and $\beta$ belong to ${\cal Z}(\hat{Z})$. By  (\ref{L: Z(hat(B)) is an additive precategory}), this implies that the morphisms $\alpha$ and $\beta$ lie in ${\cal Z}(\hat{Z})$. Then, the sequence $\eta$ lies in ${\cal Z}(\hat{Z})$, and so does the sequence $\overline{\eta}$. Therefore,  we have 
$\eta\rightmap{s\simeq}\overline{\eta}$. 

\medskip
We already know that $(E,\delta'_E)$ is an object of ${\cal Z}(\hat{Z})$ and that $\overline{\gamma}:X[1]\rightmap{}E$ is a strict morphism. Moreover, we have 
$$\begin{matrix}
 \delta_E\circ \overline{\gamma}+\overline{\gamma}\circ \delta_{X[1]}
&=&  
 -(\delta_X\circ \tau_X+\tau_X\circ \delta_{X[1]},0)^t=-(\hat{b}_1^{tw}(\tau_X),0)^t=0.\hfill\\                 
  \end{matrix}$$

  Then, applying  (\ref{L: variando deltas en inicios de conflaciones}), to the sequence $\overline{\eta}$, we have the commutative diagram 
  in ${\cal Z}(\hat{Z})$:                                          
$$\begin{matrix}
 \overline{\eta}:\hbox{ \ }   (E,\delta_E)&\rightmap{\overline{\alpha}}&(\overline{E}_\xi,\delta_{\overline{E}_\xi})&\rightmap{\overline{\beta}}&(X,\delta_{X})[1]\hfill\\
 \hbox{ \ \ \ \ \ }  \shortlmapdown{\hueca{I}_E}&&\shortrmapdown{s'}&&\shortrmapdown{\hueca{I}_{X[1]}}\\
    \eta_1:\hbox{ \ }(E,\delta_E)&\rightmap{\alpha_1}&(E_1,\delta_{E_1})&\rightmap{\beta_1}&(X,\delta_X)[1],\hfill\\
    \end{matrix}$$          
    where $s'$ is an isomorphism in ${\cal Z}(\hat{Z})$. Thus, we get $\overline{\eta}\rightmap{ \ s'\simeq \ }\eta_1$, as claimed.  
\end{proof}

The following lemma is similar to  (\ref{L: morfism de conflaciones con morfism medio en J(X)}). 

\begin{lemma}\label{L: morfism de conflaciones con morfism medio en J'(X)} 
Suppose that 
 $\xi:(X,\delta_X)\rightmap{f}(E,\delta_E)\rightmap{g}(Y,\delta_Y)$ is a canonical conflation in ${\cal Z}(\hat{Z})$ with 
 $E=X\oplus Y$ and 
 $\delta_E=\begin{pmatrix}
            \delta_X&\gamma\\ 0&\delta_Y\\
            \end{pmatrix}$, for some homogeneous  
$\gamma:Y\rightmap{}X$ in $\ad(\hat{Z})$ of degree $0$. Then, we have the following commutative diagram in ${\cal Z}(\hat{Z})$ 
    $$\begin{matrix}
 \xi(Y[-1],\delta_{Y[-1]}):&  (Y,\delta_Y)[-1]&\rightmap{\alpha_{Y[-1]}}&(J(Y[-1]),\delta_{J(Y[-1])})&\rightmap{\beta_{Y[-1]}}&(Y,\delta_Y)\hfill\\
   &\shortlmapdown{h^\gamma}&&\shortrmapdown{h^\xi}&&\shortrmapdown{\hueca{I}_Y}\\
   \xi:&(X,\delta_X)&\rightmap{f}&(E,\delta_E)&\rightmap{g}&(Y,\delta_Y),\hfill\\
   \end{matrix}$$ 
   where $h^\gamma=-(\sigma_X\circ \gamma)[-1]$ and $h^\xi=\begin{pmatrix}
             h^\gamma&0\\ 0&\hueca{I}_Y\\
            \end{pmatrix}$.
\end{lemma}

\begin{proof} We already know that $h_\gamma=-\sigma_X\circ \gamma:(Y,\delta_Y)\rightmap{}(X,\delta_X)[1]$ is a morphism in ${\cal Z}(\hat{Z})$, and so is $h^\gamma=h_\gamma[-1]$. 
 
From (\ref{L: pushouts para conflaciones en Z(hat(B))}) applied to  $\xi(Y[-1],\delta_{Y[-1]})$ and the morphism $h^\gamma$ of ${\cal Z}(\hat{Z})$, we have a commutative diagram in ${\cal Z}(\hat{Z})$ of the form
$$\begin{matrix}
 \xi(Y[-1],\delta_{Y[-1]}):&  (Y,\delta_Y)[-1]&\rightmap{\alpha_{Y[-1]}}&(J(Y[-1]),\delta_{J(Y[-1])})&\rightmap{\beta_{Y[-1]}}&(Y,\delta_Y)\hfill\\
   &\shortlmapdown{h^\gamma}&&\shortrmapdown{h^\xi}&&\shortrmapdown{\hueca{I}_Y}\\
   &(X,\delta_X)&\rightmap{(\hueca{I}_X,0)^t}&(E',\delta_{E'})&\rightmap{(0,\hueca{I}_Y)}&(Y,\delta_Y),\hfill\\
   \end{matrix}$$ 
where $E'=X\oplus Y$ and 
 $\delta_{E'}=\begin{pmatrix}
                             \delta_X&-h^\gamma\star\tau_{Y[-1]} \\0&\delta_Y\\
                            \end{pmatrix}=
                            \begin{pmatrix}
                            \delta_X&\gamma\\ 0&\delta_Y\\ 
                            \end{pmatrix}=\delta_E.$
 Indeed, from (\ref{R: sigmaX[-1] y tauX[-1]}),  we know that $\sigma_X[-1]=-\sigma_{X[-1]}$ and then, from (\ref{L: circ de tres sigmas y taus con morfismos}), we get 
 $
 -h^\gamma\star \tau_{Y[-1]}
 =
 (\sigma_X\circ \gamma)[-1]\circ \tau_{Y[-1]}
 =
 (\sigma_X[-1]\circ \gamma[-1])\circ \tau_{Y[-1]}
 =
 -(\sigma_{X[-1]}\circ \gamma[-1])\circ \tau_{Y[-1]}
 =
 \gamma[-1][1]=\gamma$. 
Moreover, we have  $h^\xi=\begin{pmatrix}
    h^\gamma&0\\ 0&\hueca{I}_Y\\
    \end{pmatrix}:(J(Y[-1]),\delta_{J(Y[-1])})\rightmap{}(E',\delta_{E'})=(E,\delta_E)$. 
\end{proof}

The following proposition is similar  (\ref{P: cada confl determina equiv de confl con J(X)}).  

\begin{proposition}\label{P: cada confl determina equiv de confl con J'(Y)} Suppose that 
 $\xi:(X,\delta_X)\rightmap{f}(E,\delta_E)\rightmap{g}(Y,\delta_Y)$ is a canonical conflation in  ${\cal Z}(\hat{Z})$, so $E=X\oplus Y$ and $\delta_E=
 \begin{pmatrix}
 \delta_X&\gamma\\ 0&\delta_Y\\                                                                                                                                                                                                                                                                                                            \end{pmatrix}$, for some homogeneous morphism 
 $\gamma:Y\rightmap{}X$ in $\ad(\hat{Z})$ of degree $0$. It determines the following pair of composable morphisms in ${\cal Z}(\hat{Z})$: 
 $$\eta:\hbox{ \ }(Y,\delta_Y)[-1]\rightmap{\alpha}J(Y[-1],\delta_{Y[-1]})\oplus(X,\delta_X)\rightmap{\beta}(E,\delta_E),$$
 where $\alpha=(\alpha_{Y[-1]},-h^\gamma)^t$ and $\beta=(h^\xi,f)$, with the notation of (\ref{L: morfism de conflaciones con morfism medio en J'(X)}). The composable pair of morphisms $\eta$ is a conflation, as in (\ref{D: def de conflation}).
\end{proposition}

\begin{proof}   Similar to the proof of (\ref{P: cada confl determina equiv de confl con J(X)}), now using (\ref{L: morfism de conflaciones con morfism medio en J'(X)}) and (\ref{L: variando deltas en finales de conflaciones}). 
\end{proof}

\begin{proposition}\label{P: morfismos encajados en conflaciones con J}
Any morphism $f:(X,\delta_X)\rightmap{}(Y,\delta_Y)$ in ${\cal Z}(\hat{Z})$ determines a conflation of the form 
$$(X,\delta_X)\rightmap{ \ \overline{\alpha}=(\alpha_X,f)^t \ }J(X,\delta_X)\oplus (Y,\delta_Y)\rightmap{ \ \ \overline{\beta} \ \ }(W,\delta'_W).$$
\end{proposition}

\begin{proof} We define $W:=Y\oplus X[1]$ and
$\delta'_W:=\begin{pmatrix}
            \delta_Y&f\circ \tau_X\\
            0&\delta_{X[1]}\\ 
            \end{pmatrix}$. 
  Since $\hat{b}_1^{tw}(f)=0$ and $\hat{b}_1^{tw}(\tau_X)=0$, we have that $\hat{b}_1^{tw}(f\circ \tau_X)=0$. Therefore, we have that $(W,\delta'_W)$ is an object of ${\cal Z}(\hat{Z})$.  We also have the object $(W,\delta_W)=(Y,\delta_Y)\oplus (X,\delta_X)[1]$ in ${\cal Z}(\hat{Z})$, with $\delta_W=\begin{pmatrix}
           \delta_Y&0\\ 0&\delta_{X[1]}\\                                                                                                                                                                                                              \end{pmatrix}$.
Consider the strict homogeneous morphism $\gamma:=(0,-\tau_X):W=Y\oplus X[1]\rightmap{}X$ with degree $0$. From (\ref{L: circ de tres sigmas y taus con morfismos}), we have 
  $$\delta_X\circ \tau_X=-\tau_X\circ (\sigma_X\circ (\delta_X\circ \tau_X))=-\tau_X\circ \delta_X[1]=-\tau_X\circ \delta_{X[1]}.$$
  Hence, we get 
 $\hat{b}_1^{tw}(\gamma)=\delta_X\circ \gamma+\gamma\circ \delta_W
  =
  (0,-\delta_X\circ \tau_X)+(0,-\tau_X\circ \delta_{X[1]})=0$. 
  So, we have the object $(E,\delta_E)$ in ${\cal Z}(\hat{Z})$ defined by $E=X\oplus W$ and $\delta_E=\begin{pmatrix}
\delta_X&\gamma\\ 0&\delta_W 
\end{pmatrix}$. Set $(\overline{E},\delta_{\overline{E}}):=J(X,\delta_X)\oplus (Y,\delta_Y)$.  We are interested in the following  diagram, which clearly commutes in $\tw(\hat{Z})$:
$$\begin{matrix} \overline{\eta}:\hbox{ \ }(X,\delta_X)&\rightmap{ \ \overline{\alpha}=(\alpha_X,f)^t \ }&(\overline{E},\delta_{\overline{E}})&\rightmap{ \ \overline{\beta}=(h,g) \ \ }&(W,\delta'_W)\hfill\\
 \hbox{ \ \ \ \ \ }  
 \shortlmapdown{\hueca{I}_X}&&\shortrmapdown{t}&&\shortrmapdown{\hueca{I}_W}\\
  \eta:\hbox{ \ }(X,\delta_X)&\rightmap{ \ \  \alpha=(\hueca{I}_X,-\rho)^t \ }&(E,\delta_E)&\rightmap{ \ \ \beta=(\rho,\hueca{I}_W) \ \  }&(W,\delta'_W),\hfill\\  
    \end{matrix}$$ 
where $E=X\oplus W$, $h=\begin{pmatrix}
  -f&0\\ 0&\hueca{I}_{X[1]}\\
 \end{pmatrix}:J(X)=X\oplus X[1]
  \rightmap{}Y\oplus X[1]=W$, 
$\rho=(-f,0)^t:X\rightmap{}Y\oplus X[1]=W$, 
$\alpha:X\rightmap{}X\oplus W=E$, $g=(\hueca{I}_Y,0)^t:Y\rightmap{}Y\oplus X[1]=W$,  
 and 
  $$t=\begin{pmatrix}
  \hueca{I}_X&0&0\\
  0&0&\hueca{I}_Y\\
  0&\hueca{I}_{X[1]}&0\\
  \end{pmatrix}:\overline{E}=X\oplus X[1]\oplus Y\rightmap{}X\oplus Y\oplus X[1]=E.$$  
  It is not hard to show that $t:(\overline{E},\delta_{\overline{E}})\rightmap{}(E,\delta_E)$ is in fact an isomorphism in ${\cal Z}(\hat{Z})$. So in order to show that $\overline{\eta}$ is a conflation, it will be enough to show that $\eta$ is so. In order to apply (\ref{L: variando deltas en finales de conflaciones}) to $\eta$, we use that $(W,\delta_W)$ is an object of ${\cal Z}(\hat{Z})$, 
  that $\gamma:W\rightmap{}X$ is strict and satisfies 
  $\delta_X\circ \gamma+\gamma\circ \delta'_W=0$. It only remains to show that $\alpha$ and $\beta$ are  morphisms in ${\cal Z}(\hat{Z})$. We will use (\ref{L: caractizacion de pares en Z(hat(B)) con derecha variada}), so we need to show that 
   $\rho:(X,\delta_X)\rightmap{}(W,\delta_W)$ and $\rho':(X,\delta_X)\rightmap{}(W,\delta'_W)$  are morphisms in ${\cal Z}(\hat{Z})$ 
   such that $\gamma\circ\rho=0$ and $\rho\circ\gamma=\delta'_W-\delta_W$. 

    Clearly, we have  $\gamma\circ \rho=0$ and $\rho\circ \gamma=\begin{pmatrix}
           0&f\circ \tau_X\\
           0&0
           \end{pmatrix}=\delta'_W-\delta_W$. Moreover, we have 
   $$\begin{matrix}\hat{b}_1^{tw}(\rho')&=&\sum_{i_0,i_1\geq 0}\hat{b}_{i_0+i_1+1}^{ad}((\delta'_W)^{\otimes i_1}\otimes \rho\otimes \delta_X^{\otimes i_0})\hfill\\
   &=&
      \begin{pmatrix}-\sum_{i_0,i_1\geq 0}\hat{b}_{i_0+i_1+1}^{ad}(\delta_Y^{\otimes i_1}\otimes f\otimes \delta_X^{\otimes i_0}\\
       0
      \end{pmatrix}=0,
\hfill\\
     \end{matrix}$$
     and, similarly, we have $\hat{b}_1^{tw}(\rho)=0$. So, we get  $\overline{\eta}\rightmap{t\simeq}\eta$ where $\eta$ is a conflation, and $\overline{\eta}$ is a conflation too.  
\end{proof}

\section{The triangulated category ${\cal H}(\hat{Z})$}

 Now, with the notation of the last section, we will prove that the category ${\cal H}(\hat{Z})$ is triangulated. We first recall some basic definitions.

\begin{definition}\label{D: cat pretriangulada}  Assume that ${\cal H}$ is an additive $k$-category together with an autofunctor $T:{\cal H}\rightmap{}{\cal H}$. A \emph{sextuple $t=(X,Y,U,u,v,w)$ in ${\cal H}$} is a sequence of composable morphisms in ${\cal H}$ of the form $$t:{\ }X\rightmap{u}Y\rightmap{v}U\rightmap{w}TX.$$
A morphism of sextuples $(\theta_1,\theta_2,\theta_3):(X,Y,U,u,v,w)\rightmap{}(X',Y',U',u',v',w')$ is a triple of morphisms such that the following diagram commutes:
$$\begin{matrix}
X&\rightmap{u}&Y&\rightmap{v}&U&\rightmap{w}&TX\\
\shortlmapdown{\theta_1}&&\shortlmapdown{\theta_2}&&
\shortrmapdown{\theta_3}&&\shortrmapdown{T(\theta_1)}\\
X'&\rightmap{u'}&Y'&\rightmap{v'}&U'&\rightmap{w'}&TX'.\\
\end{matrix}$$ 
The category ${\cal H}$ is called a \emph{pretriangulated category}   if it is equipped with a class ${\cal T}$ of sextuples  
$X\rightmap{u}Y\rightmap{v}U\rightmap{w}TX$, 
called \emph{the triangles of} ${\cal H}$, such that:  
\begin{enumerate}
 \item[TR1:]  
 \begin{enumerate}
 \item For any isomorphism between two sextuples  
  such that one of them is a triangle,  so is the other one.
  \item The sextuple $X\rightmap{id_X}X\rightmap{0}0\rightmap{0}TX$ is a triangle, for any  $X\in {\cal H}$.  
  \item For each morphism $u:X\rightmap{}Y$ in ${\cal H}$, there is a triangle of the form $$X\rightmap{u}Y\rightmap{v}U\rightmap{w}TX.$$
 \end{enumerate}
 \item[TR2:] The sextuple $X\rightmap{u}Y\rightmap{v}U\rightmap{w}TX$ is a triangle if and only if 
  the sextuple $Y\rightmap{v}U\rightmap{w}TX\rightmap{-T(u)}TY$ is a triangle. 
 \item[TR3:] Each commutative diagram 
 $$\begin{matrix}
 t:&X&\rightmap{u}&Y&\rightmap{v}&U&\rightmap{w}&TX\\
 &\shortlmapdown{\theta_1}&&\shortlmapdown{\theta_2}&&&&\\
 t':&X'&\rightmap{u'}&Y'&\rightmap{v'}&U'&\rightmap{w'}&TX',\\
   \end{matrix}$$
   such that the rows $t$ and $t'$ are triangles, can be completed to a morphism of triangles $(\theta_1,\theta_2,\theta_3):t\rightmap{}t'$.
\end{enumerate}
A pretriangulated category ${\cal H}$ is called \emph{triangulated} iff its triangles  furthermore satisfy the following axiom:
\begin{enumerate}
 \item[TR4:]     \emph{Octahedral Axiom:} Given triangles

    $$\begin{matrix}X&\rightmap{u}&Y&\rightmap{i}&U'&\rightmap{\hat{i}}&TX\\
       Y&\rightmap{v}&U&\rightmap{j}&X'&\rightmap{\hat{j}}&TY\\
       X& \rightmap{vu}&U&\rightmap{w}&Y'&\rightmap{\hat{w}}&TX\\
      \end{matrix}$$
 there is a triangle 
     $U'\rightmap{ \ f \ }Y'\rightmap{ \ g \ }X'\rightmap{ \ \ T(i)\hat{j} \ \ }TU'$ such that the following diagram commutes
 $$\begin{matrix}
    T^{-1}Y'&\rightmap{ \ T^{-1}\hat{w} \ }&X&\rightmap{1_X}&X&&&&\\
    \shortlmapdown{ \ T^{-1}(g) \ }&&\shortlmapdown{u}&&\shortrmapdown{vu}&&&&\\
    T^{-1}X'&\rightmap{T^{-1}(\hat{j})}&Y&\rightmap{v}&U&\rightmap{j}&X'&\rightmap{\hat{j}}&TY\\
    &&\shortlmapdown{i}&&\shortrmapdown{w}&&\shortlmapdown{1_{X'}}&&\shortrmapdown{T(i)}\\
    &&U'&\rightmap{f}&Y'&\rightmap{g}&X'&\rightmap{ \ T(i)\hat{j} \ }&TU'\\
    &&\shortlmapdown{\hat{i}}&&\shortrmapdown{\hat{w}}&&&&\\
    &&TX&\rightmap{1_{TX}}&TX.&&&&\\
   \end{matrix}$$
\end{enumerate}
\end{definition}

\begin{remark}\label{R: notacion y H(hat(B)) aditiva}
We keep the notation used in the last section and we denote by $\pi:{\cal Z}(\hat{Z})\rightmap{}{\cal H}(\hat{Z})$ the canonical projection. From (\ref{L: Z(hat(B)) is an additive precategory}), we already know  that ${\cal H}(\hat{Z})$ is an additive $k$-category. Moreover, we have endowed ${\cal H}(\hat{Z})$ with a $k$-linear autofunctor $T:{\cal H}(\hat{Z})\rightmap{}{\cal H}(\hat{Z})$ in (\ref{L: traslacion en Z(hat(B)) y en H(hat(B))}).
\end{remark}

\begin{definition}\label{D: canonical triangles and triangles}
A \emph{canonical triangle in ${\cal H}(\hat{Z})$} is a sextuple of the form
$$\tau_\xi:\hbox{ \ }(X,\delta_X)\rightmap{\pi(f)}(E,\delta_E)\rightmap{\pi(g)}(Y,\delta_Y)\rightmap{w}(X,\delta_X)[1]$$
such that 
$\xi:(X,\delta_X)\rightmap{f}(E,\delta_E)\rightmap{g}(Y,\delta_Y)$
is a canonical conflation in ${\cal Z}(\hat{Z})$ and $\Psi(w)=[\xi]$, see (\ref{P: biy :Hom(X,Y[1])<-->Ext(Y,X)}). Notice that this is equivalent to ask that the sextuple is of the form 
$$\tau_\xi:\hbox{ \ }(X,\delta_X)\rightmap{\pi(f)}(E,\delta_E)\rightmap{\pi(g)}(Y,\delta_Y)\rightmap{\pi(\sigma_X\circ \gamma)}(X,\delta_X)[1],$$
for some canonical conflation 
$(X,\delta_X)\rightmap{f}(E,\delta_E)\rightmap{g}(Y,\delta_Y)$, 
where $E=X\oplus Y$\\ and $\delta_E=\begin{pmatrix}
           \delta_X&\gamma\\ 0&\delta_Y\\                                                          
           \end{pmatrix}$. 
   By definition, a \emph{triangle} in ${\cal H}(\hat{Z})$ is any sextuple isomorphic to some canonical triangle.
\end{definition}

\begin{lemma}\label{R: conflations give rise to triangles}
Every conflation $\xi: \underline{X}\rightmap{f}\underline{E}\rightmap{g}\underline{Y}$ in ${\cal Z}(\hat{Z})$ gives rise to a triangle in ${\cal H}(\hat{Z})$: If $\xi$ transforms into the canonical conflation 
$\xi_n: \underline{X}\rightmap{f_n}\underline{E}_n\rightmap{g_n}\underline{Y},$
we have an isomorphism of triangles: 
$$\begin{matrix}\tau_\xi:&\underline{X}&\rightmap{ \ \pi(f) \  }&\underline{E}&\rightmap{ \ \pi(g) \ }&\underline{Y}&\rightmap{w}&T\underline{X}\\
&\shortlmapdown{\hueca{I}_X}&&\shortlmapdown{\cong}&&\shortrmapdown{\hueca{I}_Y}&&\shortrmapdown{\hueca{I}_{TX}}\\
\tau_{\xi_n}:&\underline{X}&\rightmap{ \ \pi(f_n) \  }&\underline{E}_n&\rightmap{ \ \pi(g_n) \ }&\underline{Y}&\rightmap{w}&T\underline{X},\\
\end{matrix}$$
where $\Psi(w)=[\xi_n]$.
\end{lemma}

\begin{proof} If we have a sequence of conflations $\xi_0,\ldots,\xi_n$ in ${\cal Z}(\hat{Z})$ and relations 
$$  \xi_0\rightmap{\simeq}\xi_1\leftmap{\simeq}\xi_2\rightmap{\simeq}\cdots \leftmap{\simeq}\xi_{n-1}\rightmap{\simeq}\xi_{n-1}\leftmap{\simeq}\xi_n,$$
where $\xi=\xi_0$ and $\xi_n$ is a canonical conflation, then we have a commutative diagram in ${\cal H}(\hat{Z})$

$$\begin{matrix}\underline{X}&\rightmap{ \ \pi(f) \  }&\underline{E}&\rightmap{ \ \pi(g) \ }&\underline{Y}&\rightmap{w}&T\underline{X}\\

\shortlmapdown{\hueca{I}_X}&&\shortlmapdown{\theta_1}&&
\shortrmapdown{\hueca{I}_Y}&&\shortrmapdown{\hueca{I}_{TX}}\\

\underline{X}&\rightmap{ \ \pi(f_1) \  }&\underline{E}_1&
\rightmap{ \ \pi(g_1) \ }&\underline{Y}&\rightmap{w}&T(\underline{X})\\

\shortlmapup{\hueca{I}_X}&&\shortlmapup{\theta_2}&&
\shortrmapup{\hueca{I}_Y}&&\shortrmapup{\hueca{I}_{TX}}\\

\vdots&&\vdots&&\vdots&&\vdots\\

\shortlmapup{\hueca{I}_X}&&\shortlmapup{\theta_n}&&
\shortrmapup{\hueca{I}_Y}&&\shortrmapup{\hueca{I}_{TX}}\\
\underline{X}&\rightmap{ \ \pi(f_n) \  }&\underline{E}_n&\rightmap{ \ \pi(g_n) \ }&\underline{Y}&\rightmap{w}&T\underline{X},\\
\end{matrix}$$
with $\theta_1,\ldots,\theta_n$ isomorphisms and $\Psi(w)=[\xi_n]$. Since the last row of the diagram is a canonical triangle,  
the first row of the diagram is a triangle.  
\end{proof}

\begin{lemma}\label{L: formulas para morfismos de triangulos}
 \begin{enumerate}
  \item Given $\gamma:(Y,\delta_Y)\rightmap{}(X,\delta_X)$ and $t_1:(X,\delta_X)\rightmap{}(X',\delta_{X'})$ homogeneous morphisms in  $\tw(\hat{Z})$ with degrees $0$ and $-1$, respectively, we have 
  $\sigma_{X'}\circ (t_1\star \gamma)=t_1[1]\star(\sigma_X\circ \gamma).$
 
 \item Given $t_3:(Y,\delta_Y)\rightmap{}(Y',\delta_{Y'})$ and $\gamma':(Y',\delta_{Y'})\rightmap{}(X',\delta_{X'})$ homogeneous morphisms in $\tw(\hat{Z})$ with degrees $-1$ and zero, respectively, we have 
  $\sigma_{X'}\circ(\gamma'\star t_3)=-(\sigma_{X'}\circ \gamma')\star t_3.$
 
 \item Consider a commutative diagram in ${\cal Z}(\hat{Z})$ with canonical conflations as rows
  $$\begin{matrix}
     (X,\delta_X)&\rightmap{}&(E,\delta_E)&\rightmap{}&(Y,\delta_Y)\\
     \shortlmapdown{t_1}&&\shortrmapdown{t_2}&&\shortrmapdown{t_3}\\
     (X',\delta_{X'})&\rightmap{}&(E',\delta_{E'})&\rightmap{}&(Y',\delta_{Y'})\\
    \end{matrix}$$
    where $E=X\oplus Y$, $E'=X'\oplus Y'$, $\delta_E=\begin{pmatrix}
                                                      \delta_X&\gamma\\ 0&\delta_Y
                                                     \end{pmatrix}$, and 
$\delta_{E'}=\begin{pmatrix}
             \delta_{X'}&\gamma'\\ 0&\delta_{Y'}\\
             \end{pmatrix}$. Then, we have  in ${\cal H}(\hat{Z})$ the equality
             $\pi(t_1)[1]\pi(\sigma_X\circ \gamma)=\pi(\sigma_{X'}\circ \gamma')\pi(t_3)$. 
 \end{enumerate}
\end{lemma}

\begin{proof}(1): From (\ref{L: como entra la sigma en medio de la hat(b)adn(... sigmaX...)}), we have 
$$\begin{matrix}
 \sigma_{X'}\circ (t_1\star \gamma)
 &=&
  \sum_{i_0,i_1,i_2\geq 0}\sigma_{X'}\circ 
  \hat{b}^{ad}_{i_0+i_1+i_2+2}(\delta_{X'}^{\otimes i_2}\otimes t_1\otimes \delta_X^{\otimes i_1}\otimes\gamma\otimes\delta_Y^{\otimes i_0})   \hfill\\
  &=&
   \sum_{i_0,i_1,i_2\geq 0} 
  \hat{b}^{ad}_{i_0+i_1+i_2+2}(\delta_{X'[1]}^{\otimes i_2}\otimes t_1[1]\otimes \delta_{X[1]}^{\otimes i_1}\otimes(\sigma_X\circ \gamma)\otimes\delta_Y^{\otimes i_0})   \hfill\\
  &=& t_1[1]\star (\sigma_X\circ \gamma).\hfill\\
  \end{matrix}$$
\noindent(2): As before, from (\ref{L: como entra la sigma en medio de la hat(b)adn(... sigmaX...)}), we get 
$$\begin{matrix}
 \sigma_{X'}\circ (\gamma'\star t_3)
 &=&
  \sum_{i_0,i_1,i_2\geq 0}\sigma_{X'}\circ 
  \hat{b}^{ad}_{i_0+i_1+i_2+2}(\delta_{X'}^{\otimes i_2}\otimes\gamma'\otimes \delta_{Y'}^{\otimes i_1}\otimes t_3\otimes\delta_Y^{\otimes i_0})   \hfill\\
  &=&
   \sum_{i_0,i_1,i_2\geq 0} 
  \hat{b}^{ad}_{i_0+i_1+i_2+2}(\delta_{X'[1]}^{\otimes i_2}\otimes(\sigma_{X'}\circ \gamma')\otimes \delta_{Y'}^{\otimes i_1}\otimes t_3\otimes\delta_Y^{\otimes i_0})   \hfill\\
  &=& -(\sigma_{X'}\circ \gamma')\star t_3.\hfill\\
  \end{matrix}$$

\noindent(3): We have 
$t_2=\begin{pmatrix}
      v_{1,1}&v_{1,2}\\ v_{2,1}&v_{2,2}\\
     \end{pmatrix}:X\oplus Y\rightmap{}X'\oplus Y'$. 
     From the commutativity of the diagram, we have
  $t_2\star(\hueca{I}_X,0)^t=t_2\circ(\hueca{I}_X,0)^t=(\hueca{I}_{X'},0)^t\circ t_1$  and, therefore,  $v_{1,1}=t_1$  and  $v_{2,1}=0$;
and
$(0,\hueca{I}_{Y'})\star t_2=(0,\hueca{I}_{Y'})\circ t_2=t_3\circ (0,\hueca{I}_Y)$ and, therefore, 
$v_{2,2}=t_3$.
Since $t_2$ is a morphism in ${\cal Z}(\hat{Z})$, we have 
$$\begin{matrix}
  0&=&\sum_{i_0,i_1\geq 0}\hat{b}^{ad}_{i_0+i_1+1}\left( \begin{pmatrix}
             \delta_{X'}&\gamma'\\ 0&\delta_{Y'}\\
             \end{pmatrix}^{\otimes i_1}\otimes \begin{pmatrix}
             t_1&v_{1,2}\\ 0&t_3\\
             \end{pmatrix}\otimes\begin{pmatrix}
             \delta_{X}&\gamma\\ 0&\delta_{Y}\\
             \end{pmatrix}^{\otimes i_0}\right)\\
             &=&
             \begin{pmatrix}
              \hat{b}_1^{tw}(t_1)&t_1\star \gamma+\gamma'\star t_3+\hat{b}_1^{tw}(v_{1,2})\\
              0&\hat{b}_1^{tw}(t_3)\\
             \end{pmatrix}.\hfill\\
  \end{matrix}$$
Therefore, $\pi(t_1\star \gamma)=-\pi(\gamma'\star t_3)$. Then, from (1) and (2), we obtain
$$\pi((\sigma_{X'}\circ \gamma')\star t_3)=-\pi(\sigma_{X'}\circ (\gamma'\star t_3))=\pi(\sigma_{X'}\circ (t_1\star \gamma))=\pi(t_1[1]\star (\sigma_X\circ \gamma)).$$
\end{proof}

\begin{proposition}\label{P: de morf entre confl esp a morf entre tringulos} Suppose that the following diagram 
$$\begin{matrix}
  \xi:&(X,\delta_X)&\rightmap{f}&(E,\delta_E)&\rightmap{g}&(Y,\delta_Y)\\
  &\shortlmapdown{t_1}&&\shortrmapdown{t_2}&&\shortrmapdown{t_3}\\
   \xi':&(X',\delta_{X'})&\rightmap{f'}&(E',\delta_{E'})&\rightmap{g'}&(Y',\delta_{Y'})\\
  \end{matrix}$$
commutes in ${\cal Z}(\hat{Z})$ and that its rows are special conflations. Then, we have the following commutative diagram in ${\cal H}(\hat{Z})$ 
$$\begin{matrix}
  \tau_\xi:&(X,\delta_X)&\rightmap{\pi(f)}&(E,\delta_E)&\rightmap{\pi(g)}&(Y,\delta_Y)&\rightmap{w}&(X,\delta_X)[1]\\
  &\shortlmapdown{\pi(t_1)}&&\shortlmapdown{\pi(t_2)}&&\shortrmapdown{\pi(t_3)}&&\shortrmapdown{\pi(t_1)[1]}\\
   \tau_{\xi'}:&(X',\delta_{X'})&\rightmap{\pi(f')}&(E',\delta_{E'})&\rightmap{\pi(g')}&(Y',\delta_{Y'})&\rightmap{w'}&(X',\delta_{X'})[1].\\
  \end{matrix}$$  
\end{proposition}

\begin{proof} If the rows of the first diagram are canonical conflations, then our statement follows from (\ref{L: formulas para morfismos de triangulos})(3). In the general case, we have a commutative diagram in ${\cal H}(\hat{Z})$ of the form 
$$\begin{matrix}
  (X,\delta_X)&\rightmap{\pi(f_1)}&(E_1,\delta_{E_1})&\rightmap{\pi(g_1)}&(Y,\delta_Y)&\rightmap{w_1}&(X,\delta_X)[1]\\
  \shortlmapdown{\pi(\hueca{I}_X)}&&\shortlmapdown{\pi(h_1)}&&\shortrmapdown{\pi(\hueca{I}_Y)}&&\shortrmapdown{\pi(\hueca{I}_{X[1]})}\\
 (X,\delta_X)&\rightmap{\pi(f)}&(E,\delta_E)&\rightmap{\pi(g)}&(Y,\delta_Y)&\rightmap{w}&(X,\delta_X)[1]\\
  \shortlmapdown{\pi(t_1)}&&\shortlmapdown{\pi(t_2)}&&\shortrmapdown{\pi(t_3)}&&\shortrmapdown{\pi(t_1)[1]}\\
   (X',\delta_{X'})&\rightmap{\pi(f')}&(E',\delta_{E'})&\rightmap{\pi(g')}&(Y',\delta_{Y'})&\rightmap{w'}&(X',\delta_{X'})[1]\\
    \shortlmapdown{\pi(\hueca{I}_{X'})}&&\shortlmapdown{\pi(h_2)}&&\shortrmapdown{\pi(\hueca{I}_{Y'})}&&\shortrmapdown{\pi(\hueca{I}_{X'[1]})}\\
    (X',\delta_{X'})&\rightmap{\pi(f'_1)}&(E'_1,\delta_{E'_1})&\rightmap{\pi(g_1')}&(Y',\delta_{Y'})&\rightmap{w'_1}&(X',\delta_{X'})[1]\\
  \end{matrix}$$
 where the first and the last rows are canonical triangles and $h_1$, $h_2$ are special isomorphisms. Therefore, we have the equality $\pi(t_1)[1]w=w'\pi(t_3)$. 
\end{proof}

\begin{lemma}\label{L: TR2 ''-->'' para H(hat(B))}
For any triangle in ${\cal H}(\hat{Z})$
$$\tau:\hbox{ \ }(X,\delta_X)\rightmap{u}(E,\delta_E)\rightmap{v}(Y,\delta_Y)\rightmap{w}(X,\delta_X)[1]$$
we have the triangle in ${\cal H}(\hat{Z}):$
$$\tau':\hbox{ \ }(E,\delta_E)\rightmap{v}(Y,\delta_Y)\rightmap{w}(X,\delta_X)[1]\rightmap{ \ -u[1] \ \ }(E,\delta_E)[1].$$
\end{lemma}

\begin{proof} We may assume that $\tau$ is a canonical triangle. Then, it has the form 
$$(X,\delta_X)\rightmap{\pi(f)}(E,\delta_E)\rightmap{\pi(g)}(Y,\delta_Y)\rightmap{\pi(\sigma_X\circ \gamma)}(X,\delta_X)[1],$$
for some canonical conflation 
$(X,\delta_X)\rightmap{f}(E,\delta_E)\rightmap{g}(Y,\delta_Y)$, 
where $E=X\oplus Y$\\ and $\delta_E=\begin{pmatrix}
           \delta_X&\gamma\\ 0&\delta_Y\\                                                          
           \end{pmatrix}$. 
 From (\ref{P: cada confl determina equiv de confl con J(X)}) and its proof, we have a commutative diagram in ${\cal H}(\hat{Z})$ of the form  
 $$\begin{matrix}
 (E,\delta_E)&\rightmap{\pi(h_\xi,g)^t}&J(X,\delta_X)\oplus (Y,\delta_Y)&\rightmap{\pi(\beta_X,\sigma_X\circ \gamma)}&(X,\delta_X)[1]\hfill\\
  \hbox{ \ \ \ \ \ }  \shortlmapdown{\hueca{I}_E}&&\shortrmapdown{\pi(s's)}&&\shortrmapdown{\hueca{I}_{X[1]}}\\
  (E,\delta_E)&\rightmap{ \ \pi(\alpha_1) \ }&(E_1,\delta_{E_1})&\rightmap{ \ \pi(\beta_1) \ }&(X,\delta_{X})[1],\hfill\\
    \end{matrix}$$
  where $s$ and $s'$ are isomorphisms, $(E,\delta_E)\rightmap{\alpha_1}(E_1,\delta_{E_1})\rightmap{\beta_1}(X,\delta_X)[1]$ is a canonical conflation with $E_1=E\oplus X[1]$ and
  $\delta_{E_1}=
  \begin{pmatrix}
         \delta_E&\gamma_1\\ 0& \delta_{X[1]}                                                                                                                                                                                                                                                                        \end{pmatrix}$, where $\gamma_1=(-\tau_X,0)^t$. Notice that 
         $$\sigma_E\circ \gamma_1=\begin{pmatrix}
                                   \sigma_X&0\\ 0&\sigma_Y\\
                                  \end{pmatrix}\circ 
                                  \begin{pmatrix}
                                   -\tau_X\\ 0\\
                                  \end{pmatrix}=
                                  \begin{pmatrix}
                                   -\sigma_X\circ \tau_X\\ 0\\
                                  \end{pmatrix}=
                                  \begin{pmatrix}
                                   -\hueca{I}_{X[1]}\\ 0
                                  \end{pmatrix}=-f[1].$$
So, we have the canonical triangle
$$(E,\delta_E)\rightmap{ \ \pi(\alpha_1) \ }(E_1,\delta_{E_1})\rightmap{ \ \pi(\beta_1) \ }(X,\delta_X)[1]\rightmap{ \ -\pi(f)[1] \ }(E,\delta_E)[1].$$
Therefore, since $J(X,\delta_X)$ is homologically trivial,  we have the triangle
$$(E,\delta_E)\rightmap{ \ \pi(g) \ }(Y,\delta_Y)\rightmap{ \ \pi(\sigma_X\circ \gamma) \ }(X,\delta_X)[1]\rightmap{ \ -\pi(f)[1] \ }(E,\delta_E)[1].$$
\end{proof}

\begin{lemma}\label{L: TR2 ''<--'' para H(hat(B))}
For any triangle in ${\cal H}(\hat{Z})$
$$\tau:\hbox{ \ }(X,\delta_X)\rightmap{u}(E,\delta_E)\rightmap{v}(Y,\delta_Y)\rightmap{w}(X,\delta_X)[1]$$
we have the triangle in ${\cal H}(\hat{Z}):$
$$\tau':\hbox{ \ }(Y,\delta_Y)[-1]\rightmap{ \ -w[-1] \  \ }(X,\delta_X)\rightmap{u}(E,\delta_E)\rightmap{v}(Y,\delta_Y).$$
\end{lemma}

\begin{proof} We may assume that $\tau$ is a canonical triangle. Then, it has the form 
$$(X,\delta_X)\rightmap{\pi(f)}(E,\delta_E)\rightmap{\pi(g)}(Y,\delta_Y)\rightmap{\pi(\sigma_X\circ \gamma)}(X,\delta_X)[1],$$
for some canonical conflation 
$(X,\delta_X)\rightmap{f}(E,\delta_E)\rightmap{g}(Y,\delta_Y)$, 
where $E=X\oplus Y$\\ and $\delta_E=\begin{pmatrix}
           \delta_X&\gamma\\ 0&\delta_Y\\                                                          
           \end{pmatrix}$. 
           From (\ref{P: cada confl determina equiv de confl con J'(Y)}) and its proof, we have a commutative diagram in ${\cal H}(\hat{Z})$ of the form
 $$\begin{matrix}
  (Y,\delta_Y)[-1]&\rightmap{\pi(\alpha_{Y[-1]},-h^\gamma)^t}&J(Y[-1],\delta_{Y[-1]})\oplus(X,\delta_X)&\rightmap{\pi(h^\xi,f)}&(E,\delta_E)\\
  \shortlmapdown{\hueca{I}_{Y[-1]}}&&\shortrmapdown{\pi(s's)}&&\shortrmapdown{\hueca{I}_E}\\
  (Y,\delta_Y)[-1]&\rightmap{ \ \ \pi(\alpha_1) \ \ }&(E_1,\delta_{E_1})&\rightmap{ \ \ \pi(\beta_1) \ \ }&(E,\delta_E),\hfill\\
    \end{matrix}$$  
  where $s$ and $s'$ are isomorphisms, $(Y,\delta_Y)[-1]\rightmap{\alpha_1}(E_1,\delta_{E_1})\rightmap{\beta_1}(E,\delta_E)$ is a canonical conflation with $E_1=Y[-1]\oplus E$ and 
  $\delta_{E_1}=
  \begin{pmatrix}
         \delta_{Y[-1]}&\gamma_1\\ 0& \delta_E                                                                                                                                                                                                                                                                        \end{pmatrix}$, where $\gamma_1=(0,-\tau_{Y[-1]})$. 
Notice that  $\sigma_{Y[-1]}\circ \gamma_1=(0,-\sigma_{Y[-1]}\circ\tau_{Y[-1]})=(0,-\hueca{I}_Y)=-g.$
So, we have the canonical triangle
$$(Y,\delta_Y)[-1]\rightmap{\pi(\alpha_1)}
(E_1,\delta_{E_1})\rightmap{\pi(\beta_1)}
(E,\delta_E)\rightmap{-\pi(g)}(Y,\delta_Y).$$
Therefore, since $J(Y[-1],\delta_{Y[-1]})$ is homologically trivial,  
we have the triangle
$$(Y,\delta_Y)[-1]\rightmap{ \ \pi((\sigma_X\circ \gamma)[-1]) \ }
(X,\delta_X)\rightmap{\pi(f)}(E,\delta_E)\rightmap{-\pi(g)}(Y,\delta_Y).$$
But we have the following commutative diagram in ${\cal H}(\hat{Z})$
$$\begin{matrix}(Y,\delta_Y)[-1]&\rightmap{ \ \pi((\sigma_X\circ \gamma)[-1]) \ }&
(X,\delta_X)&\rightmap{\pi(f)}&(E,\delta_E)&\rightmap{-\pi(g)}&(Y,\delta_Y)\\
\shortlmapdown{-\hueca{I}_{Y[-1]}}&&\shortlmapdown{\hueca{I}_X}&&\shortlmapdown{\hueca{I}_E}&&\shortlmapdown{-\hueca{I}_Y}\\
(Y,\delta_Y)[-1]&\rightmap{ \ -\pi(\sigma_X\circ \gamma)[-1] \ }&
(X,\delta_X)&\rightmap{\pi(f)}&(E,\delta_E)&\rightmap{\pi(g)}&(Y,\delta_Y),\\
\end{matrix}$$
so the lower row is a triangle.
\end{proof}

\begin{proposition}\label{T: H(hat(B)) es cat pretriangulada}
 The category ${\cal H}(\hat{Z})$ is a pretriangulated category with the class of triangles defined in (\ref{D: canonical triangles and triangles}). 
\end{proposition}

\begin{proof} The condition TR1(a) follows from the definition of triangle in ${\cal H}(\hat{Z})$. The condition TR1(b) is also satisfied because, for any  object $(X,\delta_X)$ in ${\cal Z}(\hat{Z})$, we have the canonical conflation in ${\cal Z}(\hat{Z})$ 
$$(X,\delta_X)\rightmap{\hueca{I}_X}(X,\delta_X)\rightmap{0}(0,0),$$
which gives rise to the triangle 
 $(X,\delta_X)\rightmap{ \ \pi(\hueca{I}_X) \ }(X,\delta_X)\rightmap{0}(0,0)\rightmap{0}(X,\delta_X)[1].$
 
Let us show that TR1(c) holds. Given any morphism 
$u:(X,\delta_X)\rightmap{}(Y,\delta_Y)$ in ${\cal H}(\hat{Z})$, we have $u=\pi(f)$, for some morphism $f:(X,\delta_X)\rightmap{}(Y,\delta_Y)$ in ${\cal Z}(\hat{Z})$. From (\ref{P: morfismos encajados en conflaciones con J}), we have a conflation  of the form 
$$\eta:{ \ }(X,\delta_X)\rightmap{ \ (f',f)^t \ }J(X,\delta_X)\oplus (Y,\delta_Y)\rightmap{ \ \ (g',g) \ \ }(W,\delta_W),$$
which is related to a canonical conflation 
$\eta_n:(X,\delta_X)\rightmap{f_n}(E_n,\delta_{E_n})\rightmap{g_n}(Y,\delta_Y)$ as in (\ref{R: conflations give rise to triangles}). 
If $h\in \Hom_{{\cal Z}(\hat{Z})}((W,\delta_W),(X,\delta_X)[1])$ is the morphism such that $\Psi(h)=[\eta_n]$, we have the commutative diagram 
$$\begin{matrix}
(X,\delta_X)&\rightmap{ \ \pi(f',f)^t \ }&J(X,\delta_X)\oplus (Y,\delta_Y)&\rightmap{ \ \ \pi(g',g) \ \ }&(W,\delta_W)&\rightmap{\pi(h)}&(X,\delta_X)[1]\\
\shortlmapdown{\hueca{I}_E}&&\shortlmapdown{\cong}&&\shortlmapdown{\hueca{I}_{W}}&&\shortlmapdown{\hueca{I}_{X[1]}}\\
(X,\delta_X)&\rightmap{ \ \pi(f_n) \ }&(E_n,\delta_{E_n})&\rightmap{ \ \ \pi(g_n) \ \ }&(W,\delta_W)&\rightmap{\pi(h)}&(X,\delta_X)[1],\\
       \end{matrix}$$
where the lower row is the canonical triangle associated to $\eta_n$. Since $J(X,\delta_X)$ is homologically trivial, we have the following  
triangle in ${\cal H}(\hat{Z})$: 
$$(X,\delta_X)\rightmap{ \ \pi(f) \ }(Y,\delta_Y)\rightmap{ \ \ \pi(g) \ \ }(W,\delta_W)\rightmap{\pi(h)}(X,\delta_X)[1].$$

The condition TR2 follows from (\ref{L: TR2 ''-->'' para H(hat(B))}) and (\ref{L: TR2 ''<--'' para H(hat(B))}). 

Now, we proceed to prove TR3. Given a commutative diagram in ${\cal H}(\hat{Z})$  
 $$\begin{matrix}
 \tau:\hbox{ \ } (X,\delta_X)&\rightmap{u}&(E,\delta_E)&\rightmap{v}&(Y,\delta_Y)&\rightmap{w}&(X,\delta_X)[1]\\
 \hbox{ \  \ \ \ }\shortrmapdown{\theta_1}&&\shortrmapdown{\theta_2}&&&&\\
 \tau':\hbox{ \ }(X',\delta_{X'})&\rightmap{u'}&(E',\delta_{E'})&\rightmap{v'}&(Y',\delta_{Y'})&\rightmap{w'}&(X',\delta_{X'})[1],\\
   \end{matrix}$$
   with rows which are  triangles, we want to find a morphism $\theta_3:(Y,\delta_Y)\rightmap{}(Y',\delta_{Y'})$ such that   $(\theta_1,\theta_2,\theta_3):\tau\rightmap{}\tau'$ is a morphism of triangles. We may assume that the triangles $\tau$ and $\tau'$ are canonical triangles. Then,  we have a canonical conflations in ${\cal Z}(\hat{Z})$ of the form: 
   $$\xi:(X,\delta_X)\rightmap{f}(E,\delta_E)\rightmap{g}(Y,\delta_Y)
   \hbox{ and } 
  \xi':(X',\delta_{X'})\rightmap{f'}(E',\delta_{E'})\rightmap{g'}(Y',\delta_{Y'})$$
 and morphisms $t_1:(X,\delta_X)\rightmap{}(X',\delta_{X'})$ and $t_2:(E,\delta_E)\rightmap{}(E',\delta_{E'})$ in ${\cal Z}(\hat{Z})$, such that 
 $\pi(f)=u$, $\pi(g)=v$, $\pi(f')=u'$, $\pi(g')=v'$, $\pi(t_1)=\theta_1$, and $\pi(t_2)=\theta_2$.

 Since $\pi(f'\star t_1)=u'\theta_1=\theta_2 u=\pi(t_2\star f)$, there is a homologically trivial   morphism $s:(X,\delta_X)\rightmap{}(E',\delta_{E'})$  in $\tw(\hat{Z})$ such that 
 $f'\star t_1=t_2\star f+s$. From (\ref{L: morf homol triv se fact por conflaciones esp}), we know that $s=s'\star f$, for some morphism $s':(E,\delta_E)\rightmap{}(E',\delta_{E'})$ in ${\cal Z}(\hat{Z})$. Then, if we define $t'_2:=t_2+s'$, we get   $f'\star t_1=t_2\star f+s'\star f=t'_2\star f$.
 
 Moreover,  $(g'\star t'_2)\star f=g'\star(t'_2\star f)=g'\star (f'\star t_1)=(g'\star f')\star t_1=0$. Thus, using (\ref{L: conlfacion es par exacto}), we know that  $g$ is the cokernel of $f$, and have the existence of a morphism $t_3:(Y,\delta_Y)\rightmap{}(Y',\delta_{Y'})$ such that $g'\star t'_2=t_3\star g$. So we get the following commutative diagram in ${\cal Z}(\hat{Z})$:
 $$\begin{matrix}
   \xi:\hbox{ \ } (X,\delta_X)&\rightmap{f}&(E,\delta_E)&\rightmap{g}&(Y,\delta_Y)\\
   \shortlmapdown{t_1}&&\shortlmapdown{t'_2}&&\shortlmapdown{t_3}\\
    \xi':\hbox{ \ } (X',\delta_{X'})&\rightmap{f'}&(E',\delta_{E'})&\rightmap{g'}&(Y',\delta_{Y'}).\\
   \end{matrix}$$
Then, we can apply (\ref{P: de morf entre confl esp a morf entre tringulos}) to this diagram and take $\theta_3=\pi(t_3)$, to obtain the wanted commutative diagram  in ${\cal H}(\hat{Z})$ 
$$\begin{matrix}
  (X,\delta_X)&\rightmap{u}&(E,\delta_E)&\rightmap{v}&(Y,\delta_Y)&\rightmap{w}&(X,\delta_X)[1]\\
  \shortlmapdown{\theta_1}&&\shortlmapdown{\theta_2}&&\shortlmapdown{\theta_3}&&\shortlmapdown{\theta_1[1]}\\
   (X',\delta_{X'})&\rightmap{u'}&(E',\delta_{E'})&\rightmap{v'}&(Y',\delta_{Y'})&\rightmap{w'}&(X',\delta_{X'})[1].\\
  \end{matrix}$$ 
\end{proof}

\begin{remark}\label{R: special iso of direct complements}
Given a right $\hat{S}$-module $E$ with direct sum decompositions $X\oplus W=E=X\oplus W'$, consider the canonical projections  $p_X$, $p_W$  associated to the first direct sum decomposition, and the canonical projections $p'_X$, $p_{W'}$ associated to the second decomposition. Let $s_W$, $s_X$, and $s_{W'}$ be the corresponding canonical injections. Then, $p_{W'}s_W:W\rightmap{}W'$ is an isomorphism.

Indeed,  
$(p_Ws_{W'})(p_{W'}s_W)=p_W(1_E-s_Xp'_X)s_W=p_Ws_W-p_Ws_Xp'_Xs_W=p_Ws_W=id_W.$
So, in this case, we have the corresponding special isomorphism $$\phi:=L(p_{W'}s_W):W\rightmap{}W' \hbox{  \ in  \ } \ad(\hat{Z}).$$  
\end{remark}

\begin{theorem}\label{T: H(hat(B)) es cat triangulada}
 The category ${\cal H}(\hat{Z})$ is a triangulated category with the class of triangles defined in (\ref{D: canonical triangles and triangles}). 
\end{theorem}

\begin{proof} It only remains to prove the octahedral axiom. We split this proof in two parts. 

\medskip
\noindent{\bf Part 1: The canonical case.}
\medskip

We prove the octahedral axiom for canonical triangles:
 $$\begin{matrix}\tau_\xi:&(X,\delta_X)&\rightmap{u}&(Y,\delta_Y)&\rightmap{i}&(U',\delta_{U'})&\rightmap{\hat{i}}&(X,\delta_X)[1]\\
       \tau_\eta:&(Y,\delta_Y)&\rightmap{v}&(U,\delta_U)&\rightmap{j}&(X',\delta_{X'})&\rightmap{\hat{j}}&(Y,\delta_Y)[1]\\
       \tau_\zeta:&(X,\delta_X)& \rightmap{vu}&(U,\delta_U)&\rightmap{w}&(Y',\delta_{Y'})&\rightmap{\hat{w}}&(X,\delta_X)[1],\\
      \end{matrix}$$
 which are associated respectively to canonical conflations:
 $$\begin{matrix}\xi:&(X,\delta_X)&\rightmap{u_1}&(Y,\delta_Y)&\rightmap{i_1}&(U',\delta_{U'})\\
\eta:&(Y,\delta_Y)&\rightmap{v_1}&(U,\delta_U)&\rightmap{j_1}&(X',\delta_{X'})\\
\zeta:&(X,\delta_X)&\rightmap{v_1u_1}&(U,\delta_U)&\rightmap{w_1}&(Y',\delta_{Y'}).\\
 \end{matrix}$$
Then, we have right $\hat{S}$-module decompositions. 
$$Y=X\oplus U' \hbox{ \  and \ }X\oplus Y'=U=X\oplus U'\oplus X'.$$ 
Moreover, we have $\delta_Y=\begin{pmatrix}
                                      \delta_X&\gamma\\ 0&\delta_{U'}\\
                                     \end{pmatrix}:X\oplus U'\rightmap{}X\oplus U'$, 
while $\delta_U$ has the following matrix form, associated to the  decomposition
$U=X\oplus U'\oplus X'$:
$$\delta_U=\begin{pmatrix}
           \delta_X&\gamma&\beta_1\\
           0&\delta_{U'}&\beta_2\\
           0&0&\delta_{X'}
          \end{pmatrix}.$$   
          
Now, we have that the special morphism $w_1$, which appears in the canonical conflation $\zeta$, has the form $w_1=(0,\phi):X\oplus(U'\oplus X')\rightmap{}Y'$, where $\phi:U'\oplus X'\rightmap{}Y'$ is the special isomorphism considered in (\ref{R: special iso of direct complements}),  
and the special morphism $w'_1:=(0,\phi^{-1})^t:Y'\rightmap{}X\oplus (U'\oplus X')$, which satisfy $w_1\circ w_1'=\hueca{I}_{Y'}$ and 
$w'_1\circ w_1=\begin{pmatrix}
          0&0\\ 0&\hueca{I}_{U'\oplus X'}\\
         \end{pmatrix}$. Consider also  the special morphism  $j'_1=(0,\hueca{I}_{U'},0)^t:U'\rightmap{}X\oplus U'\oplus X'$. Finally, we consider the special morphisms $f_1:=w_1\circ j'_1:U'\rightmap{}Y'$ and $g_1:=j_1\circ w'_1:Y'\rightmap{}X'$.
Then, we have the following commutative diagram in $\ad(\hat{Z})$:
$$\begin{matrix}
  X\oplus U'&\rightmap{v_1}&X\oplus U'\oplus X'&\rightmap{j_1}&X'\\
  \shortlmapdown{i_1}&&\shortlmapdown{w_1}&&\shortlmapdown{\hueca{I}_{X'}}\\
  U'&\rightmap{f_1}&Y'&\rightmap{g_1}&X'.\\
  \end{matrix}$$
We claim that the lower row determines a special conflation
$$\begin{matrix}(U',\delta_{U'})&\rightmap{f_1}&(Y',\delta_{Y'})&\rightmap{g_1}&(X',\delta_{X'}).\end{matrix}$$
Let us show first that $f_1$ and $g_1$ are morphisms in ${\cal Z}(\hat{Z})$. For this, notice that  $w_1:(U,\delta_U)\rightmap{}(Y',\delta_{Y'})$ and $j_1:(U,\delta_U)\rightmap{}(X',\delta_{X'})$ are morphisms in ${\cal Z}(\hat{Z})$, because they appear in $\zeta$ and $\eta$, respectively. The morphism $\phi$ has matrix form  $\phi=(\phi_1,\phi_2):U'\oplus X'\rightmap{}Y'$. Then, we have 
$$\begin{matrix}
  0&=&\delta_{Y'}\circ w_1+w_1\circ \delta_U
  = \delta_{Y'}\circ(0,\phi_1,\phi_2)+(0,\phi_1,\phi_2)\circ
  \begin{pmatrix}
         \delta_X&\gamma&\beta_1\\
         0&\delta_{U'}&\beta_2\\
         0&0&\delta_{X'}
  \end{pmatrix}\hfill\\
  &=&(0,\delta_{Y'}\circ \phi_1+\phi_1\circ \delta_{U'},\delta_{Y'}\circ \phi_2+\phi_1\circ \beta_2+\phi_2\circ \delta_{X'}),\hfill\\
  \end{matrix}$$
which implies that $\phi_1:(U',\delta_{U'})\rightmap{}(Y',\delta_{Y'})$ is a special morphism in ${\cal Z}(\hat{Z})$. Since   $f_1=w_1\circ j'_1=\phi_1$, we have that $f_1:(U',\delta_{U'})\rightmap{}(Y',\delta_{Y'})$  is a special morphism in ${\cal Z}(\hat{Z})$.

Since $w_1\circ \delta_U+\delta_{Y'}\circ w_1=0$ and $w_1\circ w'_1=\hueca{I}_{Y'}$, we get   $\delta_{Y'}=-w_1\circ \delta_U\circ w_1'$. Then, using that $j_1:(U,\delta_U)\rightmap{}(X',\delta_{X'})$ belongs to ${\cal Z}(\hat{Z})$, we obtain
$$\begin{matrix}
\delta_{X'}\circ g_1+g_1\circ\delta_{Y'}
&=&
\delta_{X'}\circ (j_1\circ w_1')+(j_1\circ w'_1)\circ \delta_{Y'}\hfill\\
&=&
-(\delta_{X'}\circ j_1)\circ w_1'+(j_1\circ w'_1)\circ \delta_{Y'}\hfill\\
&=&

(j_1\circ \delta_U)\circ w'_1-(j_1\circ w'_1)\circ (w_1\circ \delta_U\circ w'_1)\hfill\\
&=&
j_1\circ [\delta_U-(w'_1\circ w_1)\circ \delta_U]\circ w'_1=0,\hfill\\
\end{matrix}$$
and $g_1:(Y',\delta_{Y'})\rightmap{}(X',\delta_{X'})$ is a morphism in   ${\cal Z}(\hat{Z})$.  

Now, in order to show that the sequence 
$(U',\delta_{U'})\rightmap{f_1}(Y',\delta_{Y'})\rightmap{g_1}(X',\delta_{X'})$ is a special conflation, 
since all the morphisms we have considered in this proof are special, we abuse the language and consider them as morphisms of right $\hat{S}$-modules. So, we have to show that the sequence  $0\rightmap{}U'\rightmap{f_1}Y'\rightmap{g_1}X'\rightmap{}0$ is exact.  

Consider the morphisms of right $\hat{S}$-modules $f'_1:=p_1\circ w'_1:Y'\rightmap{}U'$ and $g'_1:=w_1\circ i'_1:X'\rightmap{}Y'$,
where 
$i'_1:=(0,0,\hueca{I}_{X'})^t:X'\rightmap{}X\oplus U'\oplus X'$ 
 and 
$p_1:=(0,\hueca{I}_{U'},0):X\oplus U'\oplus X'\rightmap{}U'.$
If we set $\phi^{-1}=(\phi'_1,\phi'_2)^t:Y'\rightmap{}U'\oplus X'$, we get $1_{Y'}=\phi\phi^{-1}=\phi_1\phi'_1+\phi_2\phi'_2$ and $w'_1=(0,\phi'_1,\phi'_2)^t:Y'\rightmap{}X\oplus U'\oplus X'$. Then, by direct computations, we obtain the equalities:
$$g_1\circ f_1=0,\hbox{\ } f'_1\circ g'_1=0,\hbox{\ } f'_1\circ f_1=\hueca{I}_{U'},\hbox{\ } g_1\circ g'_1=\hueca{I}_{X'},\hbox{\ } f_1\circ f'_1+g'_1\circ g_1=\hueca{I}_{Y'}.$$
They imply that $Y'=f_1(U')\oplus g'_1(X')$ and $\Ker g_1=f_1(U')$. So, we have the wanted split sequence.  Then, we have the commutative diagram in ${\cal Z}(\hat{Z})$ 
$$\begin{matrix}
  (Y,\delta_Y)&\rightmap{v_1}&(U,\delta_U)&\rightmap{j_1}&(X',\delta_{X'})\\
  \shortlmapdown{i_1}&&\shortlmapdown{w_1}&&\shortlmapdown{\hueca{I}_{X'}}\\
  (U',\delta_{U'})&\rightmap{f_1}&(Y',\delta_{Y'})&\rightmap{g_1}&(X',\delta_{X'}),\\
  \end{matrix}$$
with special conflations as rows. If we take $f:=\pi(f_1)$ and $g:=\pi(g_1)$, from (\ref{P: de morf entre confl esp a morf entre tringulos}), we get the following commutative diagram in ${\cal H}(\hat{Z})$ 
$$\begin{matrix}
  (Y,\delta_Y)&\rightmap{v}&(U,\delta_U)&\rightmap{j}&(X',\delta_{X'})&\rightmap{\hat{j}}&(Y,\delta_Y)[1]\\
  \shortlmapdown{i}&&\shortlmapdown{w}&&\shortlmapdown{\hueca{I}_{X'}}&&
  \shortlmapdown{i[1]}\\
   (U',\delta_{U'})&\rightmap{f}&(Y',\delta_{Y'})&\rightmap{g}&(X',\delta_{X'})&\rightmap{\hat{g}}&(U',\delta_{U'})[1]\\
  \end{matrix}$$ with triangles as rows.  Now, observe that we have the following commutative diagram in ${\cal Z}(\hat{Z})$ 
  $$\begin{matrix}
      (X,\delta_X)&\rightmap{ \ u_1v_1 \ }&(U,\delta_U)&\rightmap{w_1}&(Y',\delta_{Y'})\\
  \shortlmapdown{u_1}&&\shortlmapdown{\hueca{I}_U}&&\shortlmapdown{g_1}\\
  (Y,\delta_Y)&\rightmap{ \ v_1 \ }&(U,\delta_U)&\rightmap{j_1}&(X',\delta_{X'}),\\
    \end{matrix}$$
where the rows are canonical conflations by assumption.  
Using again (\ref{P: de morf entre confl esp a morf entre tringulos}), we get the following commutative diagram in ${\cal H}(\hat{Z})$ 
$$\begin{matrix}
  (X,\delta_X)&\rightmap{uv}&(U,\delta_U)&\rightmap{w}&(Y',\delta_{Y'})&\rightmap{\hat{w}}&(X,\delta_X)[1]\\
  \shortlmapdown{u}&&\shortlmapdown{\hueca{I}_U}&&\shortlmapdown{g}&&
  \shortlmapdown{u[1]}\\
   (Y,\delta_Y)&\rightmap{v}&(U,\delta_U)&\rightmap{j}&(X',\delta_{X'})&\rightmap{\hat{j}}&(Y,\delta_Y)[1]\\
  \end{matrix}$$ with triangles as rows. Therefore, we get $u[1]\hat{w}=\hat{j}g$, and after a shifting we obtain 
  $u\hat{w}[-1]=\hat{j}[-1]g[-1]$. 
  From the  commutative diagram in ${\cal Z}(\hat{Z})$
  $$\begin{matrix}
      (X,\delta_X)&\rightmap{ \ u_1 \ }&(Y,\delta_Y)&\rightmap{i_1}&(U',\delta_{U'})\\
  \shortlmapdown{\hueca{I}_X}&&\shortlmapdown{v_1}&&\shortlmapdown{f_1}\\
  (X,\delta_X)&\rightmap{ \ v_1u_1 \ }&(U,\delta_U)&\rightmap{w_1}&(Y',\delta_{Y'}),\\
    \end{matrix}$$
where the rows are canonical conflations by assumption,
and (\ref{P: de morf entre confl esp a morf entre tringulos}), we get the following commutative diagram in ${\cal H}(\hat{Z})$
$$\begin{matrix}
  (X,\delta_X)&\rightmap{u}&(Y,\delta_Y)&\rightmap{i}&(U',\delta_{U'})&\rightmap{\hat{i}}&(X,\delta_X)[1]\\
  \shortlmapdown{\hueca{I}_X}&&\shortlmapdown{v}&&\shortlmapdown{f}&&
  \shortlmapdown{\hueca{I}_{X[1]}}\\
   (X,\delta_X)&\rightmap{vu}&(U,\delta_U)&\rightmap{w}&(Y',\delta_{Y'})&\rightmap{\hat{w}}&(X,\delta_{X})[1]\\
  \end{matrix}$$ with triangles as rows. In particular, we have $\hat{w}f=\hat{i}$, as wanted. 

  %%%%%%%%%%%%%%%%%%%%%%%%%%%%%%%%%%%%%%%%%%%
  \medskip
  \noindent{\bf Part 2: The general case.}
  \medskip

 Assume that 
  we have triangles 
 $$\begin{matrix}\tau_1:&X&\rightmap{u}&Y
 &\rightmap{i}&Z'&\rightmap{i'}&TX\\
       \tau_2:&Y&\rightmap{v}&Z
       &\rightmap{j}&X'&\rightmap{j'}&TY\\
       \tau_3:&X& \rightmap{v u}&Z
       &\rightmap{k}&Y'&\rightmap{k'}&TX
       .\\
      \end{matrix}$$
For the sake of notational simplicity, in this part of the proof, the objects of 
${\cal Z}(\hat{Z})$ are written without making explicit their differential. 
We will choose appropriately  some canonical triangles isomorphic to the preceding ones, then apply the octahedral axiom to them, and from  there we show the octahedral axiom for $\tau_1,\tau_2,\tau_3$. 
We start with any isomorphism of triangles $(\theta_1,\theta_2,\theta_3)$ from $\tau_1$ to a canonical triangle $\tau_{\xi_1}$, that gives us a commutative diagram 
 $$\hbox{\hskip1cm}\begin{matrix}
 \tau_1:&X&\rightmap{u}&Y
 &\rightmap{i}&Z'&\rightmap{i'}&TX\\
 &\shortlmapdown{\theta_1}&&\shortlmapdown{\theta_2}&&\shortrmapdown{\theta_3}&&\shortrmapdown{T(\theta_1)}\\
 \tau_{\xi_1}:&A&\rightmap{\pi(a)}&B
 &\rightmap{\pi(a')}&C'&\rightmap{\pi(a'')}&TA,\\
\end{matrix}\hbox{\hskip1cm}:(D_1)$$
where $\xi_1:A\rightmap{a}B\rightmap{a'}C'$ is a canonical conflation.  
Consider the morphism $v\theta_2^{-1}:B\rightmap{}Z$ and a morphism $h:B\rightmap{}Z$ in ${\cal Z}(\hat{Z})$ with $\pi(h)=v\theta_2^{-1}$. 
From (\ref{P: morfismos encajados en conflaciones con J}), using $h$, we obtain a conflation of the form 
$$\eta_1: { \ }B\rightmap{ \ (\alpha_B,h)^t \ }J(B)\oplus Z\rightmap{ \ d \ }A'.$$
 Then, if we denote by $\sigma:Z\rightmap{}J(B)\oplus Z$ the canonical injection in ${\cal Z}(\hat{Z})$, by (\ref{R: conflations give rise to triangles}), we have a 
commutative diagram 
 $$\begin{matrix}
    \tau_2:&Y&\rightmap{v}&Z
       &\rightmap{j}&X'&\rightmap{j'}&TY
       \\
      &\shortlmapdown{\theta_2}&&\shortlmapdown{\pi(\sigma)}&&&&\\  
     \tau_{\eta_1}& B&\rightmap{ \ \pi(\alpha_B,h)^t \ }&J(B)\oplus Z&\rightmap{\pi(d)}&A'&\rightmap{}&TB\\
      &\parallel&&\shortlmapdown{\zeta_2}&&\parallel&&\parallel\\  
      \tau_{\xi_2}:&B&\rightmap{\pi(b)}&C&\rightmap{\pi(b')}&A'&\rightmap{\pi(b'')}&TB\\
      \end{matrix}$$
where $\xi_2:B\rightmap{b}C\rightmap{b'}A'$ is a canonical conflation and $\zeta_2$ is an isomorphism in ${\cal H}(\hat{Z})$. Since $\pi(\sigma)=\hueca{I}_Z$, the following diagram commutes:
$$\begin{matrix}
    \tau_2:&Y&\rightmap{v}&Z
       &\rightmap{j}&X'&\rightmap{j'}&TY\\
      &\shortlmapdown{\theta_2}&&\shortlmapdown{\zeta_2}&&&&\\  
 \tau_{\xi_2}:&B&\rightmap{\pi(b)}&C&\rightmap{\pi(b')}&A'&\rightmap{\pi(b'')}&TB\\
      \end{matrix}$$
      which by TR3, can be completed to a commutative diagram
    $$\begin{matrix}
    \tau_2:&Y&\rightmap{v}&Z
       &\rightmap{j}&X'&\rightmap{j'}&TY\\
      &\shortlmapdown{\theta_2}&&\shortlmapdown{\zeta_2}&&\shortlmapdown{\beta_3}&&\shortrmapdown{T(\theta_2)}\\  
 \tau_{\xi_2}:&B&\rightmap{\pi(b)}&C&\rightmap{\pi(b')}&A'&\rightmap{\pi(b'')}&TB.\\
      \end{matrix}\hbox{\hskip1cm}:(D_2)$$  
  Since ${\cal H}(\hat{Z})$ is pretriangulated and $\theta_2$ and $\zeta_2$ are isomorphisms, so is $\beta_3$.
  
  From (\ref{R: caract de infl (defl) canonic, compos de infl (defl) canonic es infl (defl) canonic}), there is a canonical conflation of the form 
    $$\xi_3:\hbox{ \ }A\rightmap{b\star a}C\rightmap{p}B'.$$
    Then, we have a commutative diagram 
 $$\begin{matrix}
    \tau_3:&X&\rightmap{vu}&Z
       &\rightmap{k}&Y'&\rightmap{k'}&TX\\
      &\shortlmapdown{\theta_1}&&\shortlmapdown{\zeta_2}&&&&\\  
     \tau_{\xi_3}& A&\rightmap{\pi(b\star a)}&C&\rightmap{\pi(p)}&B'&\rightmap{\pi(p')}&TA,\\
      \end{matrix}$$
      which by TR3, can be completed to a commutative diagram
    $$\begin{matrix}
    \tau_3:&X&\rightmap{vu}&Z
       &\rightmap{k}&Y'&\rightmap{k'}&TX\\
      &\shortlmapdown{\theta_1}&&\shortlmapdown{\zeta_2}&&\shortlmapdown{\zeta_3}&&\shortrmapdown{T(\theta_1)}\\  
 \tau_{\xi_3}:&A&\rightmap{\pi(b\star a)}&C&\rightmap{\pi(p)}&B'&\rightmap{\pi(p')}&TA.\\
      \end{matrix}\hbox{\hskip1cm}:(D_3)$$  
  Since $\theta_1$ and $\zeta_2$ are isomorphisms, so is $\zeta_3$.

Apply the octahedral axiom to the canonical triangles 
$$\begin{matrix}
\tau_{\xi_1}:&A&\rightmap{\pi(a)}&B
 &\rightmap{\pi(a')}&C'&\rightmap{\pi(a'')}&TA\\
 \tau_{\xi_2}:&B&\rightmap{\pi(b)}&C&\rightmap{\pi(b')}&A'&\rightmap{\pi(b'')}&TB      \\
   \tau_{\xi_3}:&A&\rightmap{\pi(b\star a)}&C&\rightmap{\pi(p)}&B'&\rightmap{\pi(p')}&TA \\
      \end{matrix}$$
  to obtain the triangle 
$\hbox{ \ }C'\rightmap{ \ f \ }B'\rightmap{ \ g \ }A'\rightmap{ \ T(\pi(a'))\pi(b'') \ }TC'$
and the commutative diagram 
 $$\begin{matrix}
    T^{-1}B'&\rightmap{ \ T^{-1}\pi(p') \ }&A&\rightmap{1_A}&A&&&&\\
    \lmapdown{ \ T^{-1}(g) \ }&\scriptstyle{(4)}&\lmapdown{\pi(a)}&&\rmapdown{\pi(b)\pi(a)}&&&&\\
    T^{-1}A'&\rightmap{T^{-1}(\pi(b''))}&B&\rightmap{\pi(b)}&C&\rightmap{\pi(b')}&A'&\rightmap{\pi(b'')}&TB\\
    &&\lmapdown{\pi(a')}&\scriptstyle{(1)}&\rmapdown{\pi(p)}&\scriptstyle{(2)}&\rmapdown{1_{A'}}&&\rmapdown{T(\pi(a'))}\\
    &&C'&\rightmap{f}&B'&\rightmap{g}&A'&\rightmap{T(\pi(a'))\pi(b'')}&TC'\\
    &&\lmapdown{\pi(a'')}&\scriptstyle{(3)}&\rmapdown{\pi(p')}&&&&\\
    &&TA&\rightmap{1_{TA}}&TA.&&&&\\
   \end{matrix}$$
Define $\bar{f}:=\zeta_3^{-1}f\theta_3$ and $\bar{g}:=\beta_3^{-1}g\zeta_3$, 
then we have the diagram 
$$\begin{matrix}
\overline{\tau}:&\hbox{ \ }Z'&\rightmap{\bar{f}}&Y'&\rightmap{\bar{g}}&X'&\rightmap{ \ \ T(i)j' \ \ }&TZ'\\
 &\shortlmapdown{\theta_3}&&\shortlmapdown{\zeta_3}&&\shortlmapdown{\beta_3}&&\shortrmapdown{T(\theta_3)}\\
\tau:&\hbox{ \ }C'&\rightmap{f}&B'&\rightmap{g}&A'&\rightmap{T(\pi(a'))\pi(b'')}&TC'.\\
\end{matrix}$$
The first two squares commute by definition of $\bar{f}$ and $\bar{g}$. The third one commutes because, from the commutativity of $(D_2)$ and $(D_1)$, we have
$$T(\pi(a'))\pi(b'')\beta_3=T(\pi(a'))T(\theta_2)j'=T(\pi(a')\theta_2)j' =T(\theta_3 i)j'=T(\theta_3)T(i)j'.$$
It follows that $\overline{\tau}$ is indeed a triangle in ${\cal H}(\hat{Z})$.   Now, we show the commutativity of the diagram 
$$\begin{matrix}
    T^{-1}Y'&\rightmap{ \ T^{-1}(k') \ }&X&\rightmap{1_X}&X&&&&\\
    \shortlmapdown{ \ T^{-1}(\bar{g}) \ }&&\shortlmapdown{u}&&\shortrmapdown{vu}&&&&\\
    T^{-1}X'&\rightmap{T^{-1}(j')}&Y&\rightmap{v}&Z&\rightmap{j}&X'&\rightmap{j'}&TY\\
    &&\shortlmapdown{i}&&\shortrmapdown{k}&&\shortrmapdown{1_{X'}}&&\shortrmapdown{T(i)}\\
    &&Z'&\rightmap{\bar{f}}&Y'&\rightmap{\bar{g}}&X'&\rightmap{T(i)j'}&TZ'\\
    &&\shortlmapdown{i'}&&\shortrmapdown{k'}&&&&\\
    &&TX&\rightmap{1_{TX}}&TX.&&&&\\
   \end{matrix}$$
Use successively the commutativity of $(D_1)$, $(1)$, and $(D_3)$,  $(D_2)$, to obtain
$$\bar{f}i=\zeta_3^{-1}
f\theta_3i=\zeta^{-1}_3f\pi(a')\theta_2
=
\zeta_3^{-1}\pi(p)\pi(b)\theta_2=
(\zeta_3^{-1}\pi(p)\zeta_2)(\zeta^{-1}_2\pi(b)\theta_2)=kv.$$
Use successively the commutativity of $(D_3)$, $(2)$, and $(D_2)$ to obtain
$$\bar{g}k=\beta_3^{-1}g\zeta_3k=\beta_3^{-1}g\pi(p)\zeta_2=\beta_3^{-1}\pi(b')\zeta_2=j.
$$
Use successively the commutativity of $(D_3)$, $(3)$, and $(D_1)$ to obtain
$$k'\bar{f}=(T(\theta_1)^{-1} \pi(p') \zeta_3)(\zeta_3^{-1}
f\theta_3)=T(\theta_1)^{-1}\pi(p') f\theta_3=T(\theta_1)^{-1} \pi(a'')\theta_3=i'.
$$
Finally, use successively the commutativity of $(D_1)$, $(D_3)$, $(4)$, and $(D_2)$ to obtain
$$\begin{matrix}
  T(u)k'
  &=& T(\theta_2^{-1}\pi(a)\theta_1)T(\theta_1)^{-1}\pi(p')\zeta_3\hfill\\
  &=&
  T(\theta_2)^{-1}T(\pi(a))\pi(p')\zeta_3\hfill\\
  &=&
  T(\theta_2)^{-1}\pi(b'')g\zeta_3=j'\beta_3^{-1}g\zeta_3=j'\bar{g}\hfill.\\
  \end{matrix}$$
  \end{proof}

%%%%%%%%%%%%%%%%%%%%%%%%%%%%%%%%%%%%%%%%%%%%%%%%%%%%%%%%%%%%%%%
%%%%%%%%%%%%%%%%%%%%%%%%%%%%%%%%%%%%%%%%%%%%%%%%%%%%%%%%%%%%%%%
%%%%%%%%%%%%%%%%%%%%%%%%%%%%%%%%%%%%%%%%%%%%%%%%%%%%%%%%%%%%%%%
%%%%%%%%%%%%%%%%%%%%%%%%%%%%%%%%%%%%%%%%%%%%%%%%%%%%%%%%%%%%%%%
%%%%%%%%%%%%%%%%%%%%%%%%%%%%%%%%%%%%%%%%%%%%%%%%%%%%%%%%%%%%%%%
%%%%%%%%%%%%%%%%%%%%%%%%%%%%%%%%%%%%%%%%%%%%%%%%%%%%%%%%%%%%%%%
%%%%%%%%%%%%%%%%%%%%%%%%%%%%%%%%%%%%%%%%%%%%%%%%%%%%%%%%%%%%%%%

\hskip2cm
%%%%%%%%%%%%%%%%%%%%%%%%%%%%%%%%%%%%%%%%%%%%%%%%%%%%%%%

\vbox{\noindent R. Bautista\\
Centro de Ciencias Matem\'aticas\\
Universidad Nacional Aut\'onoma de M\'exico\\
Morelia, M\'exico\\
raymundo@matmor.unam.mx\\}

\vbox{\noindent E. P\'erez\\
Facultad de Matem\'aticas\\
Universidad Aut\'onoma de Yucat\'an\\
M\'erida, M\'exico\\
jperezt@correo.uady.mx\\}

\vbox{\noindent L. Salmer\'on\\
Centro de Ciencias  Matem\'aticas\\
Universidad Nacional Aut\'onoma de M\'exico\\
Morelia, M\'exico\\
salmeron@matmor.unam.mx\\}
\end{document}